%% file: main.tex
\documentclass[3p]{elsarticle}

\usepackage{lineno,hyperref}
\modulolinenumbers[5]

\journal{Journal of Computational Mathematics and Data Science}

\usepackage{amsmath}
\usepackage{amsfonts}
\usepackage{graphicx}
\usepackage{amsthm}
\usepackage{moreverb}
\usepackage{amsfonts}
\usepackage{mathrsfs}
\usepackage{float}
\usepackage{fixltx2e}
\usepackage{gensymb}
\usepackage{xfrac}
\usepackage{lmodern}
\usepackage{color}
\usepackage{caption}
\usepackage{subcaption}
\usepackage{epstopdf}
\usepackage{listings}
\usepackage{todonotes}
\usepackage{booktabs}
\usepackage{array}
\usepackage{tabularx}
\usepackage{placeins}
\usepackage{xcolor}

\newtheorem{theorem}{Theorem}
\newtheorem{definition}{Definition}

\makeatletter
\newcommand{\thickhline}{%
    \noalign {\ifnum 0=`}\fi \hrule height 1pt
    \futurelet \reserved@a \@xhline
}
\newcolumntype{"}{@{\hskip\tabcolsep\vrule width 1pt\hskip\tabcolsep}}
\makeatother

\usepackage{array,multirow,graphicx}
\usepackage[ruled, lined, linesnumbered, commentsnumbered, longend]{algorithm2e}

\SetKwInOut{KwIn}{Input}
\SetKwInOut{KwOut}{Output}

\begin{document}

\begin{frontmatter}

\title{Markov Chain Generative Adversarial Neural Networks for Solving Bayesian Inverse Problems in Physics Applications}

\author[1,2]{Nikolaj T. M{\"u}cke \corref{mycorrespondingauthor}}
\cortext[mycorrespondingauthor]{Corresponding author}
\ead{nikolaj.mucke@cwi.nl}

\author[1]{Benjamin Sanderse}
\ead{b.sanderse@cwi.nl}

\author[1,3,4]{Sander M. Bohté}
\ead{s.m.Bohte@cwi.nl}

\author[2]{Cornelis W. Oosterlee}
\ead{c.w.oosterlee@uu.nl}

\address[1]{Centrum Wiskunde \& Informatica, Science Park 123, 1098 XG Amsterdam, Netherlands}
\address[2]{Mathematical Institute, Utrecht University, Utrecht, Netherlands}
\address[3]{Swammerdam Institute of Life Sciences (SILS), University of Amsterdam, Amsterdam 1098XH, The Netherlands}
\address[4]{Bernoulli Institute, Rijksuniversiteit Groningen, Groningen 9747AG, The Netherlands}

\begin{abstract}
In the context of solving inverse problems for physics applications within a Bayesian framework, we present a new approach, Markov Chain Generative Adversarial Neural Networks (MCGANs), to alleviate the computational costs associated with solving the Bayesian inference problem. GANs pose a very suitable framework to aid in the solution of Bayesian inference problems, as they are designed to generate samples from complicated high-dimensional distributions. By training a GAN to sample from a low-dimensional latent space and then embedding it in a Markov Chain Monte Carlo method, we can highly efficiently sample from the posterior, by replacing both the high-dimensional prior and the expensive forward map. We prove that the proposed methodology converges to the true posterior in the Wasserstein-1 distance and that sampling from the latent space is equivalent to sampling in the high-dimensional space in a weak sense. The method is showcased on two test cases where we perform both state and parameter estimation simultaneously. The approach is shown to be up to two orders of magnitude more accurate than alternative approaches while also being up to two orders of magnitude computationally faster, in multiple test cases, including the important engineering setting of detecting leaks in pipelines.
\end{abstract} 


\end{frontmatter}


\input{Introduction}

\input{preliminaries}

\input{Markov_Chain_GAN}

\input{Results}

\input{Conclusion.tex}

\section*{Acknowledgment}
This work is supported by the Dutch National Science Foundation NWO under the grant number 629.002.213. which is a cooperative project with IISC Bangalore and Shell Research as project partners. The authors furthermore acknowledge fruitful discussions with Dr. W. Edeling from CWI Amsterdam. 

\section*{CRediT authorship contribution statement}
\textbf{N. M{\"u}cke:} Conceptualization, Methodology, Software, formal analysis, Writing - original draft. \textbf{B. Sanderse:} Writing -review \& editing, formal analysis. \textbf{S. Boht{\'e}:}  Funding acquisition, Writing -review \& editing, supervision. \textbf{C. Oosterlee:} Funding acquisition, formal analysis, Writing -review \& editing, supervision, project administration.

\section*{Declaration of competing interest}
The authors declare that they have no competing financial or personal interests that have influenced the work presented in this paper.

\appendix

\bibliography{references}

\appendix

\input{Appendix.tex}

\end{document}

%% file: Introduction.tex
\section{Introduction} \label{introduction}

The Bayesian inference approach is popular for solving inverse problems in various fields including physics and engineering \cite{asch2016data, harlim2018data, stuart2010inverse, kaipio2006statistical}, mainly due to the fact that it does not only provide an estimate of the solution but also quantifies the uncertainty of the estimate. Information about the distribution of a computed quantity is important, for example, for digital twins \cite{kapteyn2021probabilistic}. 

\begin{figure}[ht]
    \centering
    \includegraphics[width=\textwidth]{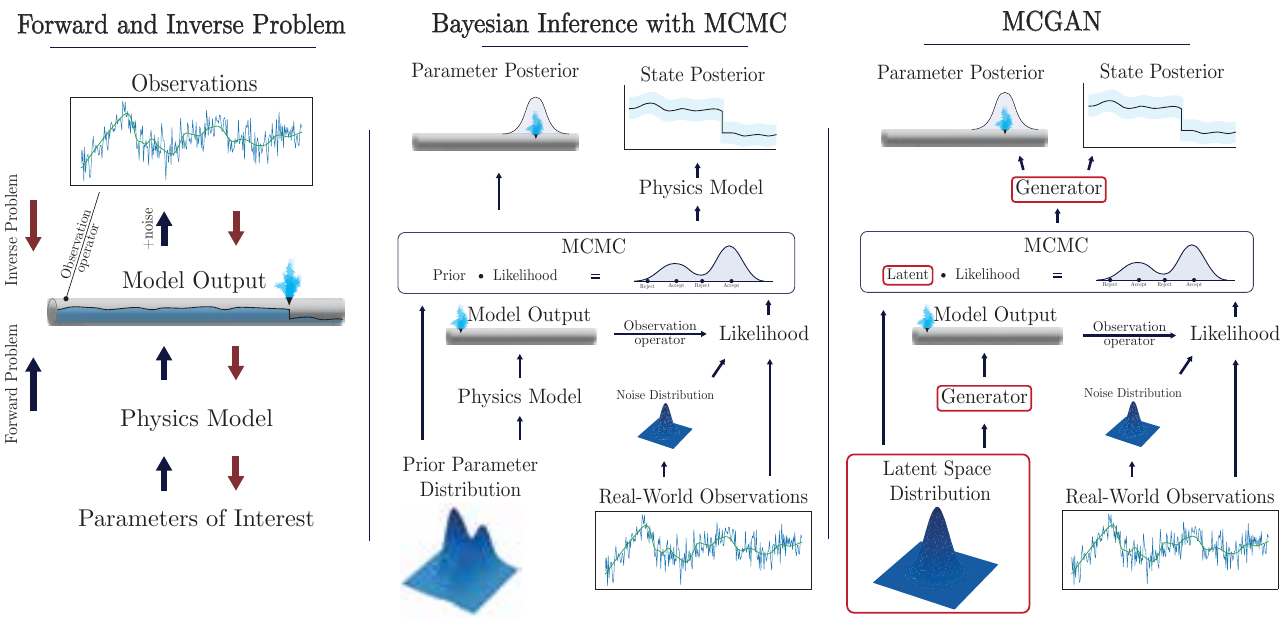}
    \caption{Left: overview of the forward problem and the inverse problem. The parameters of interest are typically boundary and/or initial conditions, or physical parameters. The physics model depends on the system at hand and is here a PDE modeling pipe flow. The model output is the result obtained from a numerical simulation, such as pressure or velocity in the case of fluid dynamics. Observations are either observed from a set of sensors or created synthetically from the model output through the observation operator. Middle: a typical approach for doing Bayesian inference with MCMC (see Section \ref{preliminaries}). Right: our proposed method, the MCGAN approach as explained in section \ref{markov_chain_GAN}. Note that the complicated prior distribution is replaced with a simple latent distribution. Furthermore, the physical model is replaced with a generator that enables us to evaluate the full forward problem, more or less, instantaneously.}
    \label{fig:overview_figure}
\end{figure}

The general idea of Bayesian inference is to use observations to update a given prior distribution towards a resulting posterior distribution over the parameters of interest. The observations and parameters are linked through a forward map and a noise distribution that make up the likelihood function. The main task in the Bayesian approach is to connect the prior and the likelihood in order to compute the posterior distribution. Since the posterior is typically not analytically tractable, one must use numerical sampling techniques such as Monte Carlo methods to approximate the distribution. However, for each sample, it is necessary to compute the likelihood which in turn requires the evaluation of the forward map. For nontrivial problems, such as high-dimensional or nonlinear partial differential equation (PDE) problems, this becomes a computational bottleneck and often results in unacceptable computation times. In Figure \ref{fig:overview_figure}, the general schematics of an inverse problem are shown.

The two most common approaches for overcoming this problem are to either minimize the required number of samples by making certain assumptions about the posterior or to reduce the computational complexity associated with the forward map by approximating it with a surrogate model. The first approach includes methods such as Kalman filters \cite{harlim2018data} and Markov Chain Monte Carlo (MCMC) methods \cite{brooks2011handbook, gamerman2006markov}. With Kalman filters, one minimizes the number of necessary samples by assuming Gaussian distributions. While this is efficient, it is often quite restrictive when it comes to highly nonlinear problems. MCMC methods, while being quite efficient, are based on fewer assumptions but still require many samples. See Figure \ref{fig:overview_figure} for a visualization of a common workflow for Bayesian inference using Markov chain Monte Carlo methods (MCMC). The surrogate modeling approach includes methods such as reduced basis methods \cite{quarteroni2015reduced}, polynomial chaos expansion (PCE) \cite{xiu2010numerical}, and Gaussian processes \cite{wang2018adaptive}. While a surrogate model enables fast likelihood evaluations, it requires a forward map that can be approximated by a low-order representation. This is however not easy for high-dimensional problems and problems involving discontinuities in either the parameters or the state. 

In this paper, we consider an approach that overcomes the above mentioned challenges (high-dimensionality, nonlinearity, discontinuities, expensive sampling) by utilizing machine learning. Specifically, we will make use of neural networks which have already been recognized as promising tools in scientific computing, especially for the case of high-dimensional and nonlinear problems that we wish to address  \cite{brunton2020machine, baker2019workshop, gribonval2021approximation, mucke2021reduced, hesthaven2018non, li2020fourier, kadeethum2021framework}. While there exist several types of neural networks, each aiming at solving specific problems, we focus on generative models in this paper. Generative models aim to learn a distribution from data in order to enable sampling from it at later times \cite{ruthotto2021introduction}. Such models include generative adversarial networks (GANs) \cite{goodfellow2014generative}, variational autoencoders (VAEs) \cite{kingma2013auto}, diffusion models \cite{dhariwal2021diffusion}, and Normalizing Flow models \cite{rezende2015variational}.

Examples where generative models have been successfully used for solving Bayesian inverse problems already exist. In \cite{goh2019solving} and \cite{whang2021composing}, a VAE and Normalizing Flow, respectively, are embedded into a variational Bayesian inference approach and in \cite{patel2020bayesian} and \cite{xia2022bayesian} a GAN and a VAE, respectively, are used as the prior distribution in MCMC sampling. While  \cite{patel2020bayesian} and \cite{xia2022bayesian} combine generative models for parameter prior approximation with MCMC sampling in order to get samples from the posterior distribution, they do not achieve significant speed-ups. Since the generative models are only used to approximate the parameter prior they still need the expensive forward problem being solved to match synthetic observations with the real observations. 

In this paper we focus on GANs due to their success in learning complicated high-dimensional distributions. Specifically, GANs learn a target distribution by training a generator to map latent space samples to samples that mimic a nontrivial high-dimensional target distribution. So, GANs provide a way to represent a complicated high-dimensional distribution by means of a low-dimensional latent space distribution.

We here present the novel Markov Chain Generative Adversarial Network (MCGAN) method, visualized in figure \ref{fig:overview_figure}. In short, we train a GAN to approximate the prior distribution for the states and parameters and thereby obtain a corresponding latent representation. By using an MCMC method, we can then efficiently sample from a latent space posterior instead of the high-dimensional posterior. As a result, we achieve dimensionality reduction, due to the approximation of the desired posterior, and furthermore the forward map is replaced by the generator. In practice, this gives significant computational speed-ups as the computational bottleneck is significantly reduced. The methodology presented draws inspiration from \cite{patel2020bayesian}, but utilizes the GAN in a different manner. Our extension is hence well-suited for both state and parameter estimation in real-time\footnote{Note that ``real-time'' is dependent on the specific problem at hand}. Furthermore, we prove that sampling in the latent space is the same as sampling in the high-dimensional space in a weak sense and we provide a proof of convergence of the posterior distribution in the Wasserstein-1 distance. 

The paper's outline is as follows. In Section \ref{preliminaries}, we explain the setting of Bayesian inverse problems as well as the MCMC methods and GANs. Then, in Section \ref{markov_chain_GAN}, we present the details of our proposed methodology, the MCGAN methodology, including the theoretical findings. In Section \ref{Results}, we show the MCGAN performance on two problems: a stationary Darcy flow and leakage localization in a pipe flow. The results are compared to ensemble Kalman filters and MCMC methods with PCE as the surrogate model. Finally, in Section \ref{Conclusion}, we conclude this work.

%% file: preliminaries.tex
\section{Notation, Problem Setting, and Preliminaries} \label{preliminaries}
Throughout the paper, we will make use of the following notation: capital letters will denote random variables, e.g.\ $X$ and $Y$. The distribution of $X$ is denoted $P_x$, where $P_x(A)=P_x(X\in A)$ is the probability of observing $X\in A$. Similarly, the probability of $x$ is denoted $P_x(x)=P_x(X=x)$. For each probability distribution, there is an associated probability measure, $\mu_x$. Furthermore, we assume that all distributions have a probability density function (PDF), $\rho_x$. The distribution, measure, and PDF are related by:
\begin{align}
    P_x(X\in A) = \mu_x(A) = \int_A\mathrm{d}\mu_x(x) = \int_A\rho_x(x) \mathrm{d}x.
\end{align}
A stochastic variable, $X$, conditioned on another stochastic variable, $Y$, is denoted $X|Y$ and is distributed according to $P_{x|y}(X|Y)$ with PDF, $\rho_{x|y}$.  

\subsection{Problem Setting}
Let $\mathbf{q}\in\mathbb{R}^{N_q}$ denote the state, $\mathbf{m}\in\mathbb{R}^{N_m}$ the parameters, and $\mathbf{y}\in\mathbb{R}^{N_y}$ the available observations. Note that $\mathbf{q}$ encapsulates the state at all discrete times for time-dependent problems. Hence, the state, $\mathbf{q}$, is the full space-time state of the system at hand. $\mathbf{q}$ is computed by solving a forward problem, typically a PDE, depending on the parameters, $\mathbf{m}$. We denote the vector of combined state and parameters, $\mathbf{u}=(\mathbf{q},\mathbf{m})\in\mathbb{R}^{N_u}$, $N_u=N_q+N_m$. 

The inverse problem deals with the recovery of $\mathbf{u}$ from a vector of observations in space-time, $\mathbf{y}\in O\subset \mathbb{R}^{N_y}$. The relation between $\mathbf{u}$ and $\mathbf{y}$ is assumed to be of the form
\begin{align} \label{observations}
    \mathbf{y} = \mathbf{h}(\mathbf{u}) + \eta, \quad \eta \sim P_\eta, \quad \eta \in \mathbb{R}^{N_y},
\end{align}
where $\mathbf{h}:\mathbb{R}^{N_u}\rightarrow O\subset \mathbb{R}^{N_y}$ is referred to as the observation operator and $\eta$ is a random variable denoting the observation or measurement noise.

From Eq.\ \eqref{observations}, we can write the PDF associated with the probability of observing $\mathbf{y}$ given $\mathbf{u}$, $\rho_{y|u}(\mathbf{y} | \mathbf{u})$, as:
\begin{align} \label{likelihood}
    \rho_{y|u}(\mathbf{y} | \mathbf{u}) = \rho_\eta(\mathbf{y} - \mathbf{h}(\mathbf{u})).
\end{align}
When observations are given, one can view this as a function of $\mathbf{u}$, i.e., $\Phi(\mathbf{u})=\rho_{y|u}(\mathbf{y} | \mathbf{u})$, in which case it is referred to as the {\em likelihood} since it is not a PDF with respect to $\mathbf{u}$.
 
We assume that, before observing any data, the probability of $\mathbf{u}$ has the PDF $\rho_0$, which is referred to as the {\em prior}. 
The goal of the Bayesian inverse problem is to identify the PDF, $\rho_{u|y}(\mathbf{u}|\mathbf{y})$, i.e.\ the PDF of  $\mathbf{u}$ given observations, $\mathbf{y}$. Using Bayes theorem, we can write this as:
\begin{align} \label{bayes_theorem}
    \rho_{u|y}( \mathbf{u}  | \mathbf{y}) = \frac{\rho_{y|u}(\mathbf{y}  | \mathbf{u}) \rho_{0}(\mathbf{u})}{\int_{\mathbb{R}^{N_u}} \rho_{y|u}(\mathbf{y} | \mathbf{u}) \rho_0(\mathbf{u}) \: \text{d}\mathbf{u}} = \frac{\rho_\eta(\mathbf{y} - \mathbf{h}(\mathbf{u})) \rho_{0}(\mathbf{u})}{\int_{\mathbb{R}^{N_u}} \rho_{y|u}(\mathbf{y} | \mathbf{u}) \rho_0(\mathbf{u}) \: \text{d}\mathbf{u}}.
\end{align}
The denominator is called the {\em evidence} and serves as a normalization constant; the lefthand side is the {\em posterior}. 

However, in order to compute the likelihood in Eq.\ \eqref{bayes_theorem}, a PDE must be solved for a given set of parameters. Moreover, choosing a suitable prior is not always an easy task, and the evidence can be restrictive to compute in high dimensions.

It should be noted that the last problem is alleviated in many methods such as maximum likehood estimation and MCMC methods as we will describe below. The other two complications will be minimized using our proposed methodology.

\subsection{Markov Chain Monte Carlo Methods} \label{markov_chain_monte_carlo_methods}

MCMC methods \cite{brooks2011handbook, gamerman2006markov} form a class of algorithms for sampling from probability distributions. Stated in terms of the posterior PDF, we have:
\begin{align} \label{general_bayes}
    \rho_{u|y}(\mathbf{u}|\mathbf{y}) \propto \rho_{y|u}(\mathbf{y} |  \mathbf{u}) \rho_{0}(\mathbf{u}),
\end{align}
and we aim to generate a set of points distributed according to the PDF $\rho_{u|y}$.
The general idea is to construct a Markov chain, $\left\{\mathbf{u}_1,\ldots, \mathbf{u}_N \right\}$, with a stationary PDF, $\tilde{\rho}_{u|y}$, that approximates $\rho_{u|y}$. We then sample according to $\tilde{\rho}_{u|y}$ by computing the next element in the chain. 
For MCMC algorithms, we have the following result, under some reasonable assumptions \cite{brooks2011handbook}:
\begin{align} \label{mcmc_convergence}
    \lim_{N_{mcmc}\rightarrow \infty} \frac{1}{N_{mcmc}} \sum_{i=1}^{N_{mcmc}} f(\mathbf{u}_i) = \mathbb{E}_{\mathbf{u}\sim P_{u|y}}[f(\mathbf{u})], \quad \mathbf{u}_i\sim \tilde{P}_{u|y} 
\end{align}
where $\tilde{P}_{u|y}$ is the probability distribution associated with the density $\tilde{\rho}_{u|y}$. Eq.\ \eqref{mcmc_convergence} indicates that with enough samples from the chain, we can approximate some statistics of the true posterior arbitrarily well, i.e.\ the distribution, $\tilde{P}_{u|y}$, converges weakly to $P_{u|y}$.

The arguably most common MCMC sampler is the Metropolis-Hasting (MH) algorithm \cite{metropolis1953equation, hastings1970monte}. However, it is well-known that the MH algorithm converges very slowly in high-dimensional settings. Therefore, in this paper, we make use of the Hamiltonian Monte Carlo Method (HMC), which can be considered a special case of the MH algorithm. Instead of computing new proposals by a random walk, the HMC algorithm computes a new sample by moving in a state space defined by a Hamiltonian ODE system. 

Starting with a `momentum' vector, $\mathbf{p}$, of the same size as $\mathbf{u}$, and a joint PDF $\rho_{u,p}(\mathbf{u},\mathbf{p})$, we define a Hamiltonian as:
\begin{align}
    H(\mathbf{u},\mathbf{p}) = - \log \rho_{u,p}(\mathbf{u},\mathbf{p}) = -\log \rho_{u|p}(\mathbf{u}|\mathbf{p})-\log \rho_{u|y}(\mathbf{u}|\mathbf{y}) \propto \underbrace{\frac{1}{2}\mathbf{p}^T M^{-1}\mathbf{p}}_{=K(\mathbf{p})} - \underbrace{\log \left[ \rho_{y|u}(\mathbf{y}  | \mathbf{u}) \rho_{0}(\mathbf{u}) \right]}_{=U(\mathbf{u})},
\end{align}
where we choose the conditional distribution of the momentum given $\mathbf{u}$ to be normally distributed, $P_{p|u}(\mathbf{p}|\mathbf{u}) \sim \mathcal{N}(0,M)$. $K(\mathbf{p})$ is referred to as the kinetic energy and $U(\mathbf{u})$ the potential energy. One can compute trajectories on level sets of the Hamiltonian by solving the Hamiltonian dynamical system. A new sample is then computed by perturbing the current sample, integrating the Hamiltonian system in time and using the final state as the new sample with acceptance probability:
\begin{align} \label{hamiltonian_alpha}
    \alpha = \min \left\{ 1,  \frac{\exp \left(- H(\mathbf{u}',\mathbf{p}')  \right)}{\exp \left(- H(\mathbf{u}_i,\mathbf{p}(0))\right)}\right\},
\end{align}
where $(\mathbf{u}',\mathbf{p}')$ is the terminal state of the trajectory. Intuitively, this procedure will be biased towards sampling from level sets in the phase space that maximize the likelihood $U(\mathbf{u})$. Furthermore, $K(\mathbf{p})$ ensures that the algorithm explores other areas of the phase space to a degree decided by $M$ and the integration horizon, $T$. Compared to the standard MH algorithm, this reduces the correlation between elements in the chain by traversing long distances in the phase space while maintaining a high acceptance probability due to the energy preserving properties of Hamiltonian dynamics. 

When sampling using the HMC algorithm, a series of choices have to be made, like the number of time steps in the integration and the end time, $T$. If $T$ is too small the sampling will resemble a random walk, while $T$ too large may result in trajectories making a `U-turn' and return to their initial condition. To avoid this, we utilize the No U-Turn Sampler (NUTS) \cite{hoffman2014no}. 

The idea is to integrate backward and forward in time until a U-turn condition is satisfied. Then, a random point from the computed trajectory is chosen, and the algorithm continues from there. 


Even though HMC with NUTS is efficient, one still needs many samples to converge. With a good initial sample, the method converges significantly faster. There are several ways of computing a suitable initial guess, one of which is the maximum a posteriori (MAP) estimate \cite{kaipio2006statistical}, which we will use in this work. This is typically computed using the log PDFs:
\begin{align} \label{MAP}
\mathbf{u}_\mathrm{MAP} = \arg\max_{\mathbf{u}} \log(\rho_{u|y}(\mathbf{y} | \mathbf{u})) + \log(\rho_0(\mathbf{u})).
\end{align}
$\mathbf{u}_\mathrm{MAP}$ can be computed using standard optimization methods such as gradient descent methods. In our case, both $\rho_{y|u}$ and $\rho_u$ are known PDFs so it is easy to compute derivatives using standard software libraries such as PyTorch.

\subsection{Alternative Methods} \label{section:alternative_methods}
Here, we will comment on some well-known alternative methods that exist to speed up solving the Bayesian inverse problem, which we will use to compare our proposed method to. We will also comment on their respective shortcomings. 

\paragraph{Ensemble Kalman Filter} Ensemble Kalman filtering (EnKF) is a Kalman filter variant that is suitable for high-dimensional and nonlinear problems \cite{harlim2018data}. The general idea is to compute the sample mean and sample covariance from an ensemble and then update the prior accordingly. However, as all distributions are assumed to be Gaussian, it means that it is not directly suitable for non-Gaussian problems. In cases with very nonlinear or high-dimensional features, large ensembles are necessary which in turn makes it computationally slow. 

Alternative approaches exist, such as particle filters, that do not assume a Gaussian distribution. However, such methods are, in general, computationally very expensive and will not be further discussed.

\paragraph{Surrogate Models} Instead of replacing the sampling method, the forward computations can be done using a surrogate model, such as polynomial chaos expansion (PCE) methods \cite{xiu2010numerical}, Gaussian processes \cite{stuart2018posterior}, or reduced basis methods \cite{quarteroni2015reduced}. The idea is to approximate the parameter-to-observations map or the forward map by a low-order model that is computationally fast to evaluate. These approaches have been shown to speed up the sampling significantly. However, in high-dimensional cases the curse of dimensionality hampers the applicability of such methods. Moreover, they usually do not perform well in discontinuous and/or highly nonlinear cases unless the approach is tailored to the problem at hand. 

\subsection{Generative Adversarial Neural Networks} \label{generative_adversarial_neural_networks}

In this section, we will give a brief overview of generative adversarial networks (GANs), see \cite{goodfellow2014generative, jabbar2020survey} for more details. We will focus on a version of GANs called Wasserstein GAN (WGAN) \cite{arjovsky2017wasserstein}.

GANs deal with the problem of learning an unknown distribution from samples. Consider a probability distribution, $P_u^r$, on a data space which is a subset of $\mathbb{R}^m$. We aim to approximate $P_u^r$ with another distribution, $P_u^g$. We will refer to $P_u^r$ as the real data probability distribution or the target distribution, and $P_u^g$ the generated distribution. In order to compute $P_u^g$, we define a stochastic latent variable, $Z\in\mathbb{R}^{N_z}$, with prior distribution $P_z^g$, typically chosen to be a Gaussian. Then, we define a generator, $G_\theta:\mathbb{R}^{N_z}\rightarrow \mathbb{R}^{N_u}$, which is a neural network parameterized by its weights, $\theta$. $G_\theta$ takes in the latent variable and outputs $G_\theta(Z)\sim P_u^g$. Hence, $P_u^g = G_{\theta\#} P_z^g$ is the {\em pushforward} of the latent space distribution with PDF $\rho_u^g = \rho_z^g \circ G^{-1}$ \cite{bremaud2020probability}. By choosing $N_z \ll N_u$, we effectively get a low-dimensional representation of the $N_u$-dimensional distribution. Therefore, the variable $z$ can be considered a latent/low-dimensional representation of samples from $P_u^r$.

Next, we introduce the discriminator, $D_\omega:\mathbb{R}^{N_u} \rightarrow \mathbb{R}$. The discriminator takes in samples from either the real data probability distribution or the generated probability distribution, and returns a real number called the {\em score}. A large score means that the discriminator believes the sample comes from the real data distribution. $D_\omega$ is a neural network parameterized by its weights, $\omega$. 

In order to learn the target distribution, a zero-sum game between the generator and the discriminator is set up. The generator aims to maximize the discriminator output, while the discriminator tries to minimize the score of generated samples while simultaneously trying to maximize the score of the real samples. For the WGAN, this game is mathematically formulated as \cite{arjovsky2017wasserstein}:
\begin{align}
    \inf_\theta \sup_{\omega} \quad \mathbb{E}_{X\sim P_u^r}\left[ D_\omega(X) \right] - \mathbb{E}_{Z\sim P_z^g}\left[ D_\omega(G_\theta(Z)) \right].
\end{align}
It can be shown that this inf-sup problem is equivalent to minimizing the Wasserstein-1 distance between $P_u^r$ and $P_u^g$ due to the Kantorovich-Rubinstein duality \cite{arjovsky2017wasserstein}. The WGAN framework requires the discriminator to be Lipschitz continuous with respect to the input. Therefore, we introduce a gradient penalty term to constrain the gradient of the discriminator \cite{gulrajani2017improved}:
\begin{align} \label{WGAN_loss}
    \inf_\theta \sup_{\omega} \quad \mathbb{E}_{X\sim P_u^r}\left[ D_\omega(X) \right] - \mathbb{E}_{Z\sim P_z^g}\left[ D_\omega(G_\theta(Z)) \right] - \lambda \mathbb{E}_{\hat{\mathbf{x}}\sim P_{\hat{X}}}\left[\left( ||\nabla_{\hat{X}} D_\omega(\hat{X})|| - 1 \right)^2\right],
\end{align}
where $\lambda$ is a regularization parameter to be tuned, $\hat{X} = \epsilon X+(1-\epsilon)G_\theta(Z)$, and $\epsilon$ is a small positive number.

In practice, we do not update the weights of the generator and the discriminator at the same time. Instead, we split \eqref{WGAN_loss} into two subproblems: a generator loss that aims to minimize \eqref{WGAN_loss} and a discriminator loss that aims to maximize \eqref{WGAN_loss} by minimizing the negative value:
\begin{subequations}
\begin{align}
    L_G &= - \mathbb{E}_{Z\sim P_z^g}\left[ D_\omega(G_\theta(Z)) \right],  &\text{(Generator Loss)}\\
    L_D &= -\mathbb{E}_{X\sim P_u^r}\left[ D_\omega(X) \right] + \mathbb{E}_{Z\sim P_z^g}\left[ D_\omega(G_\theta(Z)) \right] + \lambda \mathbb{E}_{\hat{X}\sim P_{\hat{X}}}\left[\left( ||\nabla_{\hat{X}} D_\omega(\hat{X})|| - 1 \right)^2\right].  &\text{(Discriminator Loss)} 
\end{align}
\end{subequations}
For details about the training, see \ref{training_WGANS}. The WGAN is visualized in Figure \ref{fig:gan_setup}.

It has been shown that if the generator and the discriminator have sufficient capacity, the generated distribution converges to the real data probability distribution in the Wasserstein-1 distance \cite{liu2017approximation}.

%% file: Markov_Chain_GAN.tex
\section{Markov Chain GAN} \label{markov_chain_GAN}

In this section, we will outline our proposed method, the Markov Chain GAN (MCGAN) method. The general purpose is to combine MCMC methods with GANs in order to perform state and parameter estimation in a computationally fast and accurate way. As mentioned in the previous section, similar approaches exist using polynomial surrogate models \cite{sanderse2021efficient, lu2015limitations} and Gaussian processes \cite{wang2018adaptive}. However, as will be discussed, using GANs gives significant advantages over these alternatives.  

\subsection{Proposed Algorithm}
\begin{figure}
    \centering
    \includegraphics[width=.85\linewidth]{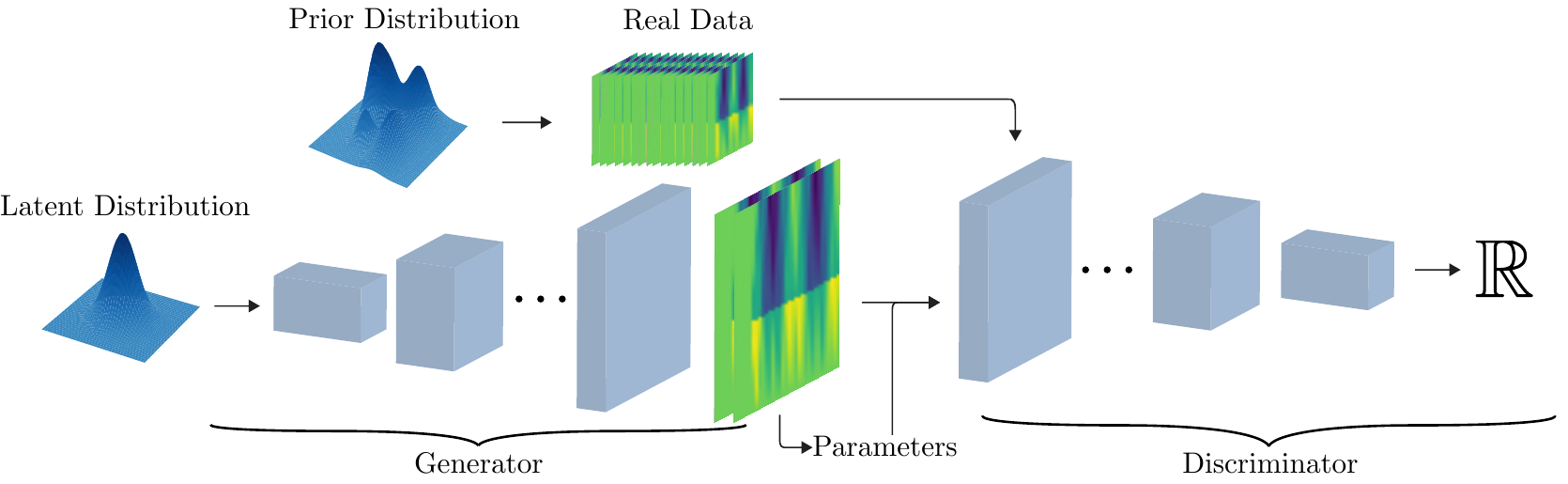}
    \caption{GAN architecture.}
    \label{fig:gan_setup}
\end{figure}

In short, the proposed algorithm aims to speed up posterior sampling without compromising too much on accuracy. The general methodology is to replace the forward model in the likelihood computation with the generator and replace the data prior with the GAN latent distribution (see Figure \ref{fig:gan_setup}). 

Firstly, in an offline stage, we train the GAN to generate discrete solutions to the PDE for the desired time span and corresponding parameters. The GAN is trained on samples from the real prior distribution, $P_0^r$, i.e.\ solutions computed through conventional numerical methods (finite elements, finite volumes, etc.) and aims to learn a generated prior distribution, $P_0^g$. Samples from $P_0^r$ are typically computed by sampling the parameters and then solving the physical model to get the state. This can be a lengthy process, but the training data can be simulated completely in parallel on several computer cores. After training, the (single) generator may generate pairs of states and parameters from a latent sample, $\mathbf{z}$:
\begin{align}
G_\theta(\mathbf{z})=(G_\theta^q(\mathbf{z}),G_\theta^m(\mathbf{z}))=(\mathbf{q}^g,\mathbf{m}^g) = \mathbf{u}^g \sim P_0^g
\end{align}
Hence, the generated distribution is approximating the real prior distribution, $P_0^g\approx P_0^r$. In the training, the discriminator will receive pairs of states and parameters, $(\mathbf{q}^r,\mathbf{m}^r)$, sampled from the real data prior and generated pairs of states and parameters, $(\mathbf{q}^g,\mathbf{m}^g)$, sampled from the generated distribution, in order to ensure that the generator learns to generate states and parameters that match.

\paragraph{Remark} It should be noted that since the GAN is trained on discrete solutions on a specific grid, it will generate samples on the same grid. Therefore, it is important to train on discrete solutions that have all the necessary properties for the online phase (sufficient resolution, etc.)\\
\\
In order to take advantage of the low-dimensional latent space, we need to be able to sample solutions and parameters from a posterior on the latent space, instead of the generated distribution in the full data space. To this end, we have the expression for the latent space posterior density:
\begin{align} \label{latent_space_posterior}
    \rho_{z|y}^g(\mathbf{z}|\mathbf{y}) = \frac{\rho_{y|u}^g(\mathbf{y}|  G_\theta(\mathbf{z})) \rho_{z}^g(\mathbf{z})}{\int_{\mathbb{R}^{N_z}} \rho_{y|u}^g(\mathbf{y} | G_\theta(\mathbf{z})) \rho_z^g(\mathbf{z}) \: \text{d}\mathbf{z}} = 
    \frac{\rho_{\eta}(\mathbf{y}  -  \mathbf{h}(G_\theta(\mathbf{z}))) \rho_{z}^g(\mathbf{z})}{\int_{\mathbb{R}^{N_z}} \rho_{y|u}^g(\mathbf{y} | G_\theta(\mathbf{z})) \rho_z^g(\mathbf{z}) \: \text{d}\mathbf{z}} \propto \rho_{\eta}(\mathbf{y}  -  \mathbf{h}(G_\theta(\mathbf{z}))) \rho_{z}^g(\mathbf{z}).
\end{align}
The forward problem is replaced by an evaluation of the generator at the sampled $\mathbf{z}$. This yields a significant speed-up in online computation time since evaluating the generator is, more or less, instantaneous. The derivation of \eqref{latent_space_posterior} is given in Section \ref{section:latent_space_sampling} in Theorem \ref{theorem:posterior_u_to_z}. 

In the online stage, we use the MCMC method, discussed earlier in Section \ref{markov_chain_monte_carlo_methods}, to sample from the latent space posterior. Given these samples, we can compute derived quantities, e.g.\ for a given quantity of interest $f$, we can compute the expected value by:
\begin{subequations}
\begin{align}
    &\mathbb{E}_{\mathbf{q}\sim P_{q|y}^r}[f(\mathbf{q})] \approx \mathbb{E}_{\mathbf{q}\sim P_{q|y}^g}[f(\mathbf{q})] =  \mathbb{E}_{\mathbf{z}\sim P_{z|y}^g}[f(G^q_{\theta}(\mathbf{z}))] \approx \frac{1}{N_{\text{MCMC}}}\sum_{i=1}^{N_{\text{MCMC}}} f(G^q_{\theta}(\mathbf{z}_i)), \quad \mathbf{z}_i \sim \tilde{P}_{z|y}^g, \\
    &\mathbb{E}_{\mathbf{m}\sim P_{m|y}^r}[f(\mathbf{m})] \approx \mathbb{E}_{\mathbf{m}\sim P_{m|y}^g}[f(\mathbf{m})] =  \mathbb{E}_{\mathbf{z}\sim P_{z|y}^g}[f(G^m_{\theta}(\mathbf{z}))] \approx \frac{1}{N_{\text{MCMC}}}\sum_{i=1}^{N_{\text{MCMC}}} f(G^m_{\theta}(\mathbf{z}_i)), \quad \mathbf{z}_i \sim \tilde{P}_{z|y}^g, 
\end{align}
\end{subequations}
By choosing an appropriate $f$, we can thereby compute various quantities of interest by sampling the latent space posterior. See Theorem \ref{theorem:posterior_u_to_z} in Section \ref{section:latent_space_sampling} for details.

For an overview of the methodology, see Algorithm \ref{MCGAN_offline} for the offline stage and Algorithm \ref{MCGAN_online} for the online stage. Here is a summary of the distinct advantages of the proposed method compared to the alternatives discussed in Section \ref{section:alternative_methods}:
\begin{itemize}
    \item The latent vector $\mathbf{z}$ is, in general, of significantly lower dimension than the state and parameters, effectively reducing the dimension of the stochastic space and making it computationally faster to sample from the distribution using MCMC methods;
    \item The computationally expensive forward problem is replaced by the generator, whose cost is computationally negligible to evaluate once it has been trained;
    \item Since the forward map is replaced by a neural network, derivatives of the log-likelihood function can be computed efficiently, which enables computationally fast MAP estimation.
\end{itemize}
While the advantages are clear, it is worth mentioning the drawbacks as well:
\begin{itemize}
    \item There is no immediate way of choosing the dimension of the latent space.   However, one can consider it a hyperparameter and perform hyperparameter optimization;
    \item Training a GAN is not always an easy task, since commonly known problems of training neural networks, such as local minima and generalization, also apply here;
    \item It is necessary to generate much training data in order to ensure accuracy of the GAN. 
\end{itemize}
Note that the drawbacks are not unique to this methodology, but general when dealing with neural networks. The first two points are a matter of hyperparameter tuning, and the last point is a matter of time in the offline stage. Furthermore, with an efficient numerical solver and the fact that the offline stage can be easily parallellized (since the training samples are independent), the generation of data is often feasible within a reasonable timeframe. If the forward solver would be too expensive to allow for this, it is recommended to first obtain a simpler forward model e.g.\ by model reduction techniques such as reduced order models. 

\paragraph{Remark} It is important to note that the purpose of the proposed methodology is to solve Bayesian inverse problems computationally fast in an online stage. As mentioned, this comes at a cost of an expensive offline stage in which the GAN is trained on simulated data. However, when the GAN is trained it can be deployed in several settings. Therefore, the method is highly suitable in settings where computational speed is crucial and offline training time is less important. This is, for example, the case for digital twins and model predictive control where repeated real-time state estimation and parameter calibration are necessities.

\begin{algorithm}[]
 \KwIn{$N_{\text{train}}$, GAN hyperparameters, GAN architecture}
 Generate training samples, $\{(\mathbf{q}_i,\mathbf{m}_i) \}_{i=1}^{N_{\text{train}}}\sim P_0^r$, by solving the forward problem\;
 Train the GAN to approximate the prior, $P_0^g\approx P_0^r$ (see Section \ref{generative_adversarial_neural_networks})\;
 \KwOut{Trained generator, $G_\theta$}
 \caption{MCGAN offline stage}
 \label{MCGAN_offline}
\end{algorithm}

\begin{algorithm}[h] \label{alg:online_stage}
 \KwIn{Generator from Algorithm \ref{MCGAN_offline}, MCMC parameters, observations $(\mathbf{y})$, }
 Compute the MAP estimate in the latent space, using a gradient descent algorithm, as initial sampling point (see Eq. \ref{MAP})\;
 Use MCMC algorithm (see Section \ref{markov_chain_monte_carlo_methods}) to sample from the latent space posterior (see Eq. \eqref{latent_space_posterior})\;
 Generate states and parameters from the posterior latent space samples:
 \begin{align*}
     \{\mathbf{q}_i, \mathbf{m}_i\} = \{G_{\theta}^q(\mathbf{z}_i), G_{\theta}^m(\mathbf{z}_i)\}, \quad i=1,\ldots,N_{\text{samples}}  \quad \mathbf{z}_i \sim \tilde{P}_{z|y}^g.
 \end{align*} \\
 Compute the relevant statistics, such as mean and variance\;
 \KwOut{$\{\mathbf{q}_i, \mathbf{m}_i \}^{N_{\text{samples}}}_{i=1}$, statistics}
 \caption{MCGAN online stage}
 \label{MCGAN_online}
\end{algorithm}

\subsection{Latent Space Sampling} \label{section:latent_space_sampling}
Here we prove that sampling from the latent space posterior yields the same results as sampling from the full data space in a weak sense. From \cite{bogachev2007measure}, we have the following results for push forward distributions:
\begin{align} \label{expected_change_of_variables}
    \mathbb{E}_{U\sim P_u}[f(U)] = \int_{E} f(\mathbf{u}) \rho_u(\mathbf{u}) \mathrm{d} \mathbf{u} = \int_{G^{-1}(E)} f(G(\mathbf{z})) \rho_z(\mathbf{z}) \mathrm{d} \mathbf{z} = \mathbb{E}_{Z\sim P_z}[f(G(Z))],
\end{align}
where $\mathbf{u}=G(\mathbf{z})\in E$, $P_u$ and $P_z$ are the distributions of $\mathbf{u}$ and $\mathbf{z}$ with PDFs $\rho_u$ and $\rho_z$, respectively, and $f$ is a measurable function on $E$. Here, $P_u=G_{\#} P_z^g$ is the push forward of $P_z$ by $G$. The derivation of Eq. \eqref{expected_change_of_variables} only requires that $G$ is measurable.

\begin{theorem} \label{theorem:posterior_u_to_z}
Let $G_\theta$ be a generator. Let $Z$ be a latent space variable distributed according to a latent space distribution, $P_z^g$, with PDF $\rho_z^g$, and let $U=G_\theta(Z)$ be distributed according to the push forward distribution of the latent space distribution, $U\sim P_u^g = G_{\theta \#} P_z^g$, with PDF $\rho_0^g$. Then, the push forward posterior distribution, conditioned on data $\mathbf{y}$, is equal to the latent space posterior distribution conditioned on the same data in a weak sense, i.e.\ for all measurable functions $f$, the following holds:
\begin{align} \label{posterior_u_to_z}
    \mathbb{E}_{U\sim P_{u|y}^g}[f(U)] = \mathbb{E}_{Z\sim P_{z|y}^g}[f(G(Z))],
\end{align}
where the associated PDFs are given by: 
\begin{align}
    \rho_{u|y}^g(\mathbf{u}|\mathbf{y}) =  \frac{\rho_{y|u}^g(\mathbf{y} |  \mathbf{u}) \rho_0^g(\mathbf{u})}{\int_{\mathbb{R}^{N_u}} \rho_{y|u}^g(\mathbf{y} | \mathbf{u}) \rho_0^{g}(\mathbf{u}) \: \text{d}\mathbf{u}}, 
    \quad
    \rho_{z|y}^g(\mathbf{z}|\mathbf{y}) =  \frac{\rho_{y|u}^g(\mathbf{y}  |  G_\theta(\mathbf{z})) \rho_{z}^g(\mathbf{z})}{\int_{\mathbb{R}^{N_z}} \rho_{y|u}^g(\mathbf{y} | G_\theta(\mathbf{z})) \rho_z^g(\mathbf{z}) \: \text{d}\mathbf{z}}.
\end{align}
\end{theorem}
 \noindent The proof can be found in \ref{appendix_posterior_u_to_z_proof}.

Using Theorem \ref{theorem:posterior_u_to_z}, we can conclude that sampling from the latent space posterior, $P_{z|y}^g$, and pushing forward using the generator, $G_\theta$, yields the same results as sampling directly from the generated posterior, $P_{u|y}^g$ in a weak sense. 
While small perturbations in the latent domain might result in different output in the full data space, Theorem \ref{theorem:posterior_u_to_z} shows that, in a weak sense, the posterior obtained from pushing forward the latent samples is equal to the full data posterior. 

\subsection{Convergence of Generated Posterior}
In this subsection, we prove that $P_{u|y}^g\approx P_{u|y}^r$ when $P_0^g\approx P_0^r$ and under some additional reasonable assumptions on the generator.  That is, the case where the prior is approximated in the Wasserstein metric and the likelihood is approximated with a surrogate forward map. While \cite{sprungk2020local} proves the cases where either the likelihood or the prior is approximated, we provide a proof where both are being approximated.

Before stating the theorem, we need to define the appropriate spaces and metrics. Let $(E, d_E)$ be a complete metric space. $E\subset \mathbb{R}^d$ is the set containing the state and parameter vectors, $\mathbf{u}\in E$ and $d_E: E\times E \rightarrow \mathbb{R_+}$ assigns non-negative distances between two elements of $E$. Furthermore, in this formulation, the observation operator, $\mathbf{h}:E\rightarrow O$ maps elements from $E$ to the observation space. 

We can then define the relevant space of probability distributions:
\begin{definition}
On a metric space, $(E,d_E)$, we define the space of probability distributions as:
\begin{align*}
    \mathcal{W}_q(E) = \left\{ P \: : \: |P|_{\mathcal{W}_q} < \infty   \right\}, \quad |P|_{\mathcal{W}_q} = \inf_{x_0\in E}\left(\int_E d_E(x,x_0)^q \rho(x) \: \mathrm{d}x\right)^{1/q}.
\end{align*}
\end{definition}
Then, we define the Wasserstein-1 distance and its dual representation \cite{panaretos2020invitation}:
\begin{definition}
For two probability distributions, $P_1,P_2 \in \mathcal{W}_1(E)$, the Wasserstein-1 distance is defined as:
\begin{align*}
    W_1(P_1,P_2) = \inf_{\gamma\in \Gamma(P_1,P_2)} \left| \int_E\int_E d_E(x,y) \gamma(x,y) \mathrm{d}x\mathrm{d}y   \right|,
\end{align*}
where $\Gamma(P_1,P_2)$ is the set of joint PDFs, $\gamma$, for combined probability distributions with $P_1$ and $P_2$ as marginal distributions, respectively. From the Kantorovich–Rubinstein duality, we can write the Wasserstein-1 distance as:
\begin{align*}
    W_1(P_1,P_2) = \sup_{\mathrm{Lip}(f)\leq 1} \left| \int_{E} f(x) \rho_1(x) \mathrm{d}x - \int_{E} f(x) \rho_2(x) \mathrm{d}x  \right|,
\end{align*}
where $f:E\rightarrow \mathbb{R}$ is a Lipschitz continuous function, $\mathrm{Lip}(f)$ is its corresponding Lipschitz constant, and $\rho_1$ and $\rho_2$ are the PDFs of $P_1$ and $P_2$, respectively.
\end{definition}
Besides the Wasserstein distance, we will also be working with the weighted norms:
\begin{align*}
    ||f||_{L^1_\rho} = \int |f(x)| \rho(x) \: \mathrm{d}x, \quad ||f||_{L^2_\rho} = \left(\int |f(x)|^2 \rho(x) \: \mathrm{d}x \right)^{1/2},
\end{align*}
and the finite dimensional norm:
\begin{align}
    ||\mathbf{v}||_{l^2} = \left(\sum_{i=1}^{N} v_i^2 \right)^{1/2}, \quad \mathbf{v}=(v_1,\ldots,v_N)\in \mathbb{R}^{N}.
\end{align}
With the proper spaces, norms, and metrics defined, we can state the following theorem inspired by \cite{sprungk2020local}:
\begin{theorem}\label{main_theorem} Let $(E,d_E)$ be a bounded metric space with $\sup_{\mathbf{x}_1,\mathbf{x}_2\in E}d_E(\mathbf{x}_1,\mathbf{x}_2) \leq D < \infty$.

Let $P_0^r\in \mathcal{W}_2(E)$ denote the prior probability distribution of real data, and let $P_0^r\in \mathcal{W}_2(E)$ denote the generated prior probability distribution of generated data. 

Let the real data and generated likelihoods satisfy
\begin{align*}
    &\rho_{y|u}^r(\mathbf{y}|\mathbf{u}) \propto \Phi^r(\mathbf{u}) = e^{-l^r(\mathbf{u})}, \quad \Phi^r:E\rightarrow \mathbb{R}_+, \quad l^r:E\rightarrow \mathbb{R}_+,\\
    &\rho_{y|u}^g(\mathbf{y}|\mathbf{u}) \propto \Phi^g(\mathbf{u}) = e^{-l^g(\mathbf{u})}, \quad \Phi^g:E\rightarrow \mathbb{R}_+, \quad l^g:E\rightarrow \mathbb{R}_+,
\end{align*}
where $l^r$ and $l^g$ are the log-likelihood functions for the real and generated data, respectively, and $\Phi^r$ and $\Phi^g$ are Lipschitz continuous functions with Lipschitz constants $\mathrm{Lip}(\Phi^r)$ and $\mathrm{Lip}(\Phi^g)$, respectively. Furthermore, let $\Phi^r, \Phi^g \in L_{\rho_0^g}^2$, where $L_{\rho_0^g}^2$ is the weighted $L^2$ space with $\rho_0^g$ as the weight function.

Assume the GAN has converged, i.e.\ 
$W_1(P_0^r, P_0^g) \leq \epsilon_1.$
Furthermore, assume that this implies convergence of the log-likelihood, as follows,
\begin{align}
||l^r(\mathbf{u}) - l^g(\mathbf{u}) ||_{L^1_{\rho_0^g}}\leq \epsilon_2, \quad ||l^r(\mathbf{u})-l^g(\mathbf{u}) ||_{L^2_{\rho_0^g}} \leq \epsilon_3,
\end{align}
where $||\cdot||_{L^1_{\rho_0^g}}$ is the weighted $L^1$-norm with $\rho_0^g$ as the weight function. Then, the Wasserstein-1 distance between the real posterior probability distribution given observations and the generated posterior probability distribution given observations satisfies:
\begin{align}
    W_1(P_{u|y}^r, P_{u|y}^g) \leq C_1 \epsilon_1  + C_2 \epsilon_2 + C_3 \epsilon_3,
\end{align}
where
\begin{align*}
    C_1 = \frac{(1+D\mathrm{Lip}(\Phi^r))}{Q^r_u(\mathbf{y})}, \quad C_2 = \frac{\max(\Phi^r,\Phi^g)}{Q^r_u(\mathbf{y})Q^g_u(\mathbf{y})} (1+D\mathrm{Lip}(\Phi^r)) |P_0^g|_{\mathcal{W}_1}, \quad C_3 = \frac{\max(\Phi^r,\Phi^g)}{Q^g_u(\mathbf{y})} |P_0^g|_{\mathcal{W}_2},
\end{align*}
where $Q^r_u$ and $Q^g_u$ are the evidence from the real and generated posterior, respectively, and $D$ denotes the maximum distance between two points in the metric space, $E$. 
\end{theorem} 
\noindent The proof can be found in \ref{appendix_main_theorem_proof}.

In short, the proof of Theorem \ref{main_theorem} helps us understand when we can expect convergence of the posterior. While the assumptions of the Theorem might seem restrictive, this is actually not the case. Firstly, we assume that the metric space, $E$, is bounded, which is typically the case in many applications. Secondly, we assume that the likelihood is of the form $e^{-l(\mathbf{u})}$, Lipschitz continuous, and is in the weighted $L^1$ and $L^2$ spaces. For simple observation operators (e.g.\ linear) this is a consequence of the negative log-likelihood function typically being an $L^2$-norm.  This leaves us with the question of the convergence of the prior which directly determines $\epsilon_1$, $\epsilon_2$, and $\epsilon_3$.

\subsection{Convergence of the Generated Prior} \label{convergence_of_generated_prior}
The convergence of the generated prior is a matter of studying convergence properties of GANs. Such studies are beyond the scope of this paper. Instead, we refer to \cite{arjovsky2017wasserstein, gulrajani2017improved, liu2017approximation} where convergence properties of GANs are discussed. In short, ensuring convergence of GANs is similar to ensuring convergence of other types of neural networks. Hence, it is a matter of having enough data and performing hyperparameter tuning. For the MCGAN, the amount of data is, in general, not a problem, as we simulate the training data. 


%% file: Results.tex
\section{Results} \label{Results}

In this section, we will present the results on two different problems using the MCGAN methodology. We show two distinct parameter and state estimation cases to highlight various advantages of using the MCGAN methodology. Firstly, we consider a Darcy flow case (stationary flow through a porous medium), with the aim of approximating the horizontal and vertical velocity, the pressure, and the permeability field. The purpose of this case is to emphasize the ability to deal with high-dimensional stochastic problems as the permeability field is spatially distributed and follows a high-dimensional distribution. Secondly, we consider the problem of leakage detection in pipe flow. Here, the challenge lies in dealing with a nonlinear hyperbolic PDE with discontinuities and a non-informative prior. 

The results will be assessed using the relative root mean squared error (RRMSE):
\begin{align}
    \text{State RRMSE} = \frac{\sqrt{\sum_{i=1}^{N_q}(\mathbf{q}^*_i - \mathbf{q}_i)^2}}{\sqrt{\sum_{i=1}^{N_q}\mathbf{q}^2_i}}, \quad \text{Parameter RRMSE} = \frac{\sqrt{\sum_{i=1}^{N_m}(\mathbf{m}^*_i - \mathbf{m}_i)^2}}{\sqrt{\sum_{i=1}^{N_m}\mathbf{m}^2_i}},
\end{align}
where $\mathbf{q}^*$ and $\mathbf{m}^*$ denote the approximated state and parameters, respectively, and $\mathbf{q}$ and $\mathbf{m}$ are the reference state and parameters, respectively. The reference values are computed using an appropriate numerical solver. These will be discussed in each test case.

Furthermore, we will look at the approximated posterior distributions resulting from the MCGAN. 

For the details on the training and hyperparameters for each of the test cases as well as GAN architectures, see \ref{training_WGANS}. Furthermore, all the training data for the GANs are generated by sampling the parameter spaces according to the chosen distribution for the test case. The number of training samples is chosen based on the performance of the resulting GAN. Note that it is, in general, a difficult problem to choose the number of necessary training samples.

The specific architectures of the generators and discriminators for each test case can be found in Figure \ref{fig:GAN_networks}. It is worth noting that we make use of convolutional neural networks in all cases due to their success in problems dealing with spatially distributed degrees of freedom \cite{mucke2021reduced, brunton2020machine}.

\begin{figure}[h]
     \centering
     \begin{subfigure}[t]{0.48\textwidth}
         \centering
         \includegraphics[width=0.75\textwidth]{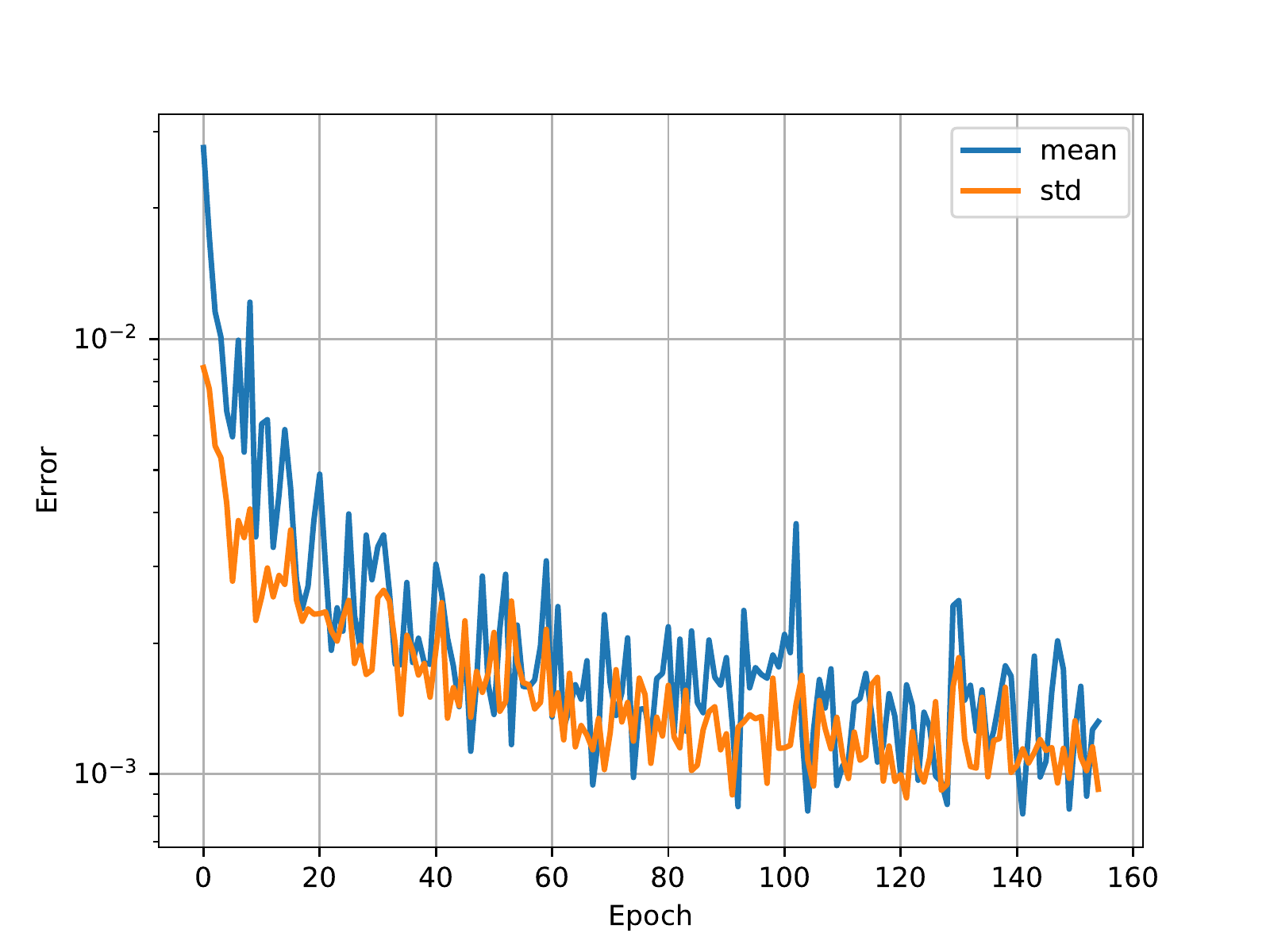}
         \caption{Darcy flow.}
         \label{fig:moment_convergence_darcy}
     \end{subfigure}
     \hfill
     \begin{subfigure}[t]{0.48\textwidth}
         \centering
         \includegraphics[width=0.75\textwidth]{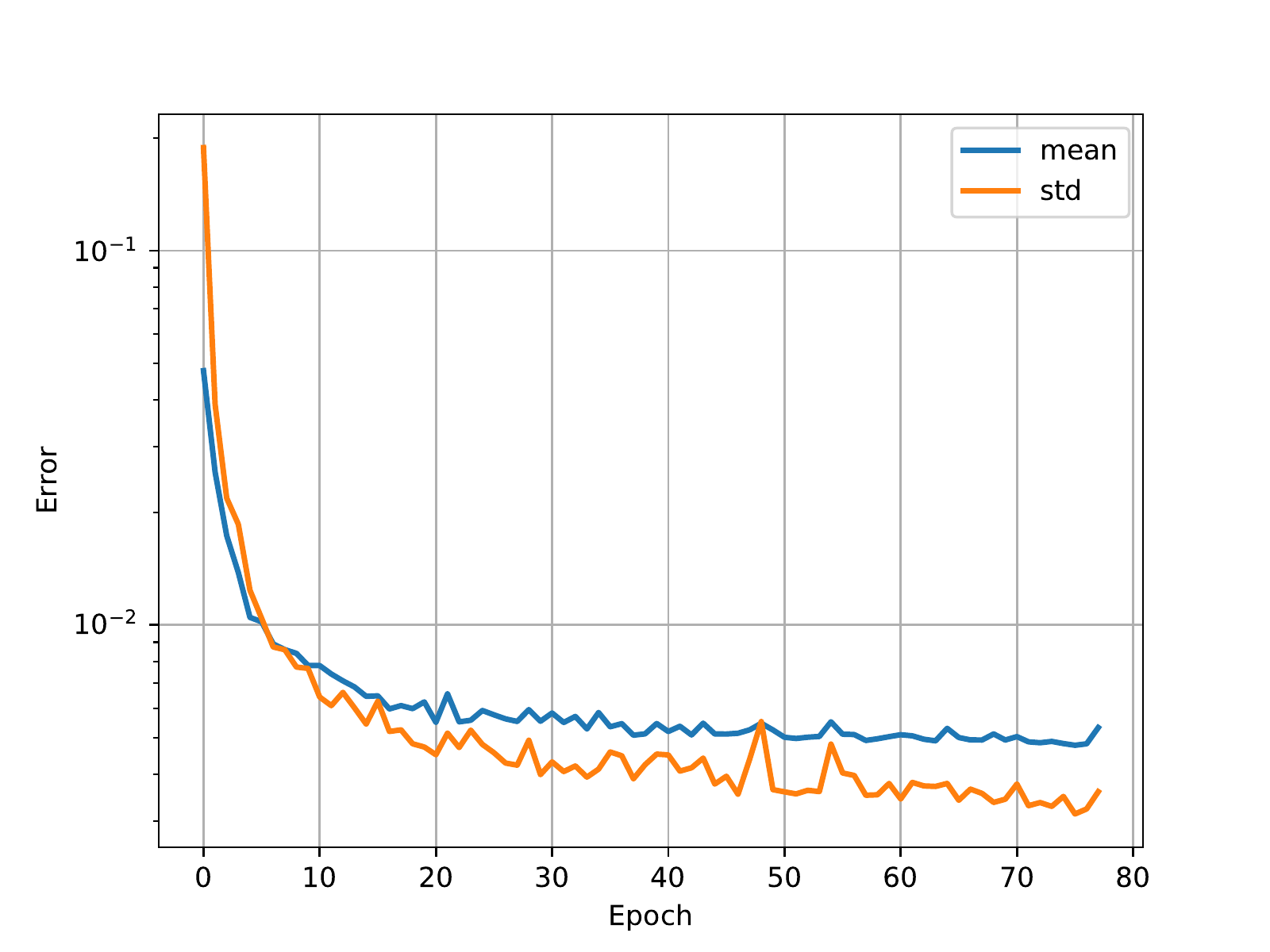}
         \caption{Pipeflow.}
         \label{fig:moment_convergence_pipeflow}
     \end{subfigure}
     \caption{Convergence of mean and variance of the generated prior towards the prior computed from simulations. The error is computed using a test dataset.}
      \label{fig:gan_convergence_figures}
\end{figure}

As mentioned in Section \ref{convergence_of_generated_prior}, it is not feasible to compute the Wasserstein distance for very high-dimensional distributions. Therefore, in order to show convergence of the generated prior, we show the convergence of the first two moments, mean and standard deviation, with the training epochs. Here, we have a value for the mean and variance at every grid point and we compute the error as the relative RMSE. The convergence plots are shown in Figure \ref{fig:gan_convergence_figures}. While this is a weaker type of convergence than convergence in the Wasserstein-1 distance, it still gives an indication that the GAN error is sufficiently small for the purpose of Bayesian inversion. 

All results are generated using synthetic observations. Therefore, all observations are simulation-based and perturbed with artificial noise. To ensure that we are not subject to inverse crime \cite{colton1998inverse}, the synthetic observations are generated with a higher resolution than what is used for the training of the GAN and PCE models, for all experiments. Furthermore, the Kalman filter results are also generated with a lower resolution. Secondly, we will use another distribution for the likelihood function than for the noise in the synthetic observations. The specifics will be discussed in each test case.

\subsection{Darcy Flow}
As a first test case, we consider stationary two-dimensional Darcy flow:
\begin{subequations} \label{darcy_equations}
\begin{alignat}{2}
    \mathbf{v} + k \nabla p &= 0, \quad &\mathbf{x}\in [0,1]^2, \\
    \nabla \cdot \mathbf{v} &= 0, &\mathbf{x}\in [0,1]^2, \\
    p &= 1, &\mathbf{x}\in 0\times[0,1], \\
    p &= 0, &\mathbf{x}\in 1\times[0,1], \\
    \mathbf{v}\cdot \mathbf{n} &= 0, &\mathbf{x}\in [0,1]\times\{0,1\}. 
\end{alignat}
\end{subequations}
$p:[0,1]^2\rightarrow \mathbb{R}$ denotes pressure, $\mathbf{v}:[0,1]^2\rightarrow \mathbb{R}^2$ denotes the velocity, $k:[0,1]^2\rightarrow \mathbb{R}$ is the spatially-dependent permeability field, and $\mathbf{x}=(x_1,x_2)$ are the spatial coordinates in the horizontal and vertical directions, respectively. The permeability field $k$ is modeled as a lognormal field, $\log k = m \sim \mathcal{N}(0,C)$. The problem of state and parameter estimation for Darcy flow is often considered in data assimilation and in uncertainty quantification, see e.g.\ \cite{domesova2017solution, ruchi2019transform}.

The covariance matrix $C$ is derived from the class of Mat{\'e}rn functions \cite{kumar2018multigrid}:
\begin{align}
    C(\mathbf{x}_i,\mathbf{x}_j) = \sigma^2\frac{2^{1-\nu}}{\Gamma(\nu)} \left(\frac{\sqrt{2\nu}}{l}d(\mathbf{x}_i,\mathbf{x}_j) \right)^{\nu} K_{\nu}\left(\frac{\sqrt{2\nu}}{l}d(\mathbf{x}_i,\mathbf{x}_j) \right).
\end{align}
$\Gamma$ is the gamma function and $K_\nu$ is the modified Bessel function of the second kind. $d(\mathbf{x}_i,\mathbf{x}_j)$ denotes the distance between two points, $\mathbf{x}_i$ and $\mathbf{x}_j$, in the domain, $\nu$ defines the smoothness, $\sigma^2>0$ is the variance, and $l>0$ is the correlation length. 

We denote by $\mathbf{m}_N$ the discretized version of $m$ defined on an $N\times N$ grid and the covariance matrix, $C_N\in\mathbb{R}^{N^2}\times\mathbb{R}^{N^2}$, has elements $(C_N)_{i,j} = C(\mathbf{x}_i,\mathbf{x}_j)$. Then, $\mathbf{m}_N$ can be sampled by computing
\begin{align}
    \mathbf{m}_N = \sum_{i=1}^{N^2} \sqrt{\lambda_i} \hat{\mathbf{m}}_i  \boldsymbol{\psi}_i, \quad \hat{\mathbf{m}}\in \mathbb{R}^{N^2}, \quad \hat{\mathbf{m}}\sim \mathcal{N}(0,I),
\end{align}
where $I\in \mathbb{R}^{N^2}\times\mathbb{R}^{N^2}$ is the identity matrix, $\lambda_i$ are the eigenvalues of $C_N$ in descending order, and $\boldsymbol{\psi}_i$ the corresponding eigenvectors. Hence, the permeability field is determined by $\hat{\mathbf{m}}_i$, $i=1,\ldots,N$. A reduced representation of the permeability field can then be computed by choosing $n<N^2$:
\begin{align}
    \mathbf{m}_N^{(n)} = \sum_{i=1}^{n} \sqrt{\lambda_i} \hat{\mathbf{m}}_i  \boldsymbol{\psi}_i.
\end{align}
Thereby, the reduced permeability field is determined by $n$, instead of $N^2$, parameters. 

For generating the training data, Eq.\ \eqref{darcy_equations} is solved using the finite element method. The velocity is discretized by discontinuous Raviart-Thomas elements of polynomial order 3 and the pressure is discretized by Lagrange elements of polynomial order 2. This is known to be a stable pairing of finite element spaces for the stationary Darcy flow \cite{cockburn2009superconvergent}. The domain is divided into $32\times 32$ squares, each divided into two triangles, resulting in 25793 degrees of freedom in total. The solutions are then evaluated on a $50\times 50$ equidistant grid. The implementation is done using the FEniCS library \cite{logg2012automated}.

The specific setting for creating the permeability field here is $n=1089$, $\nu=1.5$, $l=0.2$, and $\sigma=0.5$.

For the observations, we consider evenly distributed sensors at locations, $(\mathbf{x}_1,\ldots,\mathbf{x}_{N_y})$, measuring the horizontal velocity at $N_y=100$ discrete points, see Figure \ref{fig:darcy_true_state}. Thus, $\mathbf{h}:\mathbb{R}^{N\times N} \rightarrow \mathbb{R}^{N_y}$, and the measurements are created by:
\begin{align}
    \textbf{y} = \mathbf{h}(\mathbf{v}) + \eta, \quad \mathbf{h}(\mathbf{v}) = (v_1(\mathbf{x}_1),\ldots,v_1(\mathbf{x}_{N_y})) \quad \eta \sim \mathcal{N}(0,0.01^2 I), \quad \eta\in \mathbb{R}^{N_y}.
\end{align}
The synthetic observations are generated using $50\times 50$ squares divided into two triangles. The velocity is discretized with polynomial order 4 and the pressure with polynomial order 3. The test case is similar to the one presented in \cite{domesova2017solution}.

We compare the MCGAN method with the ensemble Kalman inversion method \cite{ding2021ensemble}. We do not compare with PCE since it is infeasible to compute a PCE model for a problem of this high dimensionality.

\subsubsection*{GAN setup}
The discriminator of the GAN consists of convolutional layers and the generator consists of transposed convolutional layers. The generator is trained to generate the velocity in the horizontal direction, $v_1$, the velocity in the vertical direction, $v_2$, the pressure, $p$, and the log-permeability field, $\log(k)$. Each quantity is considered a channel in the sense of convolutional neural networks. Thereby, the generator outputs tensors of the shape $(4,N,N)$. To avoid boundary artifacts in the generated fields originating from the transposed convolutional layers, the generator is trained to generate fields of the shape $(4,N+l,N+l)$, $l>0$, which are then cropped to the desired size. For details on the exact architecture specifications, see Figure \ref{fig:GAN_networks}. 

\subsubsection*{Results}
The MCGAN results are computed with a single chain of 20,000 samples, where the first 12,500 samples are discarded to ensure that we only use samples with a converged chain. The MAP estimate is used as the initial MCMC sample, which reduces the time until convergence for the MCMC method significantly. For the likelihood function, we use $\mathcal{N}(0,0.02^2 I)$, which is different from the distribution used to generate the observation noise. 

\begin{figure}
     \centering
     \begin{subfigure}[b]{0.32\textwidth}
         \centering
         \includegraphics[width=\textwidth]{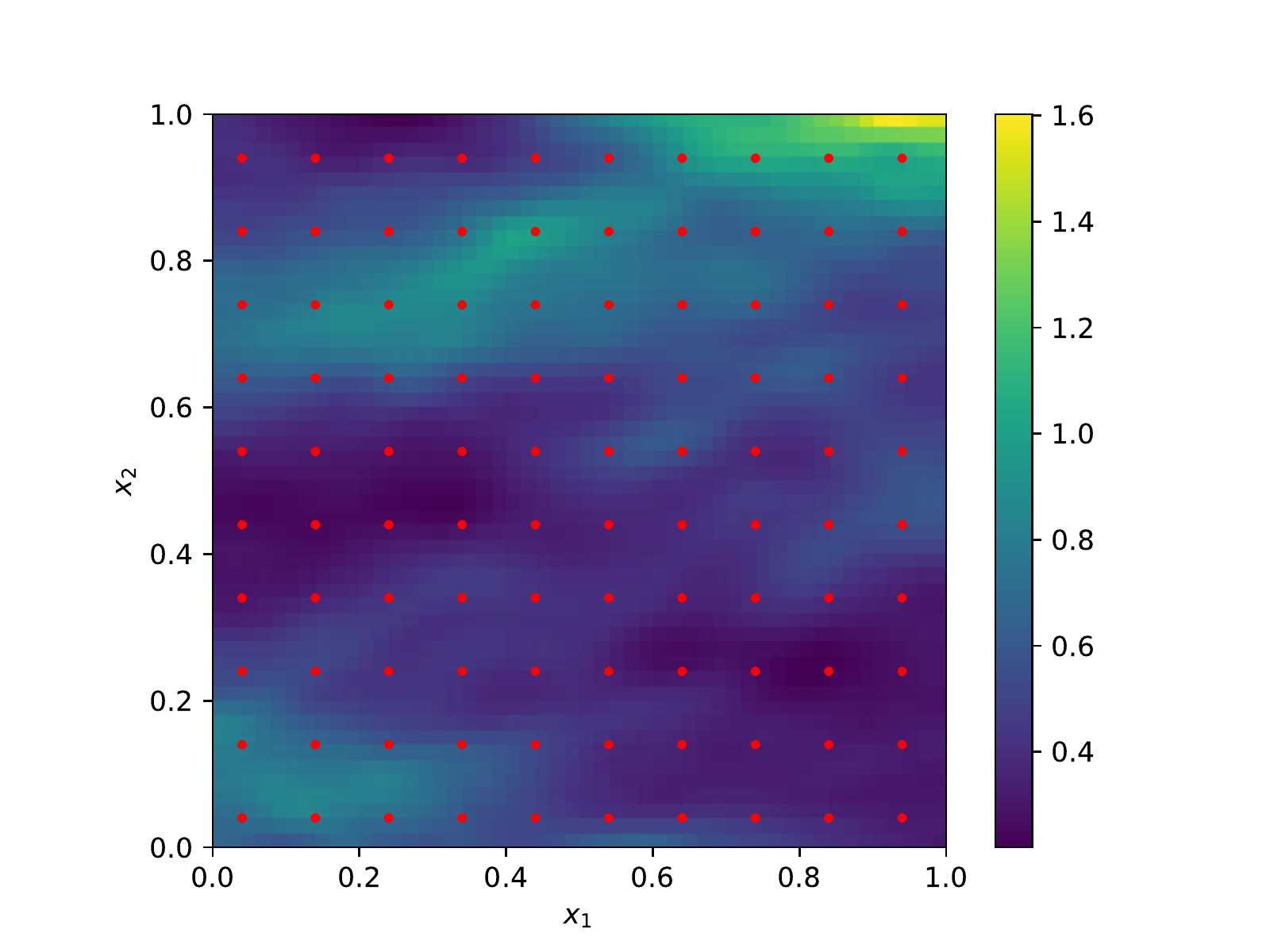}
         \caption{True $v_1$.}
         \label{fig:darcy_true_state}
     \end{subfigure}
     \hfill
     \begin{subfigure}[b]{0.32\textwidth}
         \centering
         \includegraphics[width=\textwidth]{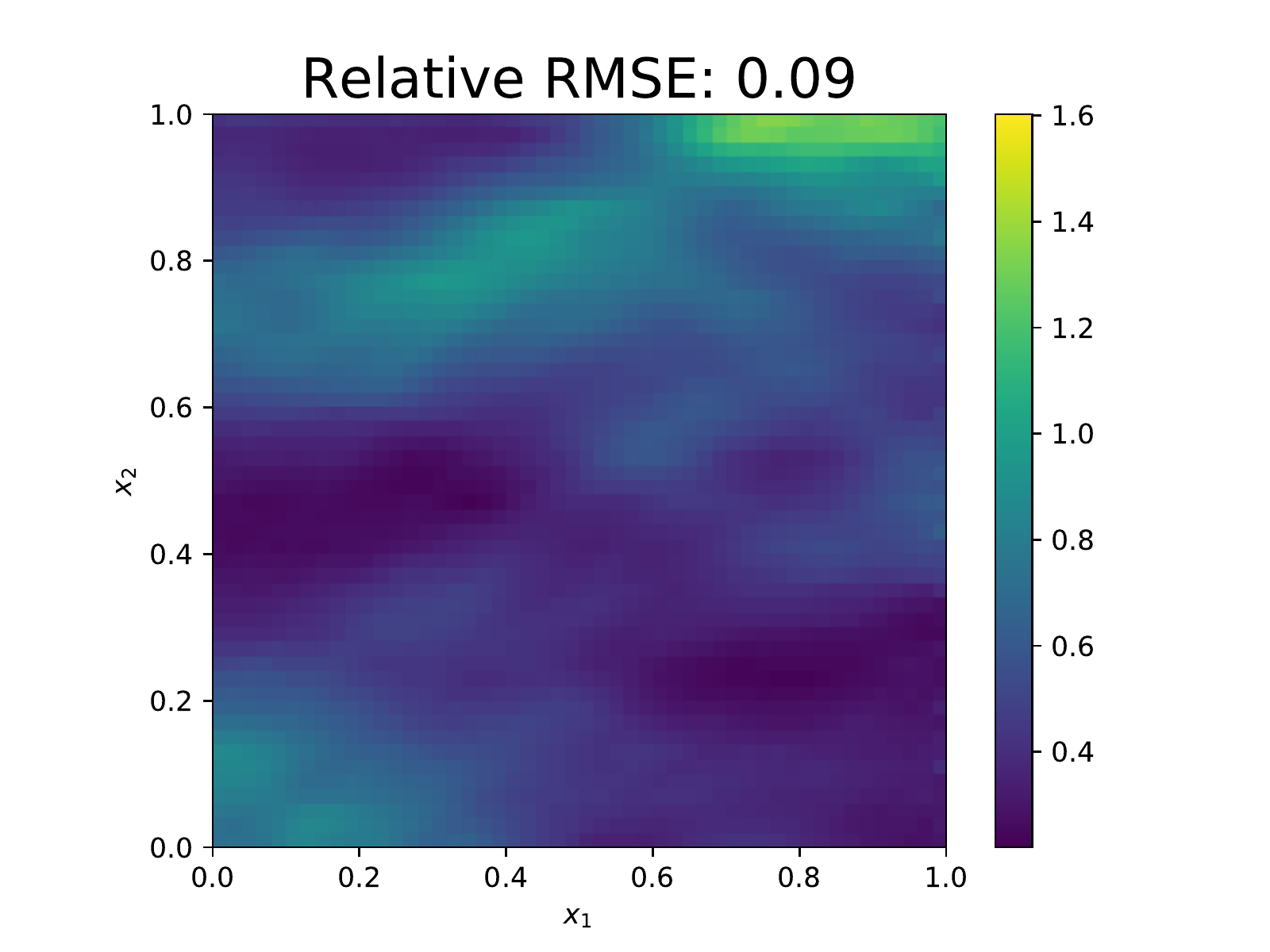}
         \caption{MCGAN approximated $v_1$.}
         \label{fig:darcy_MCGAN_state}
     \end{subfigure}
     \hfill
     \begin{subfigure}[b]{0.32\textwidth}
         \centering
         \includegraphics[width=\textwidth]{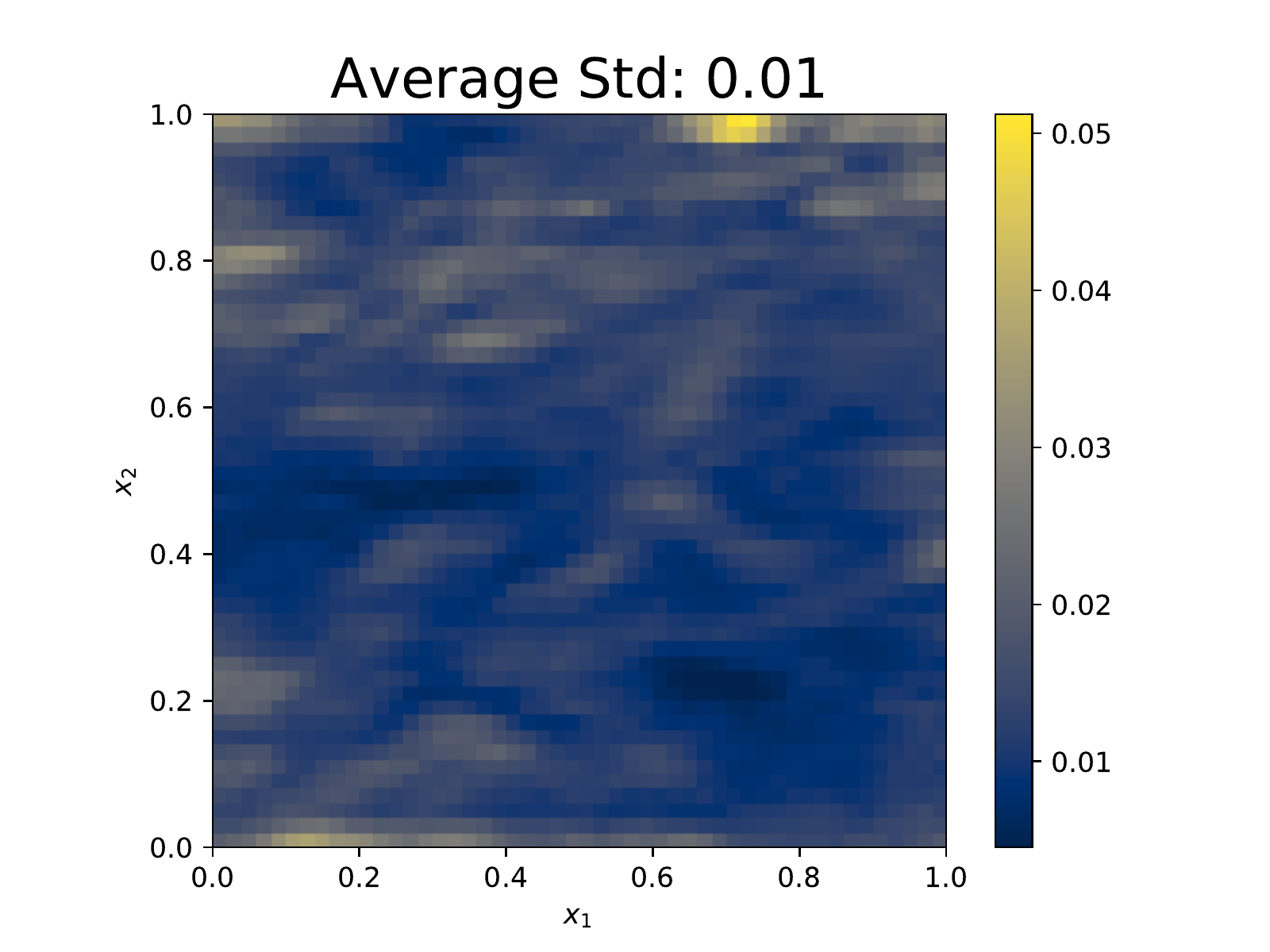}
         \caption{MCGAN standard deviation $v_1$.}
         \label{fig:darcy_MCGAN_state_std}
     \end{subfigure}
    \begin{subfigure}[b]{0.32\textwidth}
         \centering
         \includegraphics[width=\textwidth]{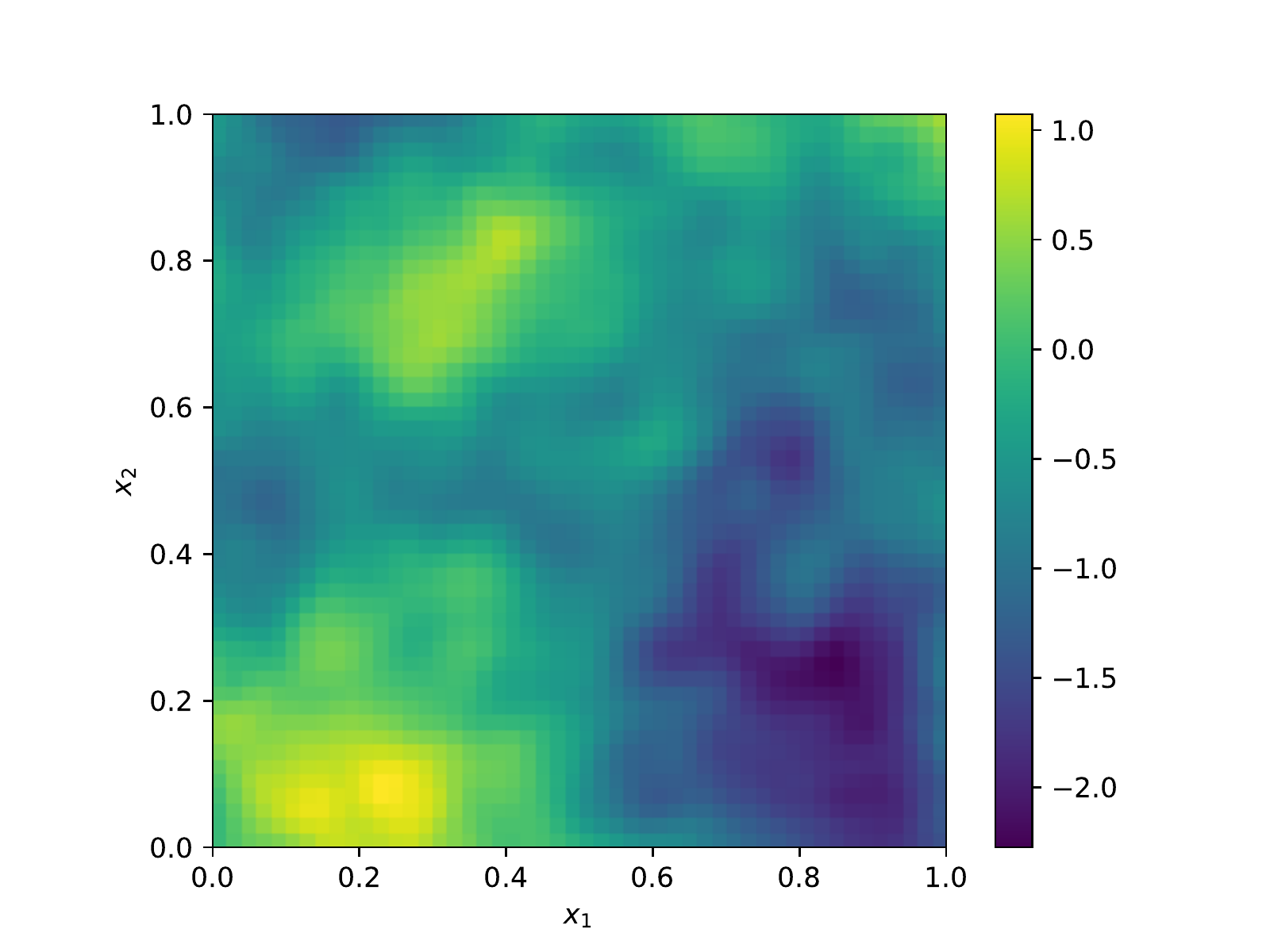}
         \caption{True $\log(k)$.}
         \label{fig:darcy_true_permeability}
     \end{subfigure}
     \hfill
     \begin{subfigure}[b]{0.32\textwidth}
         \centering
         \includegraphics[width=\textwidth]{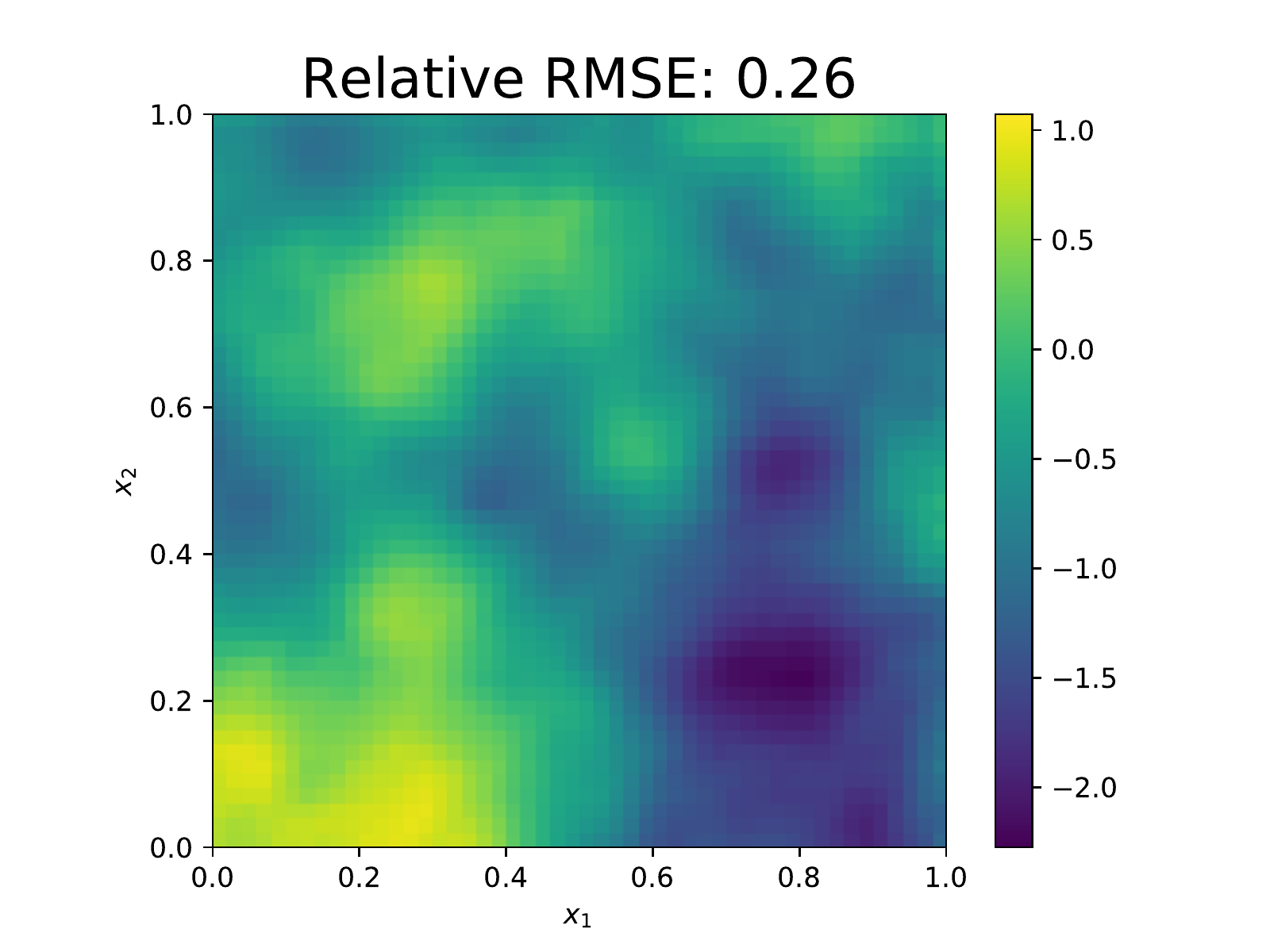}
         \caption{MCGAN approximated $\log(k)$.}
         \label{fig:darcy_MCGAN_permeability}
     \end{subfigure}
     \hfill
     \begin{subfigure}[b]{0.32\textwidth}
         \centering
         \includegraphics[width=\textwidth]{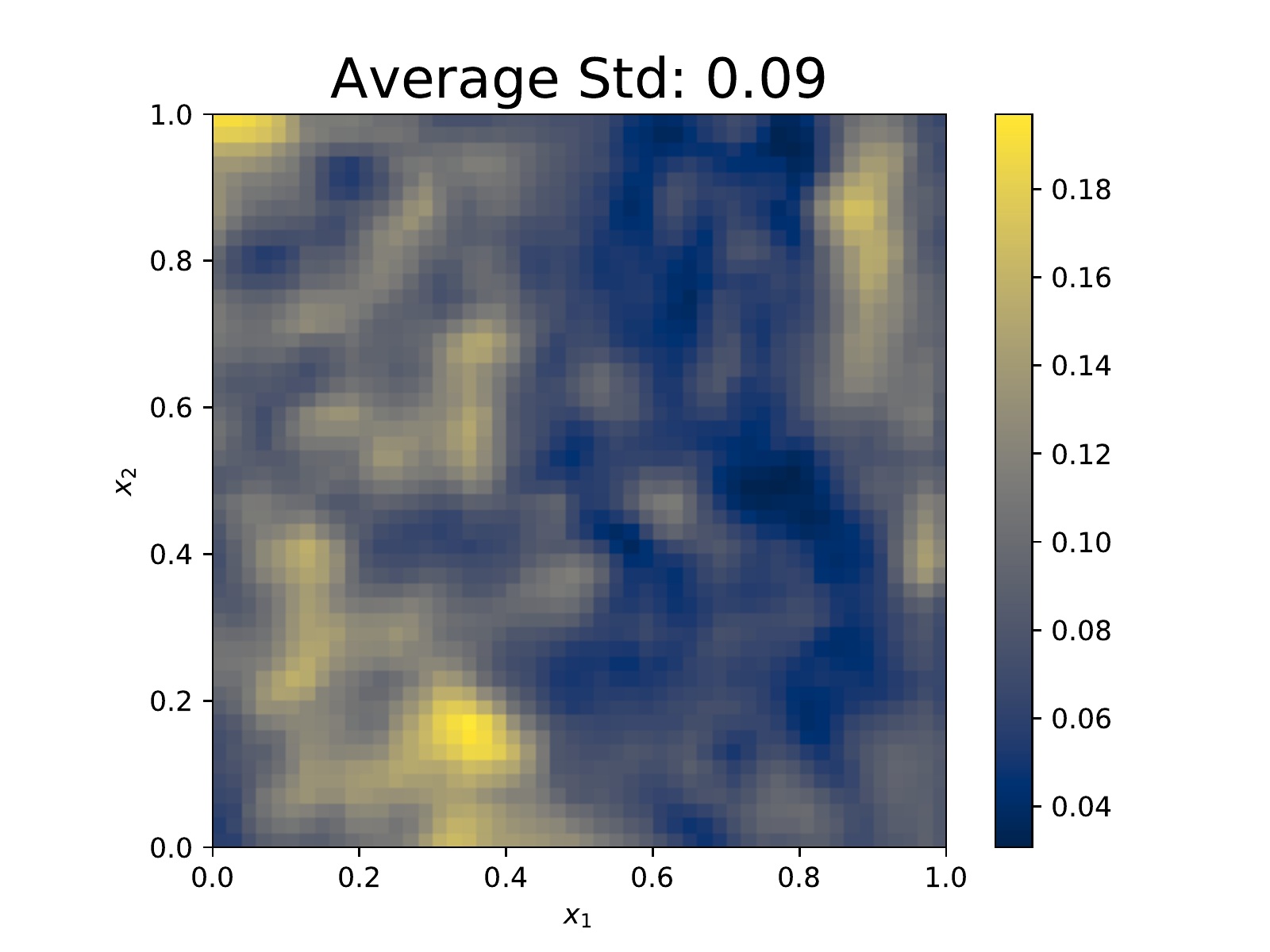}
         \caption{MCGAN standard deviation $\log(k)$.}
         \label{fig:darcy_MCGAN_parameter_std}
     \end{subfigure}
        \caption{Top row: $v_1$. Red dots are points of measurements. Bottom row: $\log(k)$.}
        \label{fig:darcy_results}
\end{figure}

In Figure \ref{fig:darcy_results}, the results from using MCGAN for the Darcy flow are shown. We see that the horizontal velocity is estimated rather accurately with a relative RMSE of $0.09$ and a relatively low standard deviation. Not surprisingly, the standard deviation seems to be largest at the upper boundary where no measurements are available. Furthermore, we see larger uncertainty in the areas of the domain where the magnitude $v_1$ is large.

Regarding the log-permeability, the MCGAN captures the structure of the true log-permeability as well as the sharp edges with a relative RMSE of 0.26.

For both the state and log-permeability, the Kalman filter gives similar, but slightly worse, accuracy and significantly smoother results than the MCGAN approach (see Figure \ref{fig:darcy_results_kalman}). Hence, the Kalman filter is not able to capture the sharper edges. Furthermore, it is an order of magnitude slower (see Table \ref{tab:computation_time}).

\subsection{Leakage Detection in Pipe Flow}
To show the method's generality, as a much different problem, we consider unsteady single phase flow through a pipeline, until suddenly (at t=10s) a leak occurs. As a consequence, pressure waves start propagating through the pipeline, and the velocity field at the leak becomes discontinuous because of the mass flow leaving through the leak. We have only two measurement locations, one close to the inlet and one close to the outlet of the pipeline measuring pressure, and the goal is to infer the leak location and size based on these measurements and a physical model of the flow in the pipeline. This is a challenging problem because of the very sparse measurement data and the discontinuity in the solution.

The governing equations are given by the one-dimensional Euler equations for mass and momentum conservation \cite{kundu2002fluid, hauge2007model}:
\begin{subequations} \label{pipe_equations}
\begin{align}
&\partial_t q_1 +  \partial_x q_2= C_d \sqrt{\rho(p(\rho)-p_{\text{amb}})} \delta(x-x_l)H(t-t_l), \label{mass_equatiosn} \\
&\partial_t q_2 +  \partial_x \left(\frac{q_2^2}{q_1} + p(\rho)A\right) = -\frac{1}{2d}\frac{q_2^2}{q_1} f_f(q), \label{momentum_equation}\\ 
&v(0,t) = v_0, \quad p(L,t) = p_L.
\end{align}
\end{subequations}
where $\rho$ is the fluid density (not to be confused with the probability density functions in previous sections), $p(\rho) = c^2 (\rho - \rho_0) + p_0$ is the pressure, $v$ is the velocity, $q_1=\rho A$, $q_2=\rho v A$, $\delta$ is the Dirac delta function, and $H$ is the Heaviside function. $v_0$ represents the boundary conditions prescribed on the velocity at the left end of the pipe and $p_L$ is the prescribed pressure at the right end of the pipe. $d$, $A$, $p_{\text{amb}}$, $c$, $\rho_0$, are all constants. The physical quantities they represent and the values we will be working with are in Table \ref{tab:pipeflow_parameters}. The righthand side in Eq.\ \eqref{mass_equatiosn} is the leakage, modeled as a discharge. $t_l$ is the time at which the leakage occurs, $x_l$ and $C_d$ are the two parameters of interest. They represent the location and size of the leakage, respectively. The righthand side of Eq.\ \eqref{momentum_equation} is the friction, where $f_f$ is the Darcy-Weisbach friction coefficient, which is given by the Haaland expression \cite{schetz1996handbook}:
\begin{align}
    \frac{1}{\sqrt{f_f}} = -\frac{1}{4}1.8\log_{10}\left[\left(\frac{\varepsilon/D}{3.7}\right)^{1.11}+\frac{6.9}{Re} \right],
\end{align}
where $Re$ is the Reynolds number, $Re = \frac{\rho v d}{\mu}$, with $\mu$ the fluid viscosity and $\varepsilon$ the pipe roughness. The values and units of all parameters in the model are in Table \ref{tab:pipeflow_parameters}. The initial condition is $(q_1,q_2)=(\rho_0 A,\rho_0 v_0 A)$.

\begin{table}[]
    \centering
    \caption{Parameters for the pipe flow equations, \eqref{pipe_equations}. Note that the discharge coefficient and the leakage location have values denoted by intervals, as they are the parameters to determine.}
    \begin{tabular}{llll} 
    \toprule
    Physical quantity & Constant &  Value & Unit\\ 
    \midrule
     Pipe length & $L$ & 2000 & m \\
     Diameter & $d$ & 0.508 & m \\
     Cross-sectional area & $A$ & 0.203 & $\text{m}^2$ \\
     Speed of sound in fluid & $c$ & 308 & m/s \\
     Ambient pressure & $p_{\text{amb}}$ & 101325 & Pa \\
     Reference pressure & $p_{\text{ref}}$ & 5016390 &Pa \\
     Reference density & $\rho_{\text{ref}}$ & 52.67 & kg/$\text{m}^3$  \\
     Inflow velocity & $v_0$ & 4.0 & m/s \\
     Outflow pressure & $p_L$ & 5016390 & Pa \\
     Pipe roughness & $\varepsilon$ & $10^{-8}$ & m \\
     Fluid viscosity & $\mu$ & $1.2\cdot 10^{-5}$ & $\text{N}\cdot\text{s}/\text{m}^2$ \\
     Leakage start time & $t_l$ & 10 & s \\
     Discharge coefficient & $C_d$ & $\left[1.0\cdot 10^{-4}, 9.0\cdot 10^{-4}\right]$ & m  \\
     Leakage location & $x_l$ & $[100,1900]$ & m \\ 
     \bottomrule
\end{tabular}
    \label{tab:pipeflow_parameters}
\end{table}

Eq.\ \eqref{pipe_equations} is solved using the nodal discontinuous Galerkin method \cite{hesthaven2007nodal}. We use Legendre polynomials for the modal representation of the local polynomials, and Lagrange polynomials for the nodal representation. The numerical flux is chosen to be the Lax-Friedrichs flux. To ensure stability and non-oscillatory behavior while ensuring high-order accuracy, a TVBM slope-limiter is applied after each time step \cite{hesthaven2007nodal}. The time-stepping is performed using the BDF2 method, with an initial implicit Euler step \cite{leveque2007finite}.

For the generation of the training data, we consider 75 elements with a local polynomial order of 3. The resulting solution is then evaluated on an equidistant grid consisting of 256 points. For the time-stepping, we consider a horizon of $T=64$ seconds with 256 time steps. Hence, $(q_1,q_2)\in\mathbb{R}^{256\times 256}\times \mathbb{R}^{256\times 256}$.

We assume a uniform prior for both the leakage location, $x_l\sim \mathcal{U}(100,1900)$, and the discharge coefficient, $C_d\sim \mathcal{U}\left( 1.0\cdot 10^{-4}, 9.0\cdot 10^{-4}\right)$. Other choices of distributions of $x_l$ and $C_d$ are subject to future studies.

For the state and parameter estimation only measurements of the pressure are observed. We consider the vector, $(x_1,\ldots,x_{N_y})$, of measurement locations, and the vector of measurement times, $(t_1,\ldots,t_{N_y})$. This gives rise to the synthetic observations:
\begin{align}
    \textbf{y} = \mathbf{h}(p) + \eta, \quad \mathbf{h}(p) = (p(x_1,t_1),\ldots,p(x_{N_y},t_{N_y})), \quad \eta \sim \mathcal{N}(0,1500^2 I), \quad \eta \in \mathbb{R}^{N_y}
\end{align}
We specifically consider the case where we only observe at $x=20m$ and at $x=1980m$ and for all time instances, i.e.\ $(x_1,\ldots,x_{N_y})=(20,\ldots,20,1980,\ldots,1980)$ and $(t_1,\ldots,t_{N_y}) = (0.25,\ldots,64,0.25,\ldots,64)$. Hence, $N_y=2\cdot 256=512$. For simulating the synthetic observations, we used 100 elements with a local polynomial order of 4.

\subsubsection*{GAN setup}
As for the above two test cases, we use convolutional layers for the discriminator and transposed convolutional layers for the generator. The GAN is trained to generate the velocity, $v$, and pressure, $p$, instead of generating the conservative variables $q_1$ and $q_2$, since $v$ and $p$ are the quantities of interest. The GAN is trained to generate full space-time solutions in the intervals, $x\in [0,L]$ and $t\in [0,T]$. $v$ and $p$ are considered channels in the sense of convolutional neural networks. Hence, the generator generates tensors of size $(2,256,256)$. 

At the location of the leakage, there will be a discontinuity in the velocity, due to a drastic drop in the velocity. We use this information to compute the leakage location by identifying the spatial location of the discontinuity, by convolving the state with an appropriate kernel. Furthermore, a dense neural network takes in the generated state and outputs the discharge coefficient. See Figure \ref{fig:GAN_networks} for a visulization of the GAN.

Due to the large differences in orders of magnitude, the velocity, pressure and discharge coefficients are scaled to have values between -1 and 1.

\subsubsection*{Results}
The MCGAN results are computed with a single chain of 15,000 samples, where the first 10,000 samples are discarded to ensure that we only use samples after the chain has converged. The MAP estimate is, again, used as the initial MCMC sample in order to speed up convergence. For the likelihood function, we use $\mathcal{N}(0,3000^2 I)$, which is different from the distribution used to generate the observation noise. 

\begin{figure}
     \centering
     \begin{subfigure}[b]{0.32\textwidth}
         \centering
         \includegraphics[width=\textwidth]{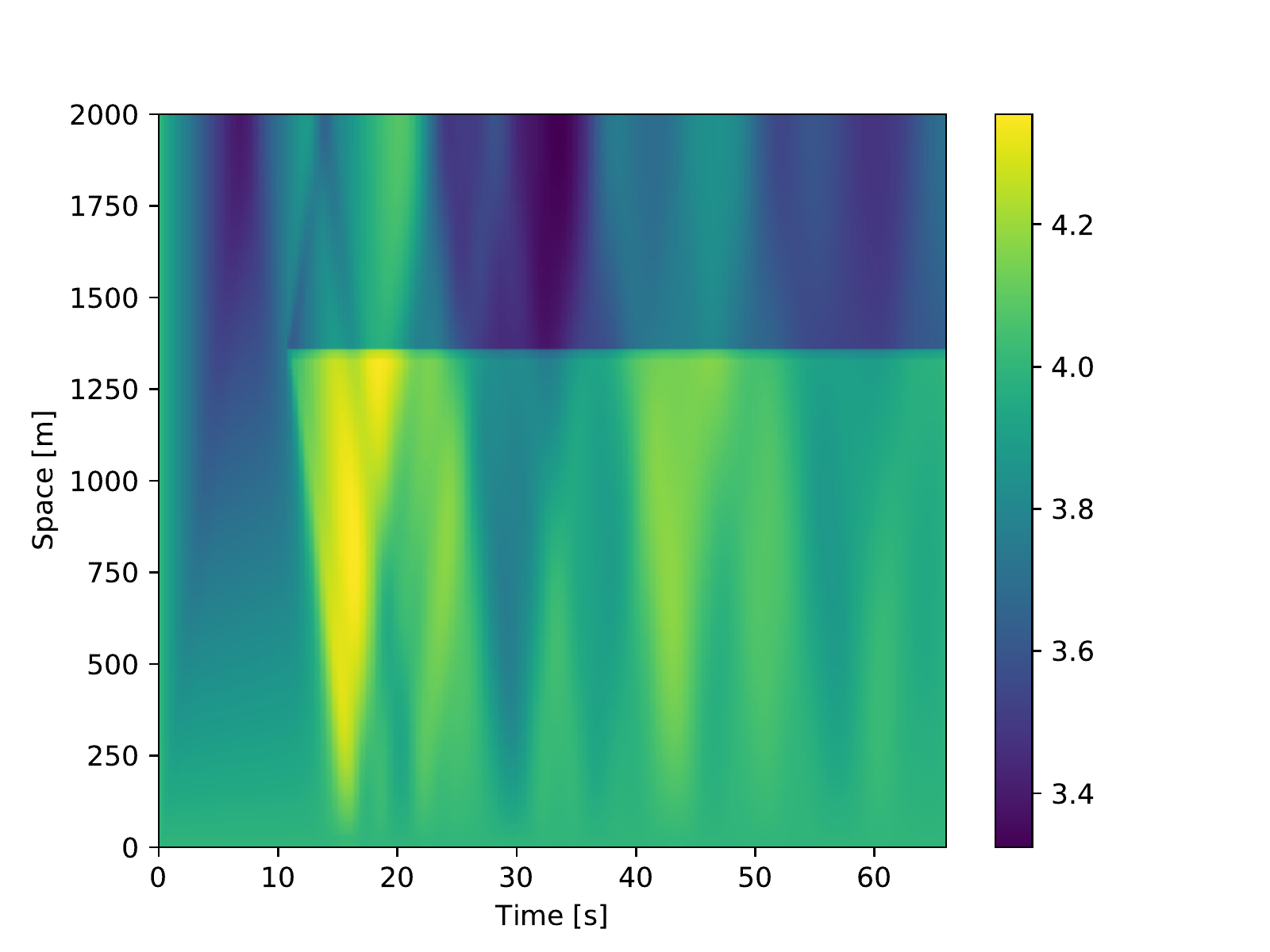}
         \caption{True velocity.}
         \label{fig:pipe_true_state}
     \end{subfigure}
     \hfill
     \begin{subfigure}[b]{0.32\textwidth}
         \centering
         \includegraphics[width=\textwidth]{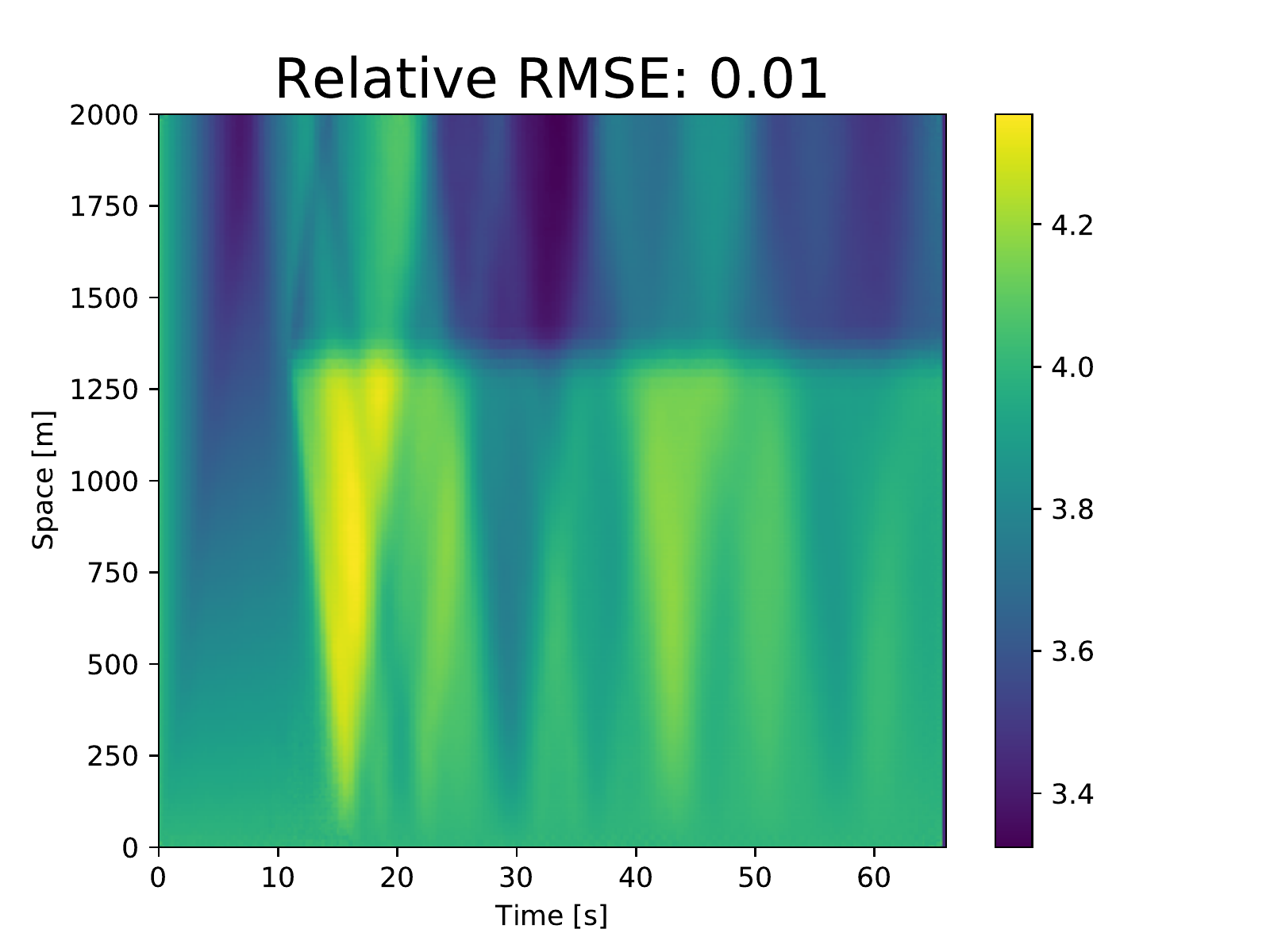}
         \caption{MCGAN approximated velocity.}
         \label{fig:pipe_MCGAN_state}
     \end{subfigure}
     \hfill
     \begin{subfigure}[b]{0.32\textwidth}
         \centering
         \includegraphics[width=\textwidth]{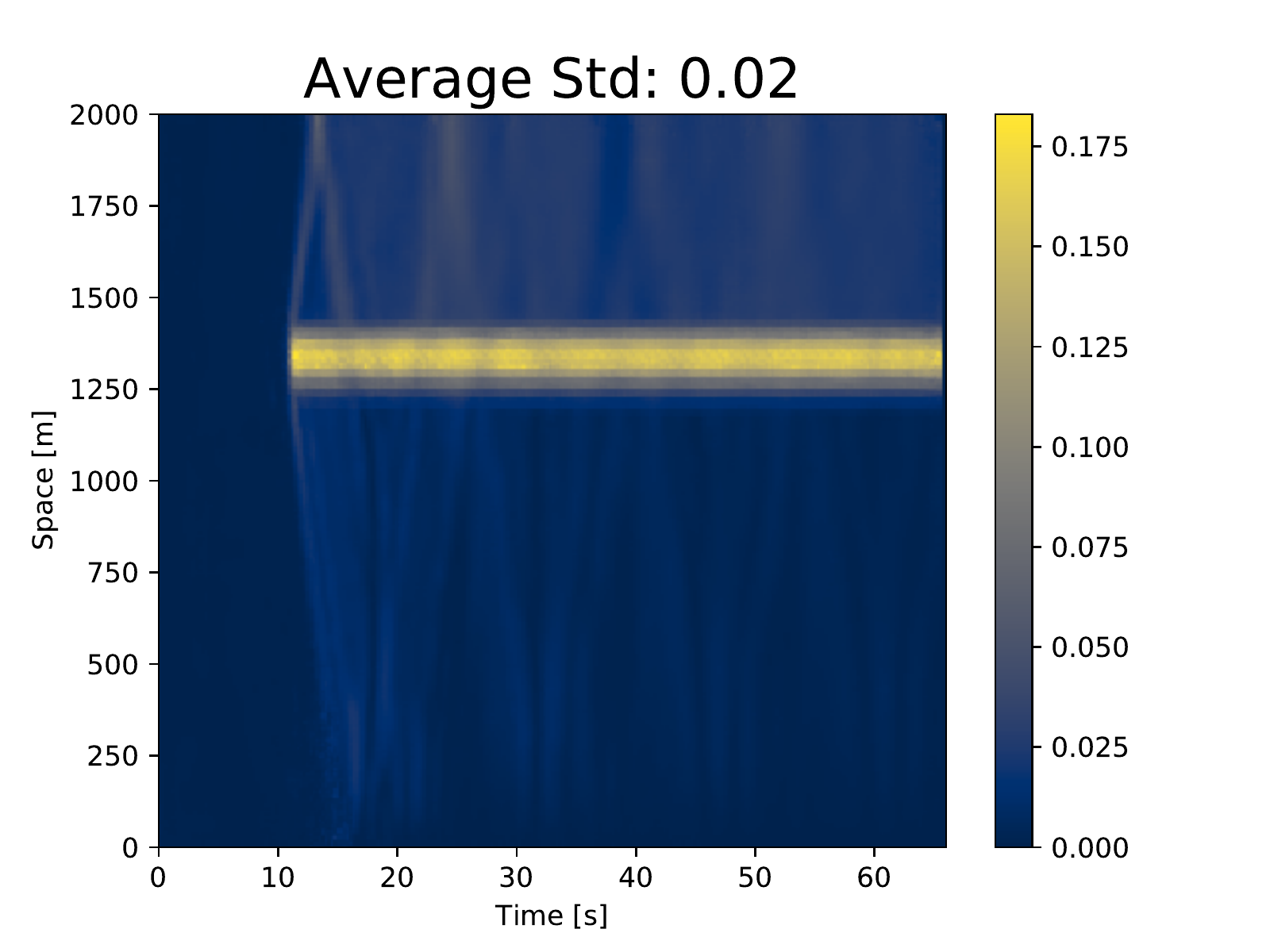}
         \caption{MCGAN velocity standard deviation.}
         \label{fig:pipe_MCGAN_state_std}
     \end{subfigure}
     \begin{subfigure}[b]{0.32\textwidth}
         \centering
         \includegraphics[width=0.95\textwidth]{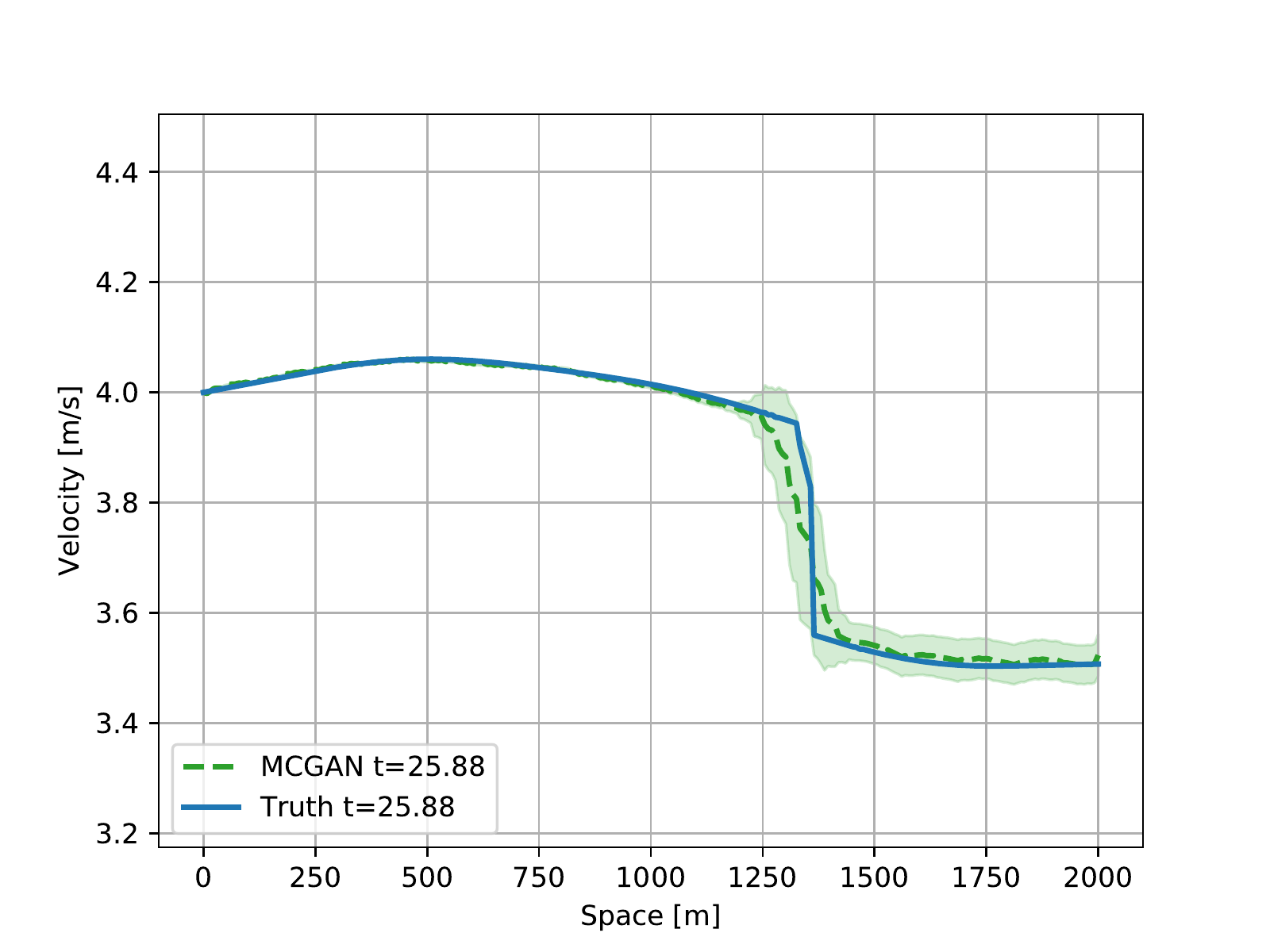}
         \caption{Time 25.88 sec.}
         \label{fig:pipe_reconstruction_at_t_0}
     \end{subfigure}
     \hfill
     \begin{subfigure}[b]{0.32\textwidth}
         \centering
         \includegraphics[width=0.95\textwidth]{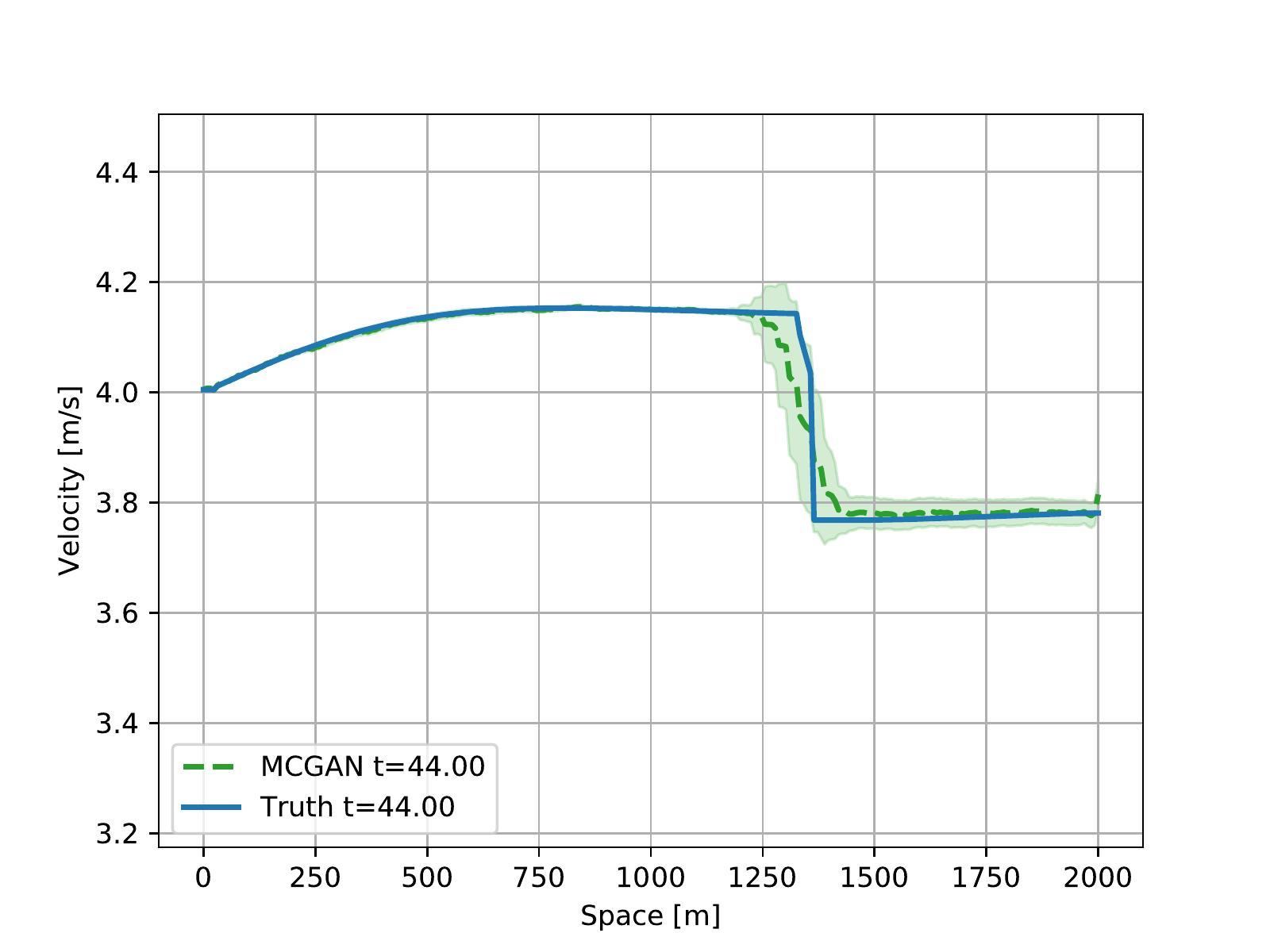}
         \caption{Time 44.00 sec.}
         \label{fig:pipe_reconstruction_at_t_05}
     \end{subfigure}
     \hfill
     \begin{subfigure}[b]{0.32\textwidth}
         \centering
         \includegraphics[width=0.95\textwidth]{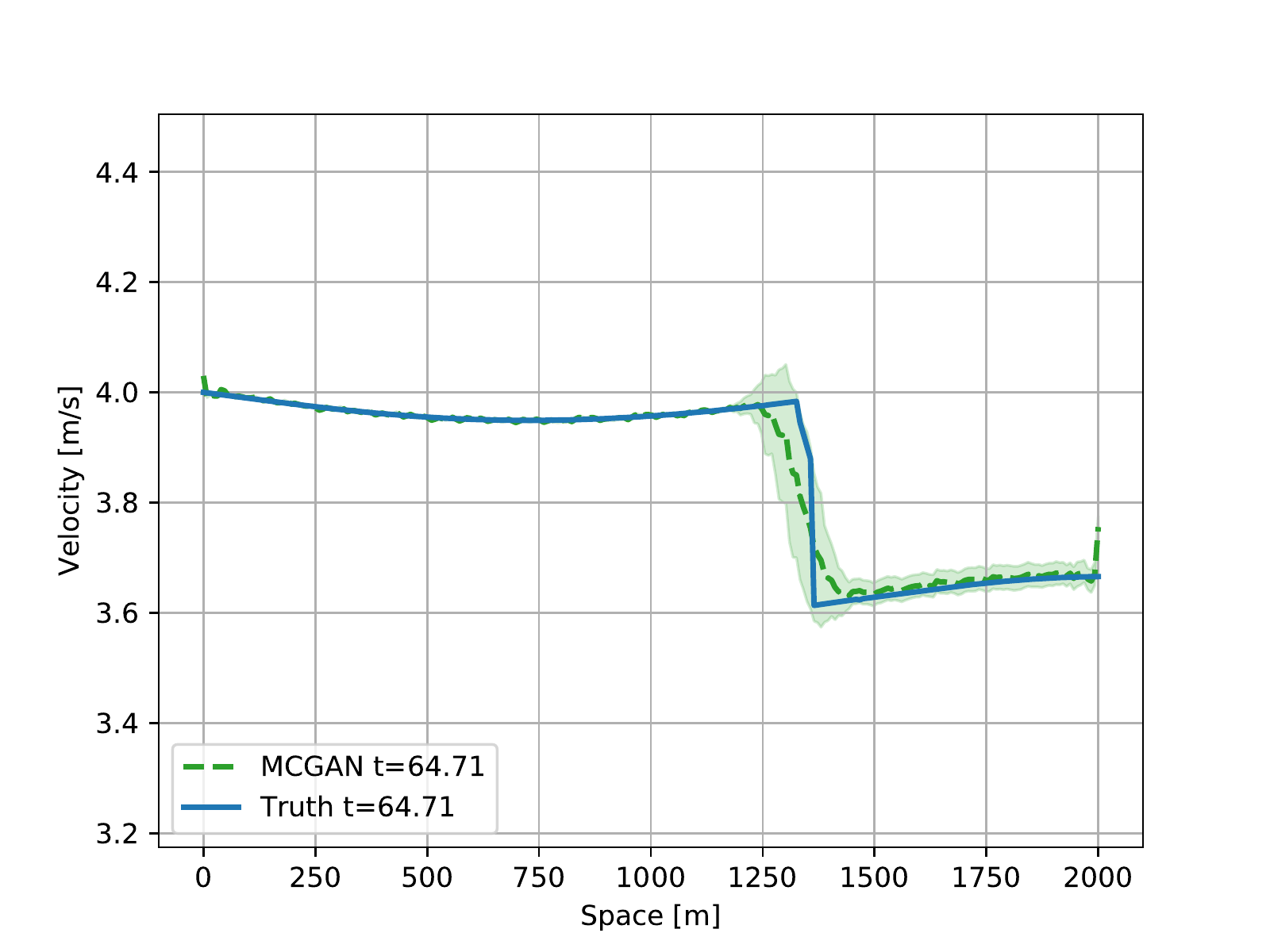}
         \caption{Time 64.71 sec.}
         \label{fig:pipe_reconstruction_at_t_1}
     \end{subfigure}
        \caption{Results for the MCGAN method applied to the pipe flow with a leakage, Eq.\ \eqref{adv_diff_equation}. (a)-(c) are space-time contour plots of the true state, the MCGAN estimated state, and the standard deviation, respectively. (d)-(f) show the state reconstruction at various instances in time with the shaded area denoting one standard deviation away from the reconstruction.}
        \label{fig:pipe_results}
\end{figure}

Figure \ref{fig:pipe_results} presents the reconstruction of the velocity. It is apparent that the velocity is reconstructed very well with a relative RMSE of 0.01. It is especially worth noting that the uncertainty is largest around the drop in velocity, i.e.\ at the location of the leakage, as expected. This uncertainty information could further be used to estimate the location of the leakage. While the state estimation is accurate, it is apparent that the velocity estimation is slightly worse in the domain to the right of the leakage ($x>x_l$). The lack of accuracy is accompanied by an increased standard deviation in that part of the domain. Hence, the uncertainty estimates provide useful information. 

While the MCGAN is performing well in the interior of the domain, it is noteworthy that the estimation at the boundary at $x=2000$ is not as good. 

\begin{figure}
     \centering
     \begin{subfigure}[b]{0.48\textwidth}
         \centering
         \includegraphics[width=0.8\textwidth]{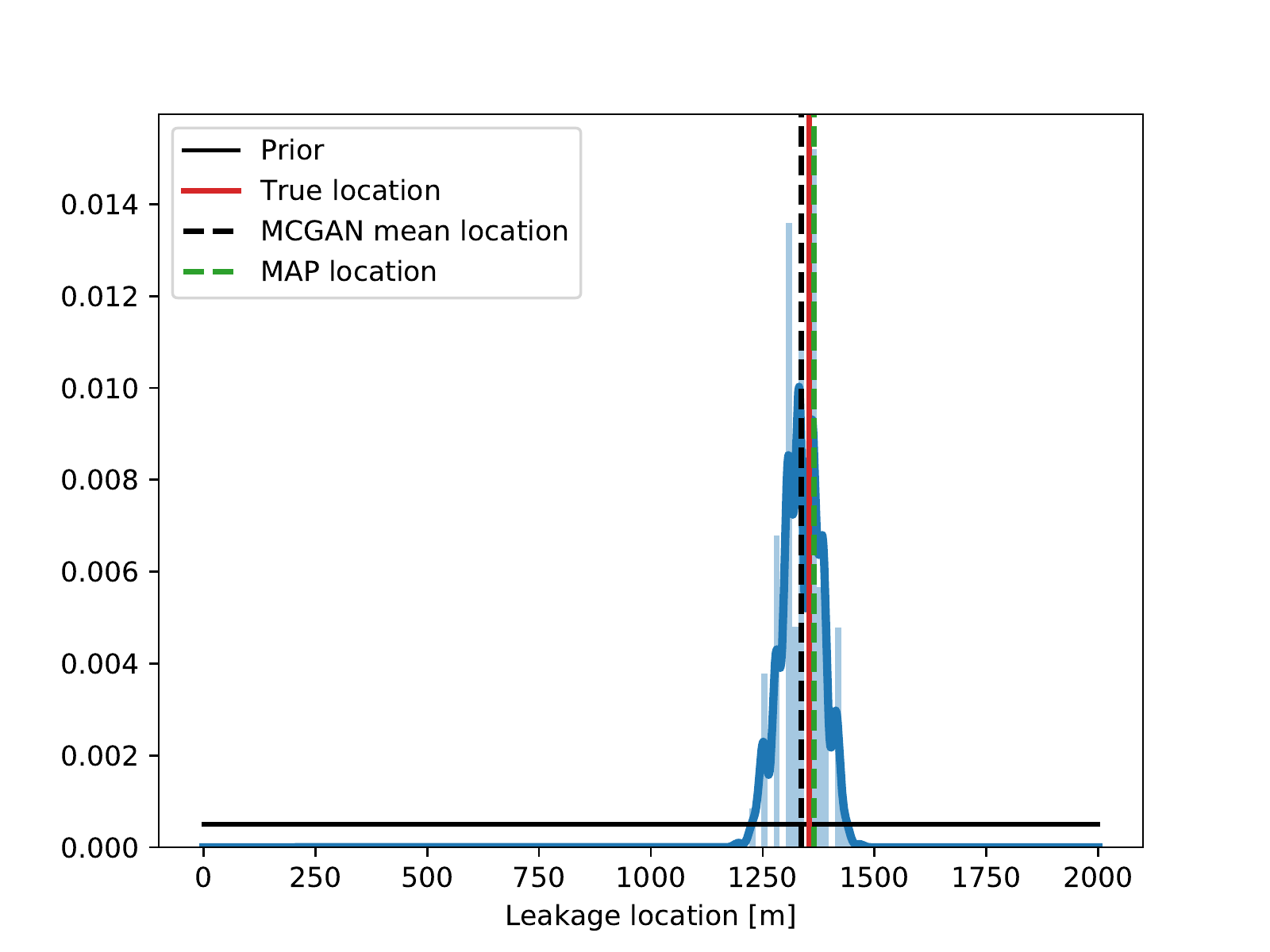}
         \caption{Leakage location.}
         \label{fig:pipe_location_histogram}
     \end{subfigure}
     \hfill
     \begin{subfigure}[b]{0.48\textwidth}
         \centering
         \includegraphics[width=0.8\textwidth]{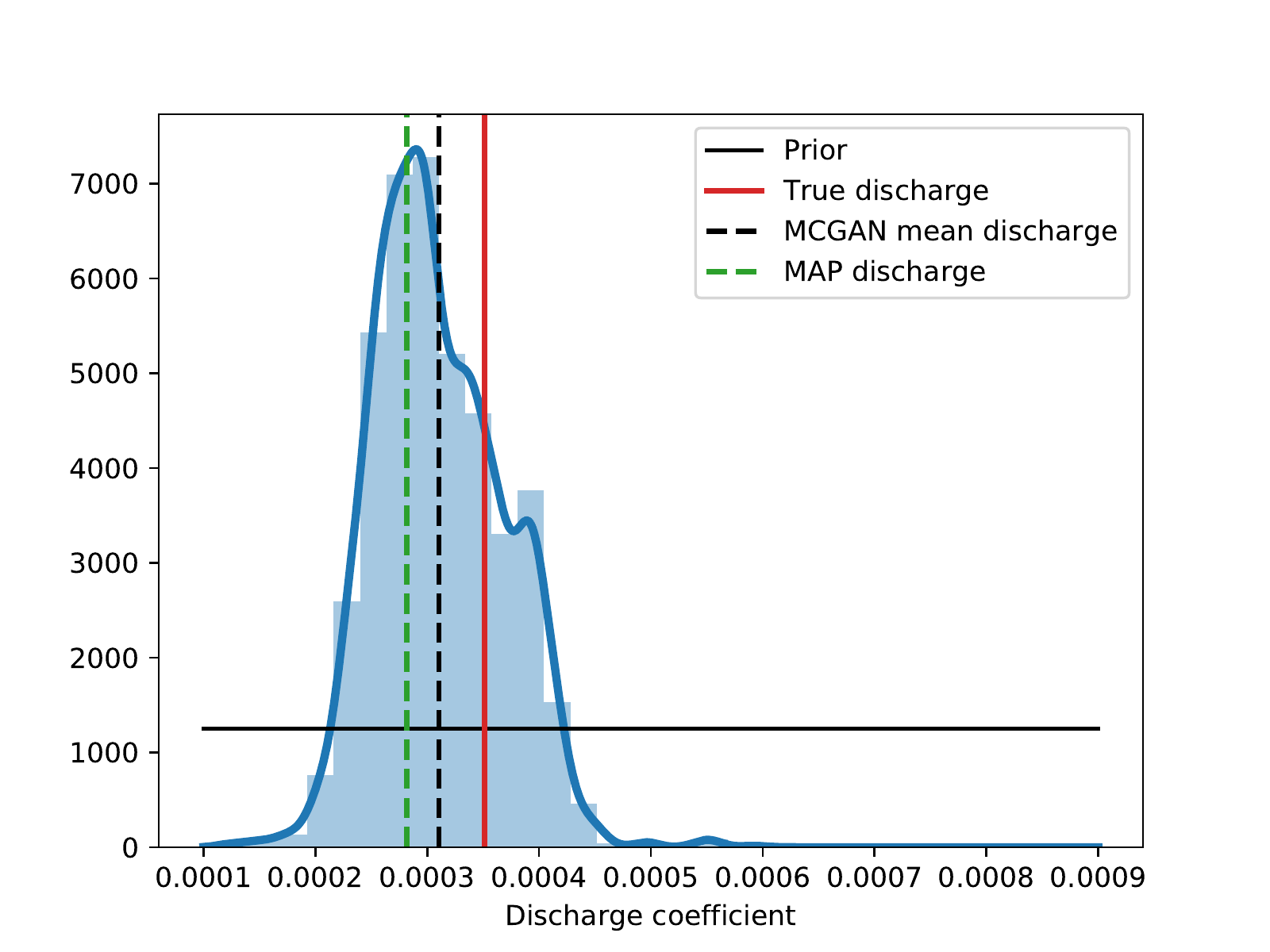}
         \caption{Discharge coefficient.}
         \label{fig:pipe_size_histogram}
     \end{subfigure}
    
        \caption{Posterior distributions of the leakage location and discharge coefficient in the pipe flow equation.}
        \label{fig:pipe_histograms}
\end{figure}

In Figure \ref{fig:pipe_histograms} we see the estimated posterior distributions of the leakage location and discharge coefficient, respectively. In the leakage location posterior, the estimated mean is close to the true mean (see also Table \ref{tab:pipe_flow_estimates}) and it is, more or less, symmetric. In the discharge coefficient posterior, on the other hand, the estimated value appears to be smaller than the true value. 
The MCGAN method significantly outperforms the PCE and EnKF methods in this case. The EnKF is initiated with $x_l=1000$ and $C_d=5\cdot 10^{-4}$ and the MCMC with PCE is initiated at the MAP estimate. In \ref{appendix:alternative_pipe_results}, the results obtained from using the EnKF and PCE method are shown. None of the two approaches manages to estimate the state nor the parameters in a satisfying manner. Table \ref{tab:pipe_flow_estimates} shows that the EnKF approach is unable to update the posterior and the PCE approach only performs marginally better. 

The pipe flow equations are highly nonlinear and the solution exhibits a discontinuity at the location of the leakage. Both phenomena are not easy to handle with the PCE nor EnKF approaches, while neural networks have been shown to be well-suited for such tasks. 

\begin{table}[ht]
\caption{Estimated parameters for the pipe flow using MCGAN, PCE, and EnKF. The best estimates are highlighted in boldface.} \label{tab:pipe_flow_estimates}
\centering
\begin{tabular}{@{}lccccccc@{}}
\toprule
& &\multicolumn{2}{c}{MCGAN}
&\multicolumn{2}{c}{PCE}
&\multicolumn{2}{c}{EnKF} \\
\cmidrule(lr){3-4} \cmidrule(lr){5-6} \cmidrule(lr){7-8} 
Parameter & True value & Mean & Std & Mean & Std & Mean & Std\\
\midrule
Leakage loc. & 1354.45 & \textbf{1336.79} & 44.80 & 1099.19 & 69.98 & 1005.03 & 10.15 \\
Discharge coef. & $3.52\cdot 10^{-4}$  &$\mathbf{3.10}\cdot \mathbf{10^{-4}}$ & $5.70\cdot 10^{-5}$ &$2.64\cdot 10^{-4}$ & $2.36\cdot 10^{-5}$ & $7.23\cdot 10^{-4}$ & $3.15\cdot 10^{-3}$\\
\bottomrule
\end{tabular}
\end{table}

\subsection{Summary of Results}
To summarize the results obtained using our proposed MCGAN method, we highlight accuracy and computation time. Firstly, in Table \ref{tab:rrmse_summary} the relative RMSE for the state and parameters are presented for both test cases. The MCGAN performs better than the two alternative approaches. Especially, in the leakage localization in the pipe flow test case, the MCGAN method outperforms the PCE and the EnKF approaches. The MCGAN results are very close to the true values with relative RMSEs that are one order of magnitudes better than the alternatives for the state and parameter estimation. In Table \ref{tab:computation_time}, the online computation times for the methods applied to the two test cases are shown. It is interesting to note that the computation time does not change much in the two test cases for the MCGAN. This is due to the fact that there are only minor differences in computation time between evaluating a small neural network and a large one. For the Darcy flow, the MCGAN runs an order of magnitude faster than the EKI. The MCGAN method excels for the pipe flow case, where it is two orders of magnitude faster than the alternative approaches. This is because the numerical PDE solution is rather costly, while a forward pass using the generator is cheap.

Lastly, we briefly comment on the offline training time. For the computationally most expensive case, the pipe flow, the most time consuming part is the generation of data. Generating 100,000 training trajectories took about 80 hours on 30 CPU cores (90 seconds per trajectory). The training of the GAN was finished in about 24 hours. In total, the offline stage took approximately 104 hours. The offline time for the Darcy flow was shorter, totaling around 50 hours. We did not experience high sensitivity to the hyperparameters, such as learning rate, batch size, etc. This might be a product of the large number of training samples.

\begin{table}[ht]
\caption{Relative RMSE for the state and parameter estimation for the various test cases. For the Darcy flow, the Relative RMSE for $(v_1, v_2, p)$ is computed. For the pipe flow, the Relative RMSE for $(u,p)$ is computed. The best performing cases are highlighted in boldface.} \label{tab:rrmse_summary}
\centering
\begin{tabular}{@{}lcccccc@{}}
\toprule
&\multicolumn{2}{c}{MCGAN}
&\multicolumn{2}{c}{PCE}
&\multicolumn{2}{c}{EnKF/EKI} \\
\cmidrule(lr){2-3} \cmidrule(lr){4-5} \cmidrule(lr){6-7} 
& State & Parameters & State & Parameters & State & Parameters\\
\midrule
Darcy flow & \textbf{0.1467} & \textbf{0.2619} & - & - & 0.2015 & 0.4240 \\ 
Pipe flow & \textbf{0.0048}  & \textbf{0.013} & 0.0128  & 0.1885 & 0.0327 & 0.2579 \\ 
\bottomrule
\end{tabular}
\end{table}

\begin{table}[ht]
\caption{Comparison of online computation time. The best performing cases are highlighted in boldface.} \label{tab:computation_time}
\centering
\begin{tabular}{@{}lccc@{}}
\toprule
& MCGAN (CPU 1 core) & PCE (CPU 20 cores) & EnKF/EKI (CPU 20 cores)\\
\midrule
Darcy flow & $\mathbf{3.165\cdot 10^{2}}$\textbf{s} & - &$4.111 \cdot 10^{3}$s \\ 
Pipe flow & $\mathbf{7.097 \cdot 10^{2}}$\textbf{s} & $1.1788 \cdot 10^{4}$s & $1.2513 \cdot 10^{4}$s\\
\bottomrule
\end{tabular}
\end{table}

%% file: Conclusion.tex
\section{Conclusion} \label{Conclusion}
We have presented a new method, named MCGAN, to efficiently and accurately solve Bayesian inverse problems in physics and engineering applications. The method combines Generative Adversarial Networks and Markov Chain Monte Carlo methods to sample from posterior distributions by utilizing a low-dimensional latent space and a push-forward map defined as a neural network.

The methodology is divided into two distinct stages, an offline stage, in which the GAN is trained on simulated training data in order to learn the prior distribution, and an online stage, in which the inverse problem is solved for a new set of observations. While the offline stage potentially takes significant computational time, the online stage is computationally very fast and efficient.

We presented a proof of theoretical convergence of the posterior distribution in the Wasserstein-1 distance, in the case where the GAN would be perfectly trained. Furthermore, we provided the insight that sampling from the latent space yields essentially the same results as sampling from the high-dimensional space, in a weak sense. 

To showcase the method's performance, we applied it to two computational engineering test cases with different characteristics and compared it to two alternative approaches. In the high-dimensional problem, the Darcy flow with uncertain permeability field, an improved accuracy was found with MCGAN, as well as a speed-up of one order of magnitude compared to the EKI method. In the second test case, the leakage localization for flow in a pipe, the MCGAN approach was the only method that was able to solve the problem accurately and within a reasonably short timeframe.

While the MCGAN approach performed well on the two test cases, there is still room for future research. We believe the offline stage can be improved by identifying optimal ways of simulating training data and determining hyperparameters for the GAN. This includes determining the optimal size of the latent space. Furthermore, the GAN can further be improved by incorporating physics knowledge either in the training or directly in the neural network architecture. This could possibly alleviate the boundary estimation problems. Lastly, possibilities of using the MCGAN framework in a sequential fashion, as is the case for Kalman filters, could be an interesting direction to explore. 

In conclusion, we believe that the MCGAN methodology can form an important piece of the puzzle towards a well performing digital twin framework, in which real-time state and parameter estimation is of crucial importance.

%% file: Appendix.tex
\section{Proof of Theorem \ref{theorem:posterior_u_to_z}} \label{appendix_posterior_u_to_z_proof}

\begin{proof}
Firstly, since neural networks with continuous activation functions are continuous, they are also measurable \cite{bogachev2007measure}. Therefore, the generator defines a push forward distribution and Eq. \eqref{expected_change_of_variables} is applicable.

Secondly, we look at the evidence. Assuming the likelihood is measurable with respect to $\mathbf{u}$, we have from Eq. \eqref{expected_change_of_variables}:
\begin{align} \label{evidence_u_to_z}
    Q_u(\mathbf{y}) = \int_{\mathbb{R}^{N_u}} \rho_{y|u}^g(\mathbf{y} | \mathbf{u}) \rho_0^g(\mathbf{u}) \: \text{d}\mathbf{u} 
    = \int_{\mathbb{R}^{N_z}} \rho_{y|u}^g(\mathbf{y} | G_\theta(\mathbf{z})) \rho_z^g(\mathbf{z}) \: \text{d}\mathbf{z} = Q_z(\mathbf{y}).
\end{align}
Consider the expected value of the likelihood times some measurable function, $f$, with respect to the prior:
\begin{align} \label{likelihood_and_prior}
\begin{split}
    \mathbb{E}_{U\sim P_{0}^g}[f(U) \rho_{y|u}^g(\mathbf{y}|U)] &= \int_{E} \underbrace{f(\mathbf{u}) \rho_{y|u}^g(\mathbf{y}|\mathbf{u})}_{=\xi(\mathbf{u})} \rho_0^g(\mathbf{u}) \: \mathrm{d} \mathbf{u} \\
    &= \int_{G_\theta^{-1}(E)} \underbrace{f(G_\theta(\mathbf{z})) \rho_{y|u}^g(\mathbf{y}|G_\theta(\mathbf{z}))}_{=\xi(G_\theta(\mathbf{z}))} \rho_z^g(\mathbf{z}) \mathrm{d} \mathbf{z} \\
    &=\mathbb{E}_{U\sim P_{z}^g}[f(G_\theta(Z)) \rho_{y|u}^g(\mathbf{y}|G_\theta(Z))].
\end{split}
\end{align}
Note that $\xi$ is the product of two measurable functions and is therefore measurable. Hence, Eq.\ \eqref{expected_change_of_variables} applies. Now using Eq.\ \eqref{evidence_u_to_z} and \eqref{likelihood_and_prior} we get:
\begin{align*}
    \mathbb{E}_{U\sim P_{u|y}^g}[f(U)] &= \frac{1}{Q_u(\mathbf{y})}\mathbb{E}_{U\sim P_{0}^g}[f(U) \rho_{y|u}^g(\mathbf{y}|U)]  \\
    &= \frac{1}{Q_z(\mathbf{y})}\mathbb{E}_{Z\sim P_{z}^g}[f(G_\theta(Z) \rho_{y|u}^g(\mathbf{y}|G_\theta(Z))] \\
    &= \mathbb{E}_{Z\sim P_{z|y}^g}[f(G(Z))].
\end{align*}
\end{proof}

\section{Proof of Theorem \ref{main_theorem}} \label{appendix_main_theorem_proof}

\begin{proof}
We write the Wasserstein-1 distance between the real prior and the generated prior in dual form:
\begin{align}
    W_1(P_0^r, P_0^r) = \sup_{\mathrm{Lip}(f)\leq 1} \left| \int_{E} f(\mathbf{u}) \rho_0^r(\mathbf{u}) \text{d}\mathbf{u} - \int_{E} f(\mathbf{u}) \rho_0^g(\mathbf{u}) \text{d}\mathbf{u}  \right|,
\end{align}
where $f:E\rightarrow \mathbb{R}$ is Lipschitz continuous with Lipschitz constant less or equal $1$ and $f(\mathbf{u}_0)=0$ for some $\mathbf{u}_0$. Note that any function, $g$, with Lipschitz constant less than or equal $1$, is a contraction and therefore admits a fixed point. Now, assuming that $\mathbf{u}_0$ is the fixed point, we can simply define $f=g-\mathbf{u}_0$, which admits $f(\mathbf{u}_0)=\mathbf{u}_0$. Therefore, assuming $f(\mathbf{u}_0)=0$ for some $\mathbf{u}_0$ is not a restriction. Furthermore, we have:
\begin{align*}
    |f(\mathbf{u})| = |f(\mathbf{u})+f(\mathbf{u}_0)-f(\mathbf{u}_0)| = |f(\mathbf{u})+f(\mathbf{u}_0)| \leq \mathrm{Lip}(f) d(\mathbf{u},\mathbf{u}_0) \leq D.
\end{align*}
The Wasserstein-1 distance between the posteriors is given by:
\begin{align*}
    W_1(P_{u|y}^r, P_{u|y}^g) &= \sup_{\mathrm{Lip}(f)\leq 1} \left| \int_{E} f(\mathbf{u}) \rho_{u|y}^r(\mathbf{u}|\mathbf{y}) \text{d}\mathbf{u} - \int_{E} f(\mathbf{u}) \rho_{u|y}^g(\mathbf{u}|\mathbf{y}) \text{d}\mathbf{u}  \right| \\
    &= \sup_{\mathrm{Lip}(f)\leq 1} \left| \int_{E} f(\mathbf{u}) (\rho_{u|y}^r(\mathbf{u}|\mathbf{y}) -\rho_{u|y}^g(\mathbf{u}|\mathbf{y})) \text{d}\mathbf{u}  \right| \\
    &= \sup_{\mathrm{Lip}(f)\leq 1} \left| \int_{E} f(\mathbf{u}) \left( \frac{\Phi^r(\mathbf{u}) \rho_0^r(\mathbf{u})}{Q^r_u(\mathbf{y})}-\frac{\Phi^g(\mathbf{u}) \rho_0^g(\mathbf{u})}{Q^g_u(\mathbf{y})} \right) \text{d}\mathbf{u}  \right|,
\end{align*}
Adding and subtracting the term $f(\mathbf{u})\frac{\Phi^r(\mathbf{u})\rho_0^g(\mathbf{u})}{Q^r_u(\mathbf{y})}$ gives
\begin{align*}
    W_1(P_{u|y}^r, P_{u|y}^g) &= 
    \sup_{\mathrm{Lip}(f)\leq 1} \left| \int_{E} f(\mathbf{u}) \left( \frac{\Phi^r(\mathbf{u}) \rho_0^r(\mathbf{u})}{Q^r_u(\mathbf{y})}-\frac{\Phi^g(\mathbf{u}) \rho_0^g(\mathbf{u})}{Q^g_u(\mathbf{y})} + \frac{\Phi^r(\mathbf{u})\rho_0^g(\mathbf{u})}{Q^r_u(\mathbf{y})} - \frac{\Phi^r(\mathbf{u})\rho_0^g(\mathbf{u})}{Q^r_u(\mathbf{y})} \right) \text{d}\mathbf{u}  \right| \\
    &\leq \sup_{\mathrm{Lip}(f)\leq 1} \left| \int_{E} f(\mathbf{u}) \frac{\Phi^r(\mathbf{u})}{Q^r_u(\mathbf{y})} \left( \rho_0^r(\mathbf{u})-\rho_0^g(\mathbf{u}) \right)  \text{d}\mathbf{u}  \right|
    + \left| \int_{E}f(\mathbf{u})\rho_0^g(\mathbf{u}) \left( \frac{\Phi^r(\mathbf{u})}{Q^r_u(\mathbf{y})} - \frac{\Phi^g(\mathbf{u})}{Q^g_u(\mathbf{y})} \right) \text{d}\mathbf{u}  \right|.
\end{align*}
Subsequently, adding and subtracting the term $f(\mathbf{u})\frac{\Phi^r(\mathbf{u})\rho_0^g(\mathbf{u})}{Q^g_u(\mathbf{y})}$ in the second integral gives
\begin{align*}
    W_1(P_{u|y}^r, P_{u|y}^g) \leq&  \sup_{\mathrm{Lip}(f)\leq 1} \left| \int_{E} f(\mathbf{u}) \frac{\Phi^r(\mathbf{u})}{Q^r_u(\mathbf{y})} \left( \rho_0^r(\mathbf{u})-\rho_0^g(\mathbf{u}) \right)  \text{d}\mathbf{u}  \right| \\
    &+ \left| \int_{E}f(\mathbf{u})\rho_0^g(\mathbf{u}) \left( \frac{\Phi^r(\mathbf{u})}{Q^r_u(\mathbf{y})} - \frac{\Phi^g(\mathbf{u})}{Q^g_u(\mathbf{y})} \right) 
    + f(\mathbf{u})\rho_0^g(\mathbf{u})\left( \frac{\Phi^r(\mathbf{u})}{Q^g_u(\mathbf{y})} -\frac{\Phi^r(\mathbf{u})}{Q^g_u(\mathbf{y})}    \right)\text{d}\mathbf{u}  \right| \\
    \leq& \sup_{\mathrm{Lip}(f)\leq 1} \underbrace{\left| \int_{E} f(\mathbf{u}) \frac{\Phi^r(\mathbf{u})}{Q^r_u(\mathbf{y})} \left( \rho_0^r(\mathbf{u})-\rho_0^g(\mathbf{u}) \right)  \text{d}\mathbf{u}  \right|}_{=I_1} \\
    &+ \underbrace{\left|\int_{E}\left( \frac{1}{Q^r_u(\mathbf{y})} - \frac{1}{Q^g_u(\mathbf{y})}\right)  f(\mathbf{u}) \rho_0^g(\mathbf{u}) \Phi^r(\mathbf{u})   \text{d}\mathbf{u}  \right|}_{=I_2} \\
     &+ \underbrace{\left| \int_{E}\frac{1}{Q^g_u(\mathbf{y})}  f(\mathbf{u}) \rho_0^g(\mathbf{u}) \left( \Phi^r(\mathbf{u})-\Phi^g(\mathbf{u})  \right) \text{d}\mathbf{u}  \right|}_{=I_3}.
\end{align*}
We will consider $I_1$, $I_2$, and $I_3$ individually. Starting with $I_3$, we use that
\begin{align} \label{exp_inequality}
    |e^{-x_1}-e^{-x_2}| \leq e^{-\min(x_1,x_2)}|x_1-x_2| \Rightarrow |\Phi^r(\mathbf{u}) -\Phi^g(\mathbf{u}) | \leq \max_\mathbf{u}(\Phi^r,\Phi^g)|l^r(\mathbf{u})-l^g(\mathbf{u})|.
\end{align}
Using Eq. \eqref{exp_inequality} and the fact that $|f(\mathbf{u})| \leq d(\mathbf{u},\mathbf{u}_0)$, together with the Cauchy-Schwartz inequality, we get:
\begin{align*}
    \sup_{\mathrm{Lip}(f)\leq 1} I_3 &\leq
   \sup_{\mathrm{Lip}(f)\leq 1}\frac{\max(\Phi^r,\Phi^g)}{Q^g_u(\mathbf{y})}   \left| \int_{E} f(\mathbf{u}) \rho_0^g(\mathbf{u}) |l^r(\mathbf{u})-l^g(\mathbf{u})| \text{d}\mathbf{u}  \right| \\
   &\leq
   \frac{\max(\Phi^r,\Phi^g)}{Q^g_u(\mathbf{y})}   \left| \int_{E} d(\mathbf{u},\mathbf{u}_0) \rho_0^g(\mathbf{u}) |l^r(\mathbf{u})-l^g(\mathbf{u})| \text{d}\mathbf{u}  \right| \\
     &\leq\frac{\max(\Phi^r,\Phi^g)}{Q^g_u(\mathbf{y})}   \left( \int_{E} d(\mathbf{u},\mathbf{u}_0)^2 \rho_0^g(\mathbf{u}) \text{d}\mathbf{u} \right)^{1/2}  ||l^r(\mathbf{u})-l^g(\mathbf{u}) ||_{L^2_{\rho_0^g}} \\
   &\leq \underbrace{\frac{\max(\Phi^r,\Phi^g)}{Q^g_u(\mathbf{y})} |P_0^g|_{\mathcal{W}_2}}_{=C_3} ||l^r(\mathbf{u})-l^g(\mathbf{u}) ||_{L^2_{\rho_0^g}}.
\end{align*}
Considering $I_2$, we use the following \cite{sprungk2020local}:
\begin{align}
    \left| \frac{1}{Q^r_u(\mathbf{y})} - \frac{1}{Q^g_u(\mathbf{y})}\right| = \frac{|Q^g_u(\mathbf{y})-Q^r_u(\mathbf{y})|}{Q^r_u(\mathbf{y})Q^g_u(\mathbf{y})},
\end{align}
and
\begin{align}
    |Q^g_u(\mathbf{y})-Q^r_u(\mathbf{y})| \leq \max(\Phi^r,\Phi^g) ||l^r(\mathbf{u}) - l^g(\mathbf{u}) ||_{L^1_{\rho_0^g}},
\end{align}
in order to get:
\begin{align*}
    I_2 &= \left|\int_{E} \left(\frac{1}{Q^r_u(\mathbf{y})} - \frac{1}{Q^g_u(\mathbf{y})}\right)  f(\mathbf{u}) \rho_0^g(\mathbf{u}) \Phi^r(\mathbf{u})   \text{d}\mathbf{u}  \right| \\
    &\leq \left|\frac{1}{Q^r_u(\mathbf{y})} - \frac{1}{Q^g_u(\mathbf{y})}\right| \left|\int_{E}  f(\mathbf{u}) \rho_0^g(\mathbf{u}) \Phi^r(\mathbf{u})   \text{d}\mathbf{u}  \right| \\
    &\leq \frac{|Q^g_u(\mathbf{y})-Q^r_u(\mathbf{y})|}{Q^r_u(\mathbf{y})Q^g_u(\mathbf{y})}\left|\int_{E}  f(\mathbf{u}) \rho_0^g(\mathbf{u}) \Phi^r(\mathbf{u})   \text{d}\mathbf{u}  \right| \\
    &\leq \frac{\max(\Phi^r,\Phi^g)}{Q^r_u(\mathbf{y})Q^g_u(\mathbf{y})} \left|\int_{E}  f(\mathbf{u}) \rho_0^g(\mathbf{u}) \Phi^r(\mathbf{u})   \text{d}\mathbf{u}  \right| ||l^r(\mathbf{u}) - l^g(\mathbf{u}) ||_{L^1_{\rho_0^g}}
\end{align*}
By defining the function, $g:E\rightarrow\mathbb{R}$, $g(\mathbf{u})= f(\mathbf{u})\Phi^r(\mathbf{u})$, one can show that $g$ is Lipschitz continuous with Lipschitz constant $\mathrm{Lip}(g)=1+D\mathrm{Lip}(\Phi^r)$ \cite{sprungk2020local}. Furthermore, we have have $|g(\mathbf{u})| \leq d(\mathbf{u},\mathbf{u}_0)$. This gives:
\begin{align*}
     \sup_{\mathrm{Lip}(f)\leq 1} I_2 &\leq \frac{\max(\Phi^r,\Phi^g)}{Q^r_u(\mathbf{y})Q^g_u(\mathbf{y})} (1+D\mathrm{Lip}(\Phi^r)) \left|\int_{E}  d_E(\mathbf{u},\mathbf{u}_2)\rho_0^g(\mathbf{u}) \text{d}\mathbf{u}  \right| ||l^r(\mathbf{u}) - l^g(\mathbf{u}) ||_{L^1_{\rho_0^g}}\\
    &\leq \underbrace{\frac{\max(\Phi^r,\Phi^g)}{Q^r_u(\mathbf{y})Q^g_u(\mathbf{y})} (1+D\mathrm{Lip}(\Phi^r)) |P_0^g|_{\mathcal{W}_1}}_{=C_2} ||l^r(\mathbf{u}) - l^g(\mathbf{u}) ||_{L^1_{\rho_0^g}}.
\end{align*}
Finally, we consider $I_1$. By using the function, $g:E\rightarrow\mathbb{R}$, $\mathbf{u}\mapsto f(\mathbf{u})\Phi^r(\mathbf{u})$, as defined above, we get:
\begin{align*}
    \sup_{\mathrm{Lip}(g)\leq 1} I_1 &\leq  \sup_{\mathrm{Lip}(f)\leq 1} \frac{(1+D\mathrm{Lip}(\Phi^r))}{Q^r_u(\mathbf{y})}  \left| \int_{E} g(\mathbf{u})  \left( \rho_0^r(\mathbf{u})-\rho_0^g(\mathbf{u}) \right)  \text{d}\mathbf{u}  \right| \\
    &= \underbrace{\frac{(1+D\mathrm{Lip}(\Phi^r))}{Q^r_u(\mathbf{y})}}_{=C_1} W_1(P_0^r, P_0^r)
\end{align*}
Combining, $I_1$, $I_2$, and $I_3$ we then get:
\begin{align*}
    W_1(P_{u|y}^r, P_{u|y}^g) &\leq \sup_{\mathrm{Lip}(f)\leq 1} I_1 + I_2 + I_3 \\
    &= C_1 W_1(P_0^r, P_0^r)  + C_2 ||l^r(\mathbf{u}) - l^g(\mathbf{u}) ||_{L^1_{\rho_0^g}} + C_3 ||l^r(\mathbf{u})-l^g(\mathbf{u}) ||_{L^2_{\rho_0^g}} \\
    &\leq C_1 \epsilon_1  + C_2 \epsilon_2 + C_3 \epsilon_3,
\end{align*}
\end{proof}

\section{Training Wasserstein GANs} \label{training_WGANS}
\begin{figure}
    \centering
    \includegraphics[width=\textwidth]{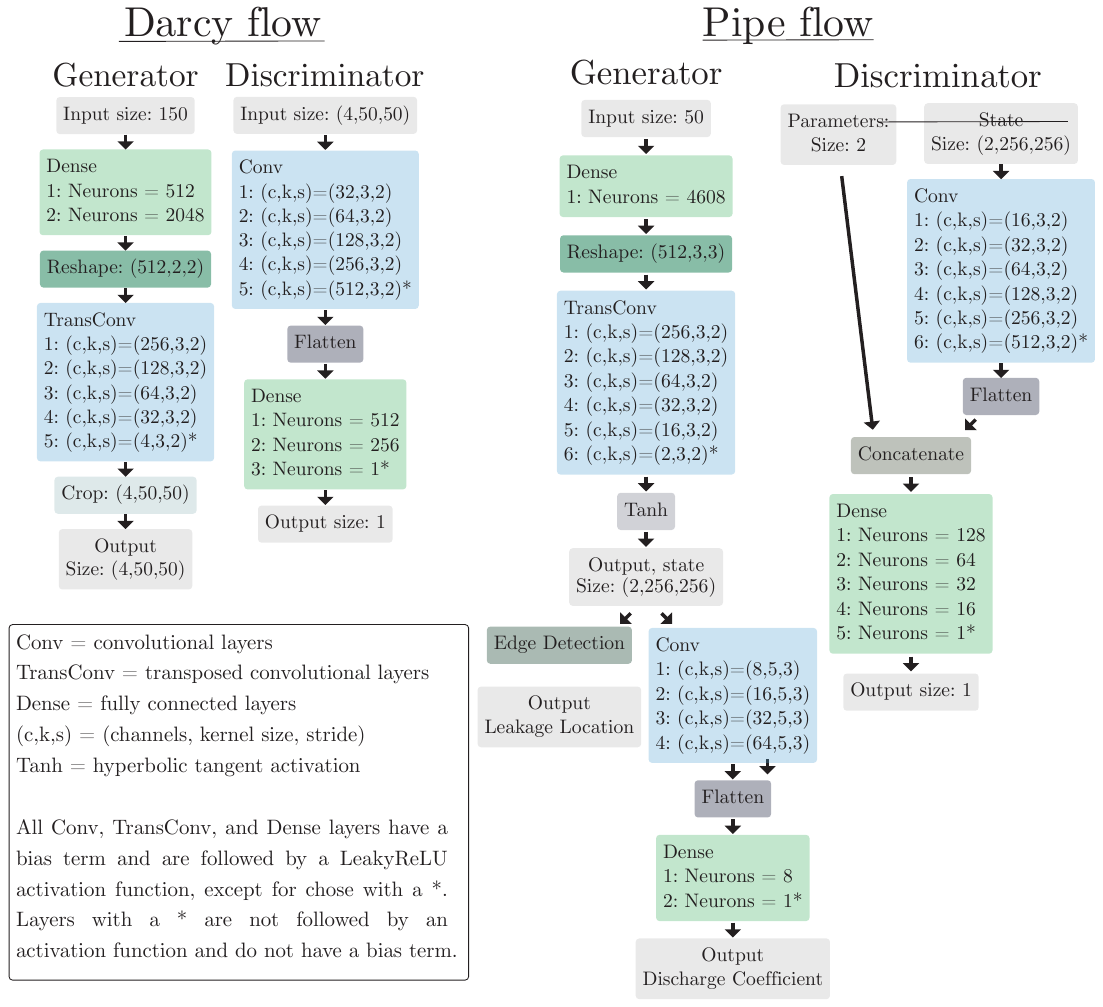}
    \caption{Generator and discriminator architectures for the two test cases.}
    \label{fig:GAN_networks}
\end{figure}

In the WGAN framework, it is important to properly train the discriminator. Therefore, it is common practice to update the discriminator parameters more frequently than the generator parameters. The number of discriminator updates, relative to those of the generator, is denoted by $n_{disc}/n_{gen}$. The hyperparameters for the training of the three test cases are shown in Table \ref{tab:hyperparameters}. The specific architectures used are shown in Figure \ref{fig:GAN_networks}.

We compared the Adam optimizer and the RMSprop optimizer for training, and found that RMSProp, in general, showed superior results in our test cases. 

\begin{table}[h]
\centering
\caption{Hyperparameters for the WGANs for the three test cases.}
\label{tab:hyperparameters}
\begin{tabular}{@{}lccc@{}}
\toprule
Hyperparameters $\backslash$ Test case & Darcy flow & Pipe flow \\
\midrule
Optimizer & RMSProp & RMSProp \\ 
Learning rate  & $10^{-4}$ & $10^{-4}$ \\
Batch size & 64 & 64 \\ 
Gradient penalty  & 5 & 5\\ 
$n_{disc}/n_{gen}$  & 1 & 2\\ 
$N_{train}$ &  300,000 & 100,000\\ 
Latent dimension ($N_z$) & 150 & 50\\ 
\bottomrule
\end{tabular}
\end{table}

\section{Alternative Methods}\label{alternative_methods}

\subsection{Ensemble Kalman Filter} \label{kalman_filter}
We make use of two variations of the ensemble Kalman filter (EnKF):
\begin{itemize}
\item The (standard) EnKF for dynamic problems, where the state and parameter distributions are computed based on previous time steps along with data availability;
\item Ensemble Kalman Inversion (EKI), used for stationary problems, where an artificial time dimension is introduced in order to iteratively update the posterior of the state and parameters.
\end{itemize}
For the pipe flow equations the standard EnKF is utilized while for the Darcy flow the EKI is used. The EnKF implementation is based on \cite{harlim2018data} and the EKI implementation is based on \cite{ding2021ensemble}.

For simultaneously estimating the parameter and state, we make use of disturbance modeling \cite{horsholt2019spatial}. Here, we define an augmented model:
\begin{subequations} \label{EnKF_equations}
\begin{align} 
    \mathbf{q}_i &= F(\mathbf{q}_{i-1},\mathbf{m}_{i-1}) + \Gamma_q \epsilon_i, \quad  \epsilon_i\sim \mathcal{N}(0,Q_q),\\
    \mathbf{m}_i &= \mathbf{m}_{i-1}+ \Gamma_m \delta_i, \quad  \delta_i\sim \mathcal{N}(0,Q_m), \\
    \mathbf{y}_i &= \mathbf{h}(\mathbf{q}_{i}) + \eta_i, \quad \eta_i\sim \mathcal{N}(0,R),
\end{align}
\end{subequations}
where $F$ is the discrete one-step time advancement model, $\epsilon_i$ is the model noise, $Q$ the model covariance, and $R$ the observation covariance. With this formulation, both the state and parameters are updated in every step of the EnKF algorithm.   

\subsection{Polynomial Chaos Expansion} \label{polynomial_chaos_expansion}
The basic idea behind polynomial chaos expansion (PCE) is to create a surrogate model that maps the stochastic parameters, $\mathbf{m}$, to a quantity of interest, $Q$ \cite{xiu2010numerical}. The surrogate model is defined by a linear expansion of orthogonal polynomials:
\begin{align}
    Q(\mathbf{m}) = \sum_{i=1}^{N} \alpha_i \phi_i(\mathbf{m}), 
\end{align}
where $\phi_i$ are the polynomials that are chosen based on the distribution of $\mathbf{m}$, and $\alpha_i$ are the generalized Fourier coefficients.

The coefficients, $\alpha_i$, are typically computed using either spectral projection methods or by least squares minimization. In both cases, the evaluations are carefully chosen according to a quadrature rule. 

In our test cases, we choose the quantity of interest to be the observations, i.e.\ $Q(\mathbf{m}) \approx \mathbf{h}(\mathbf{q}(\mathbf{m}))$ and $\alpha_i\in\mathbb{R}^{N_y}$. $\mathbf{m}$ are the parameters of interest, which are often the model parameters and/or initial and boundary conditions.  

When the PCE is computed, the posterior PDF is defined by:
\begin{align}
    \rho^{y}_m(\mathbf{m}|\mathbf{y}) = \frac{1}{\rho_y(\mathbf{y})}\rho_\eta(\mathbf{y}-Q(\mathbf{m}))\rho_{0}^m(\mathbf{m}).
\end{align}
The expected state and parameters are then computed by:
\begin{align}
    \mathbb{E}_\mathbf{q}\left[ \mathbf{q} \right] \approx 
    \frac{1}{N_{sample}}\sum_{i=1}^{N_{sample}}\mathbf{q}(\mathbf{m}_i), \quad 
    \mathbb{E}_{\mathbf{m}}\left[ \mathbf{m} \right] \approx 
    \frac{1}{N_{sample}}\sum_{i=1}^{N_{sample}} \mathbf{m}_i, \quad \mathbf{m}_i \sim  P_{m}^{y},
\end{align}
and the variance is computed in a similar manner. 

It is important to notice that the sampling is done in the parameter space and the state is thereafter computed by using the sampled parameters as input for the forward problem. Directly sampling the state is infeasible due to the high-dimensionality of the state. 

The implementation of the PCE method is done using the Python library Chaospy \cite{feinberg2015chaospy}.

\subsection{Darcy Flow}
\begin{figure}[h]
     \centering
     \begin{subfigure}[t]{0.24\textwidth}
         \centering
         \includegraphics[width=\textwidth]{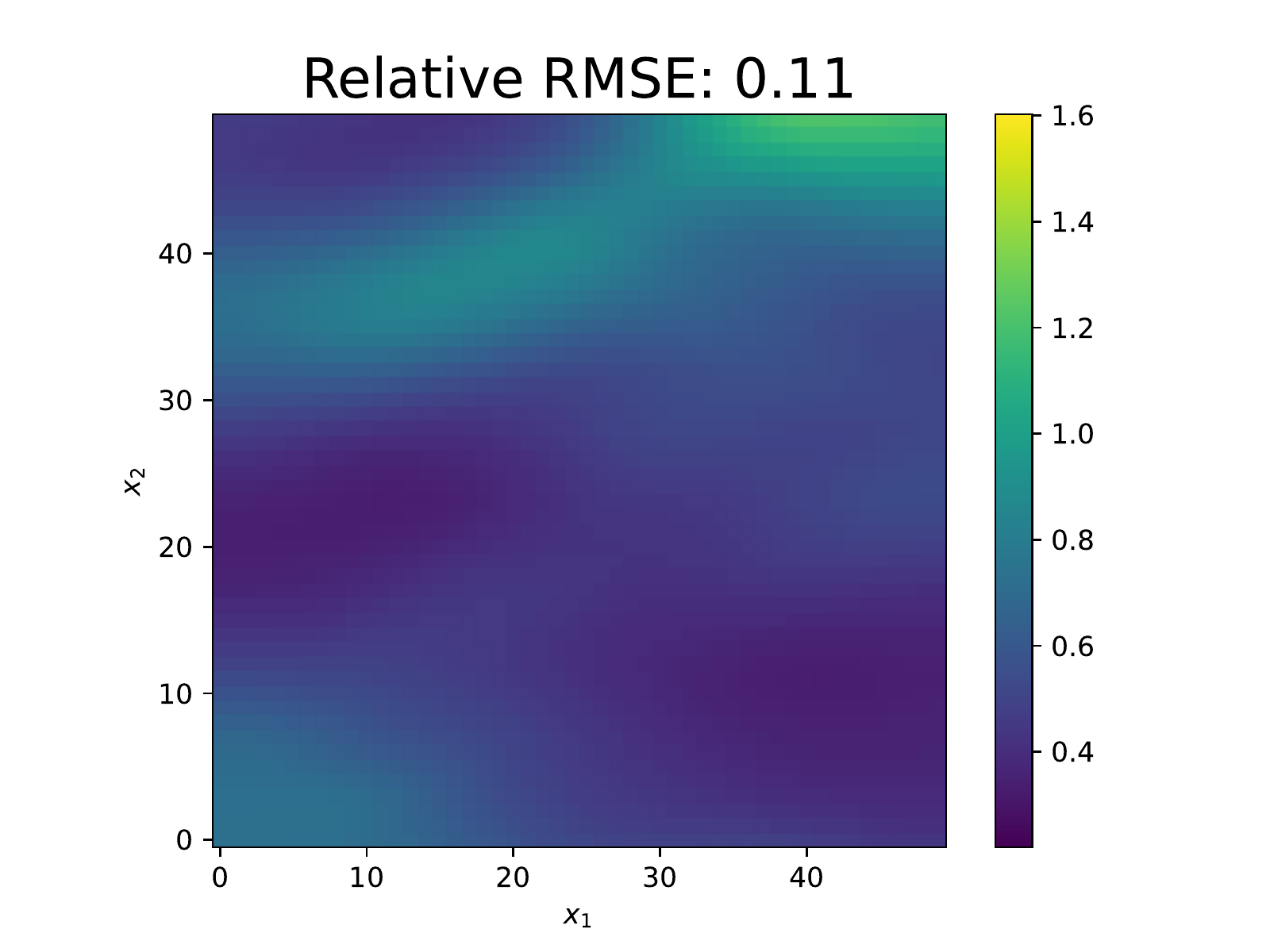}
         \caption{EKI approximated $v_1$.}
         \label{fig:darcy_kalman_state}
     \end{subfigure}
     \begin{subfigure}[t]{0.24\textwidth}
         \centering
         \includegraphics[width=\textwidth]{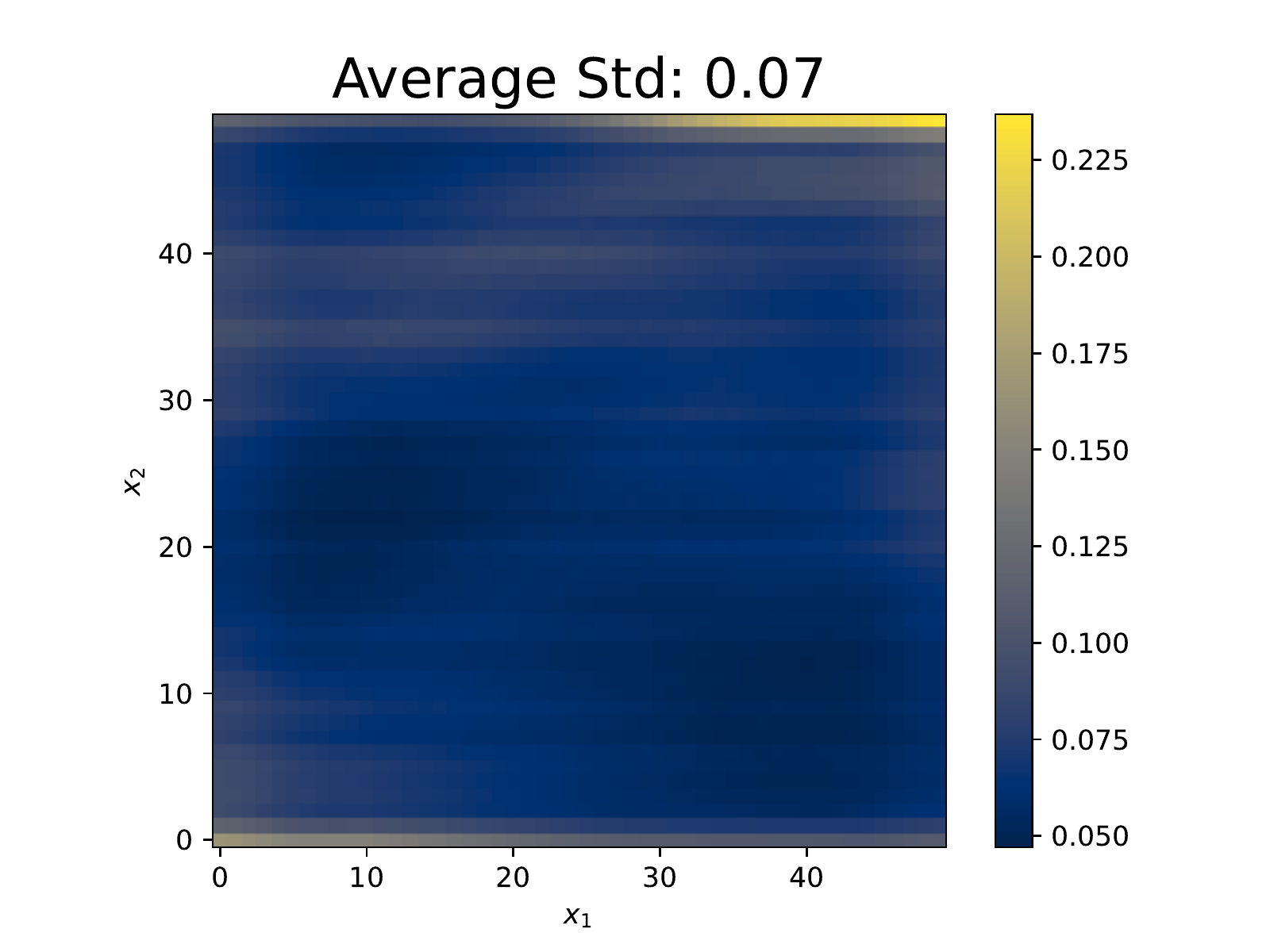}
         \caption{EKI standard deviation $v_1$.}
         \label{fig:darcy_kalman_state_std}
     \end{subfigure}
     \begin{subfigure}[t]{0.24\textwidth}
         \centering
         \includegraphics[width=\textwidth]{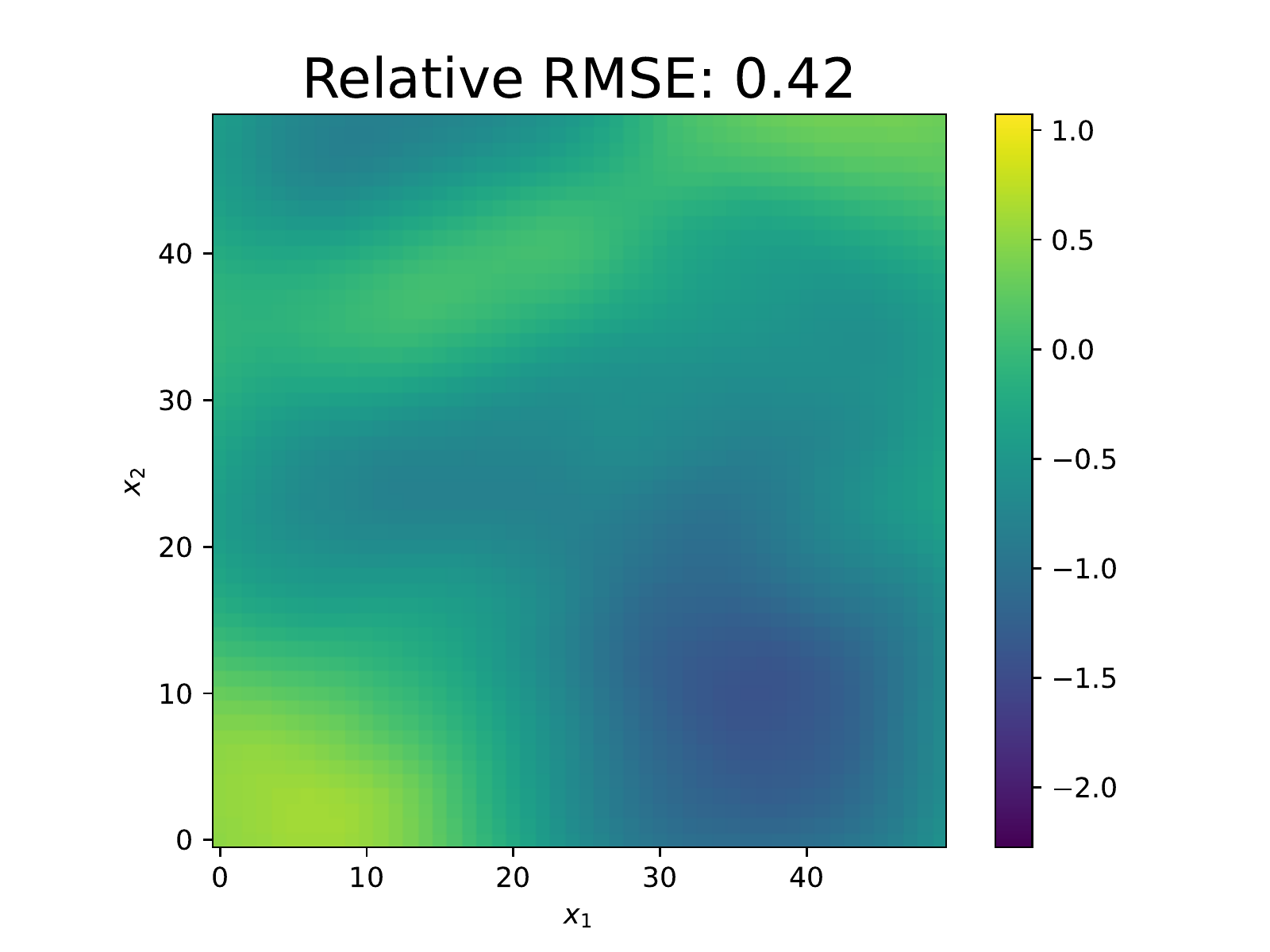}
         \caption{EKI approximated $\log(k)$.}
         \label{fig:darcy_kalman_permeability}
     \end{subfigure}
     \begin{subfigure}[t]{0.24\textwidth}
         \centering
         \includegraphics[width=\textwidth]{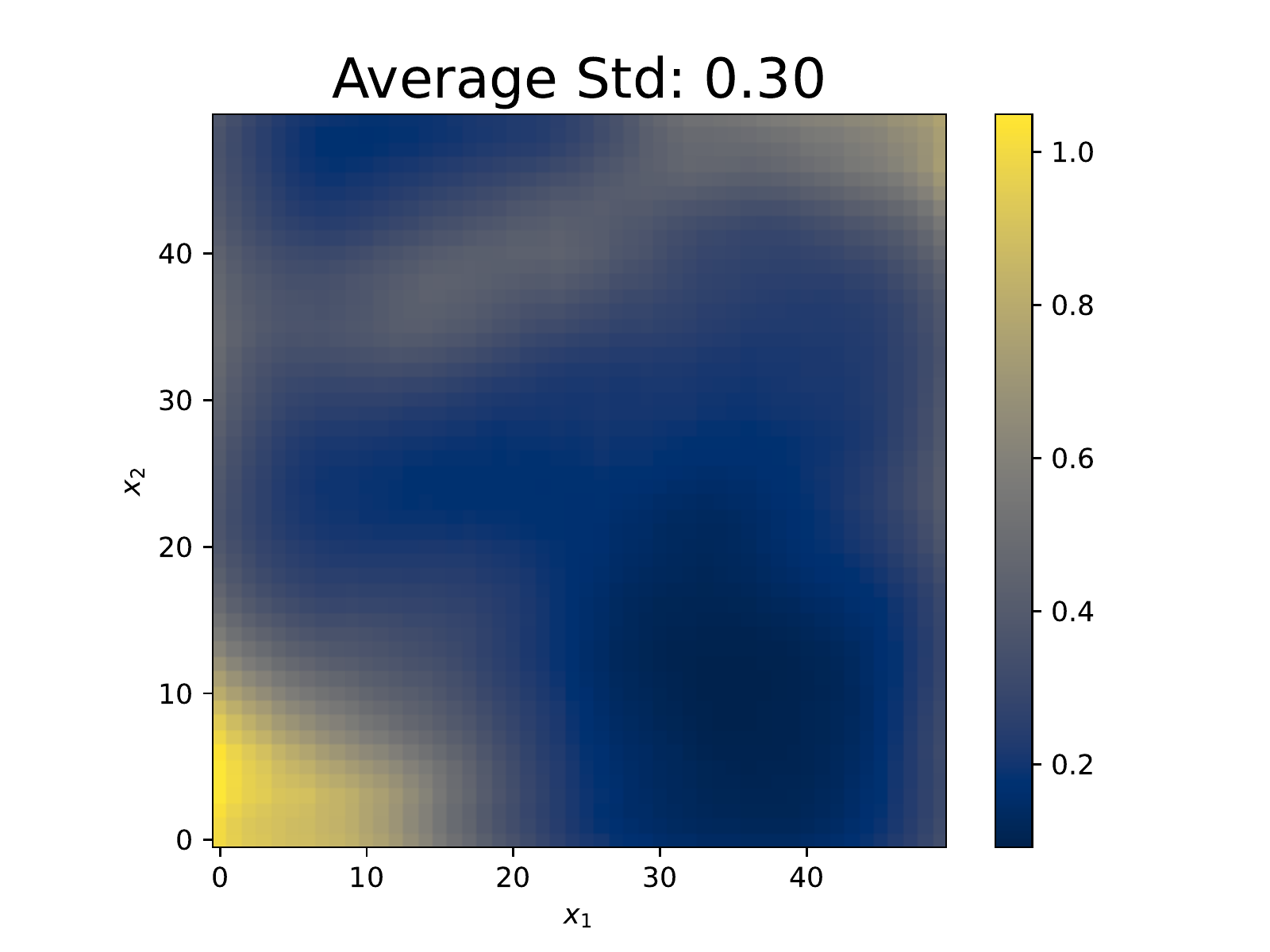}
         \caption{EKI standard deviation $\log(k)$.}
         \label{fig:darcy_kalman_permeability_std}
     \end{subfigure}
        \caption{In (a)-(b) we see the reconstruction of $v_1$ and the standard deviation of the reconstruction, respectively. In (c)-(d) we see the reconstruction of $\log(k)$ and the standard deviation of the reconstruction, respectively.}
        \label{fig:darcy_results_kalman}
\end{figure}

For the Darcy flow we compared the MCGAN results with EKI, since it is infeasible to compute  high-dimensional distributions with PCE. We compute ensembles consisting of 4000 forward computations and use 25 iterations. Note that the EKI method is parallel since each member of the ensemble can be computed independently from the other members. Therefore, we run the EKI using 20 CPU cores. For computing the permeability field, we use $n=1089$, which is the total number of degrees of freedom. See Figure \ref{fig:darcy_results_kalman} for the results.

\subsection{Leakage Detection in Pipe Flow} \label{appendix:alternative_pipe_results}
\begin{figure}[h]
     \centering
     \begin{subfigure}[t]{0.24\textwidth}
         \centering
         \includegraphics[width=\textwidth]{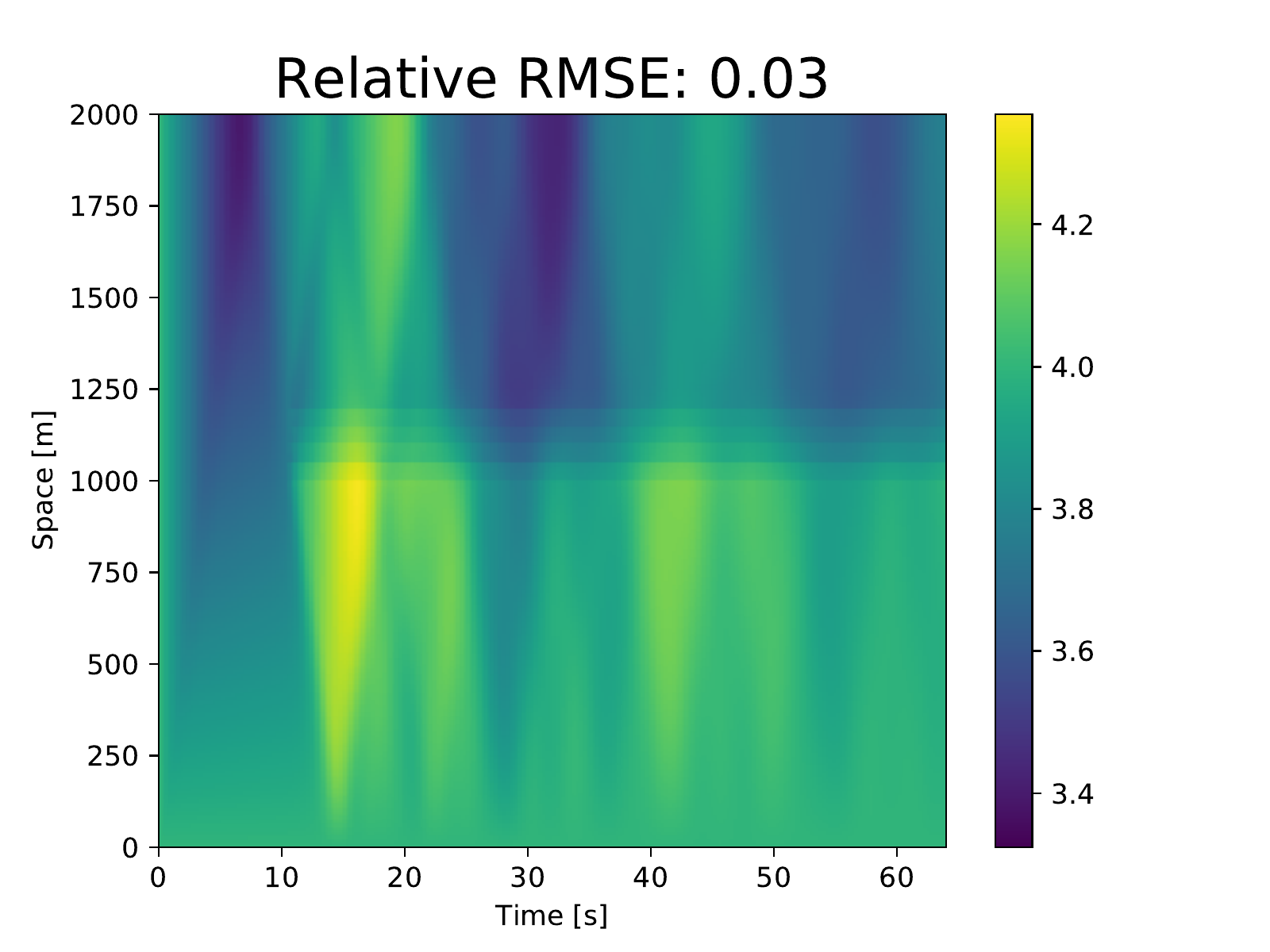}
         \caption{PCE approximated velocity.}
         \label{fig:pipe_PCE_state}
     \end{subfigure}
     \begin{subfigure}[t]{0.24\textwidth}
         \centering
         \includegraphics[width=\textwidth]{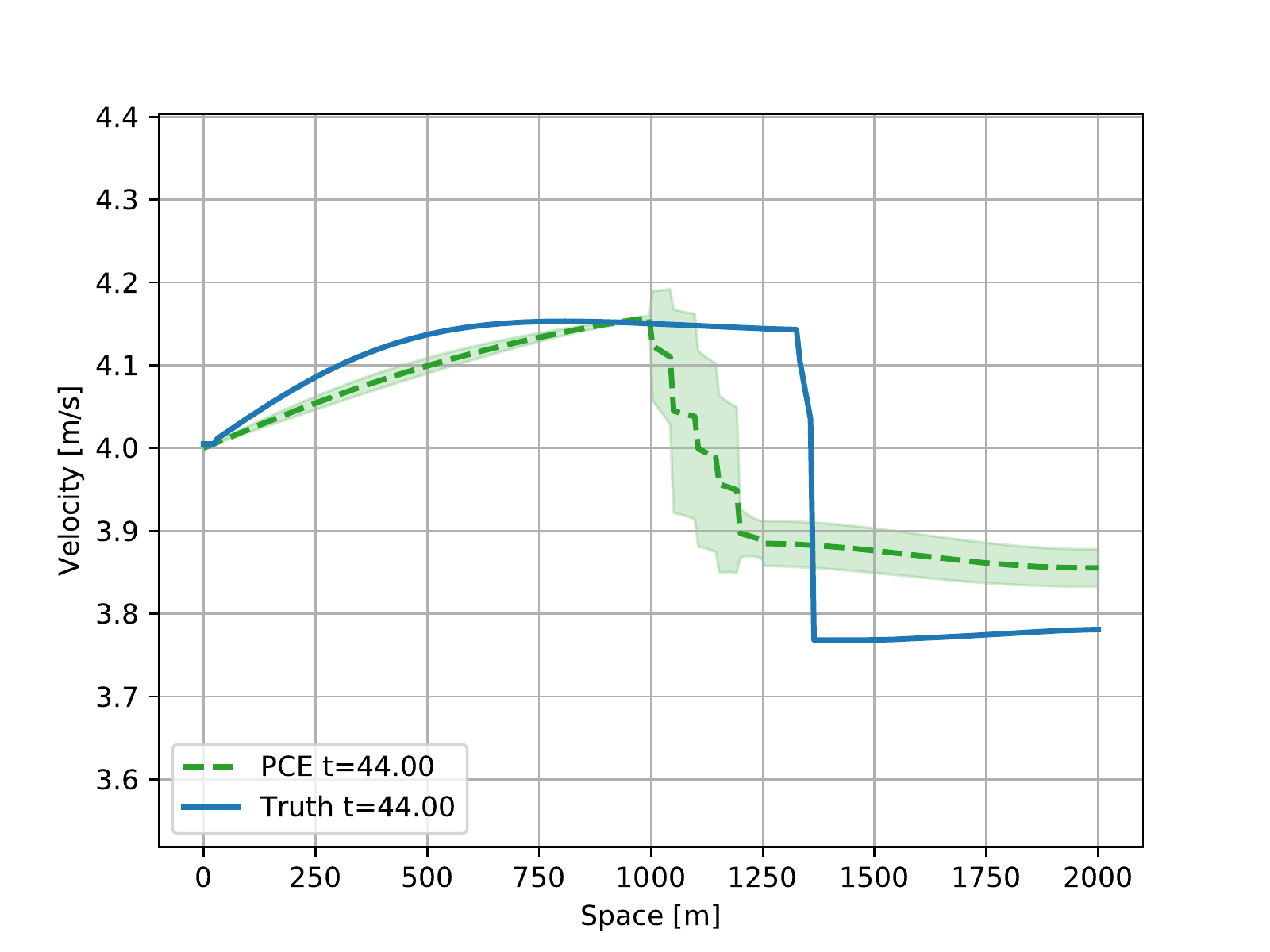}
         \caption{Reconstruction at $t=44$}
         \label{fig:pipe_reconstruction_at_t_05_PCE}
     \end{subfigure}
     \begin{subfigure}[t]{0.24\textwidth}
         \centering
         \includegraphics[width=\textwidth]{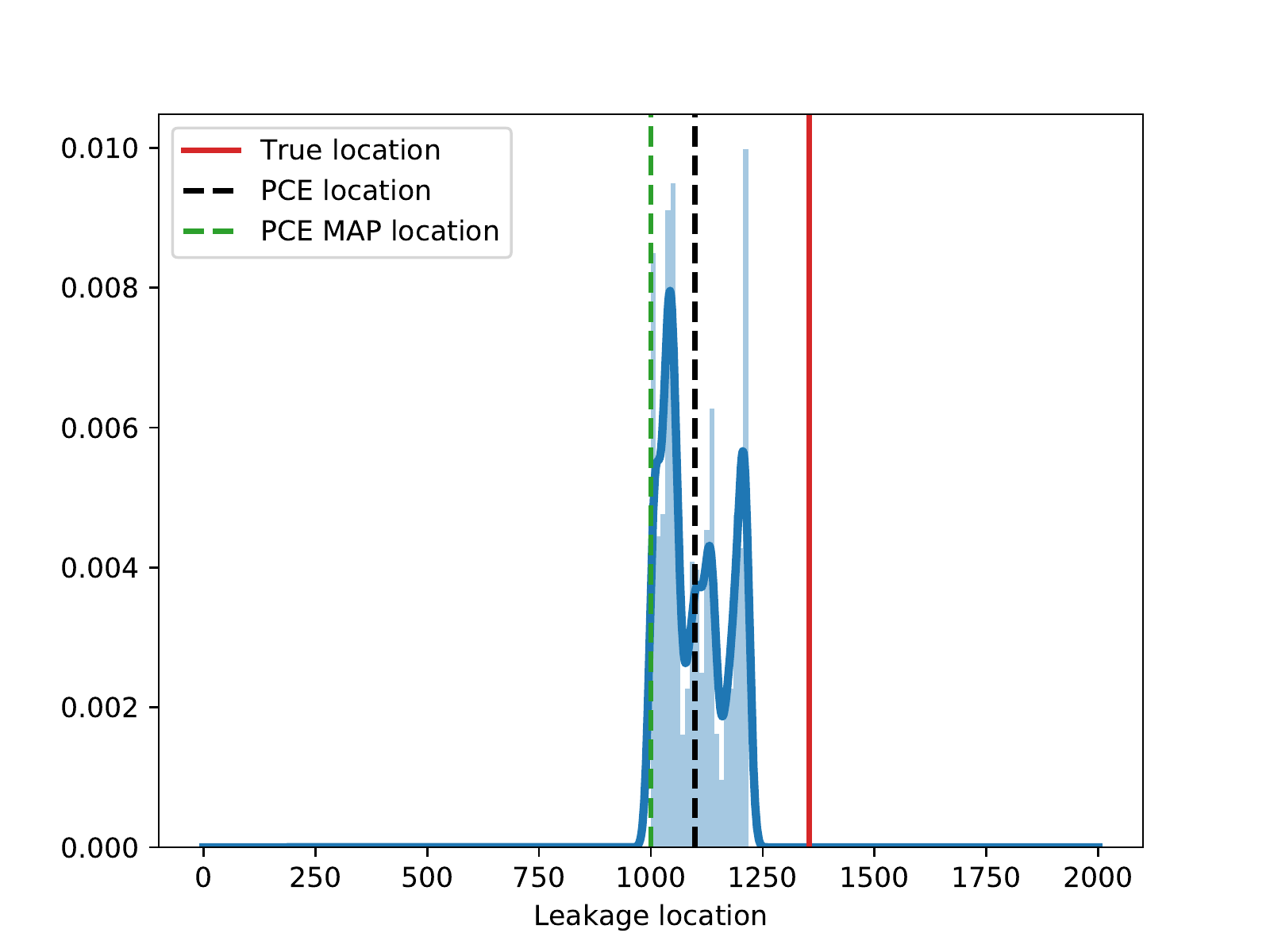}
         \caption{Leakage location.}
         \label{fig:location_histogram_PCE}
     \end{subfigure}
     \begin{subfigure}[t]{0.24\textwidth}
         \centering
         \includegraphics[width=\textwidth]{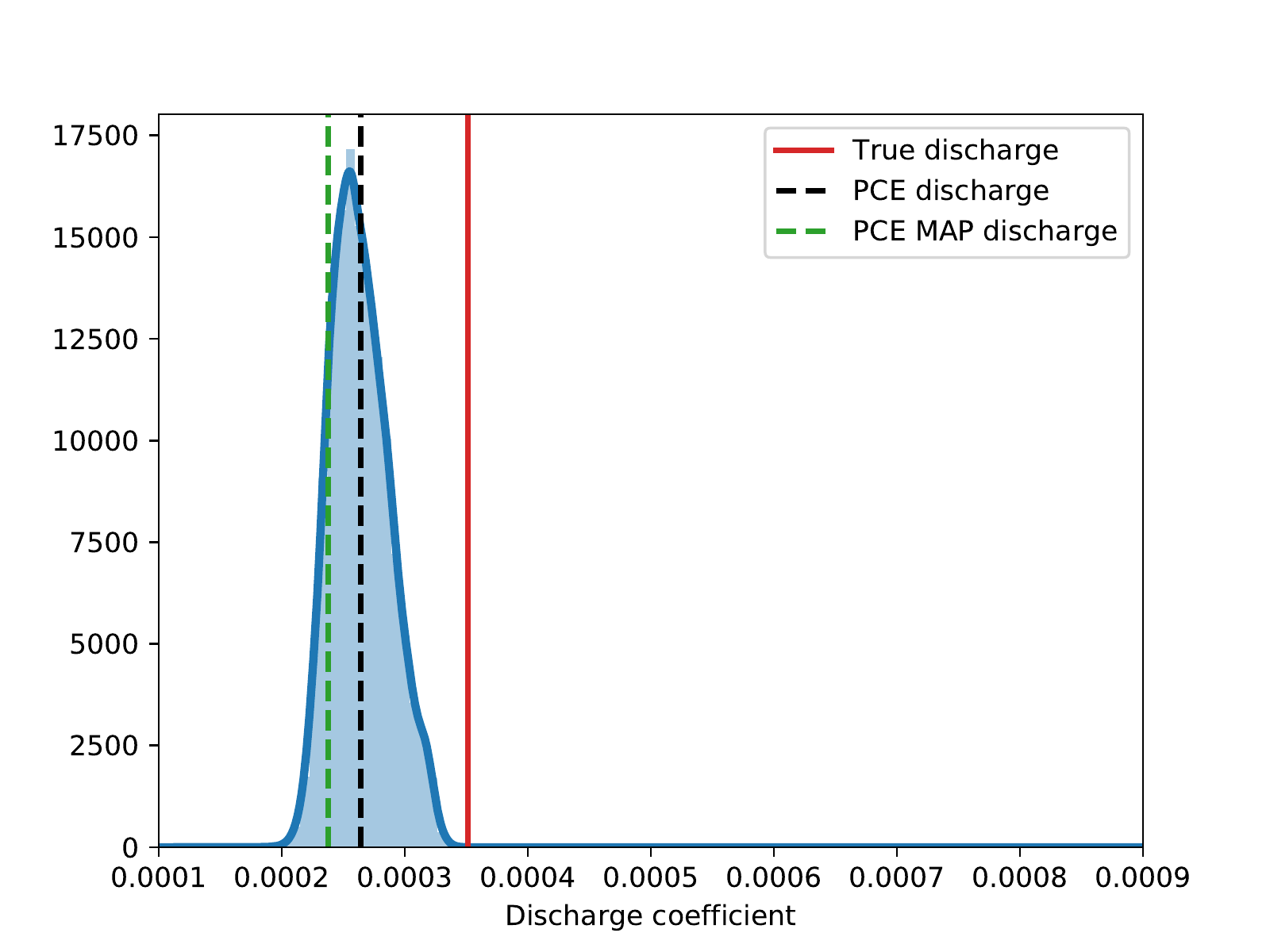}
         \caption{Discharge coefficient.}
         \label{fig:discharge_histogram_PCE}
     \end{subfigure}
    \caption{Results for the MCMC sampling with a PCE surrogate model applied to the pipe flow with a leakage, Eq. \eqref{pipe_equations}. In (a) we see the space-time contour plots of the reconstructed velocity. In (b) we see the velocity reconstruction at $t=44$ with the shaded area denoting the standard deviation. In (c)-(d) we see the posterior distributions of the leakage location and discharge coefficient.}
    \label{fig:pipe_results_PCE}
\end{figure}

For the leakage detection in the pipe, we compare our method with the PCE and the EnKF approaches. The PCE model is trained to map the leakage location, $x_l$, and discharge coefficient, $C_d$, to the observations. We achieved the highest precision with fourth-order polynomials. We performed 50,000 MCMC posterior samples and discarded the first 40,000. The state reconstruction is performed after the sampling by computing the state using the parameters samples from the MCMC sampling. The state reconstructions are computed in parallel using 20 cores. See Figure \ref{fig:pipe_results_PCE} for results.

\begin{figure}[h]
     \centering
     \begin{subfigure}[t]{0.24\textwidth}
         \centering
         \includegraphics[width=\textwidth]{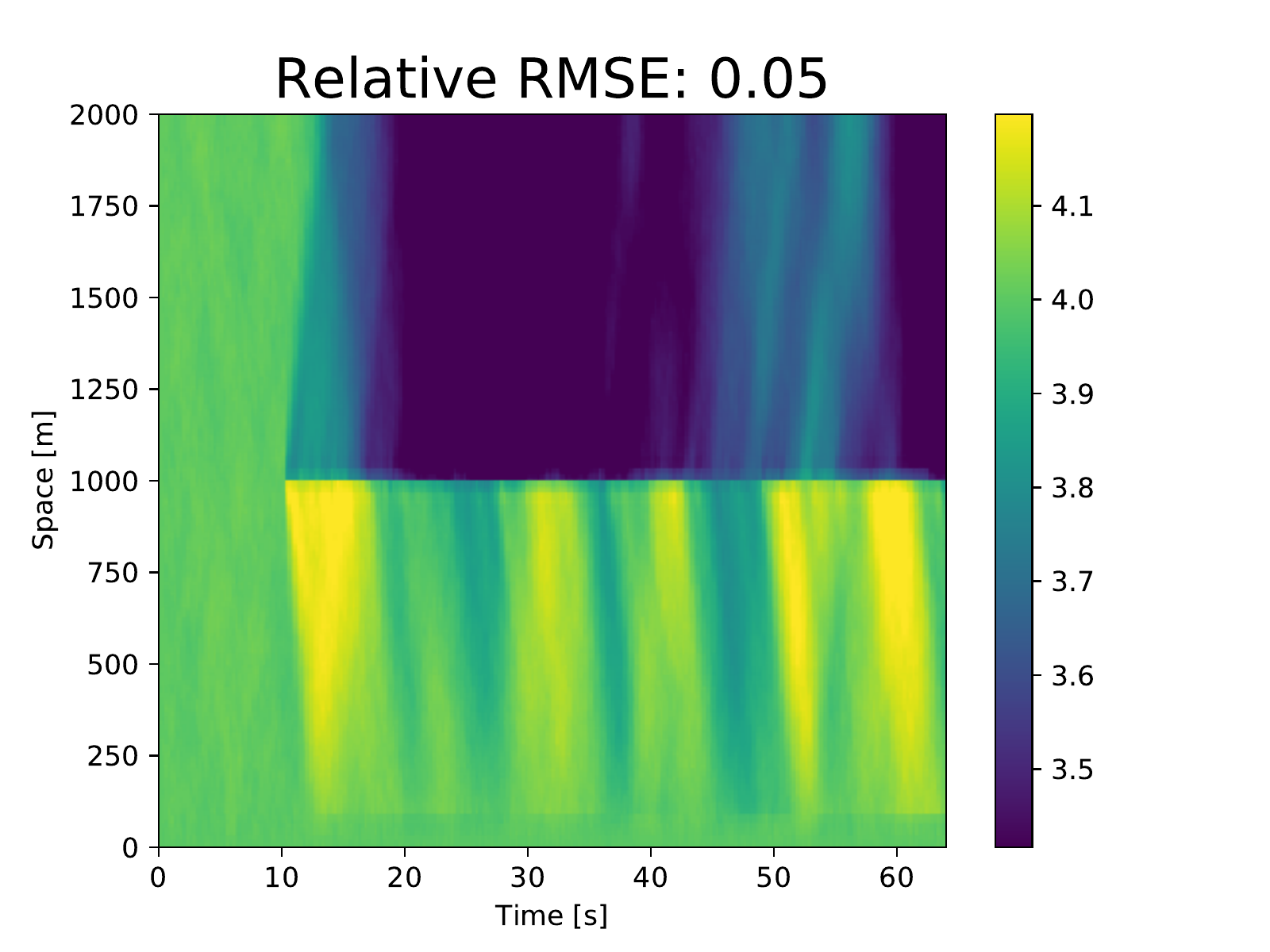}
         \caption{EnKF approximated velocity.}
         \label{pipe_kalman_state}
     \end{subfigure}
     \begin{subfigure}[t]{0.24\textwidth}
         \centering
         \includegraphics[width=\textwidth]{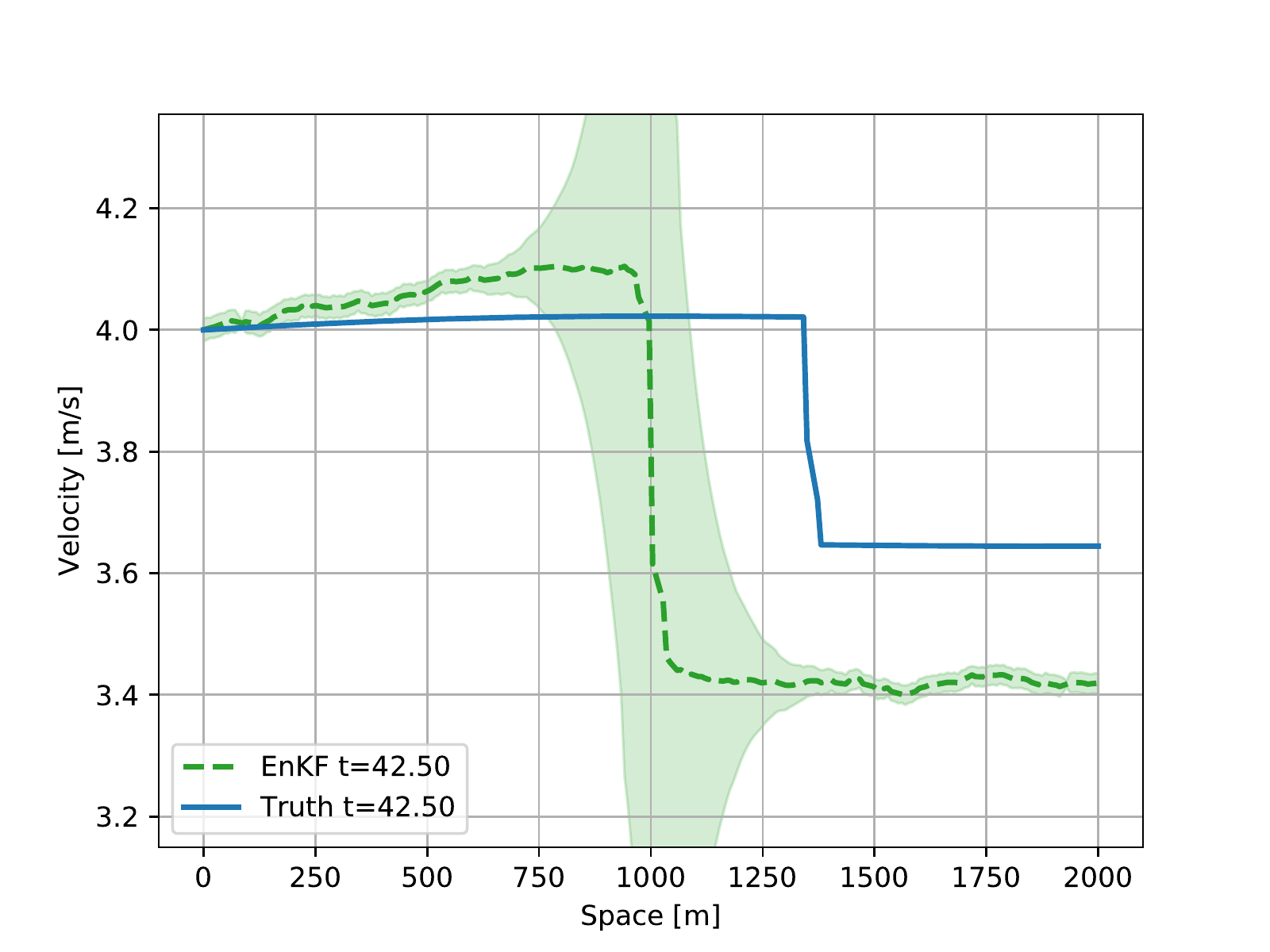}
         \caption{Reconstruction at $t=44$}
         \label{fig:pipe_reconstruction_at_t_05_kalman}
     \end{subfigure}
     \begin{subfigure}[t]{0.24\textwidth}
         \centering
         \includegraphics[width=\textwidth]{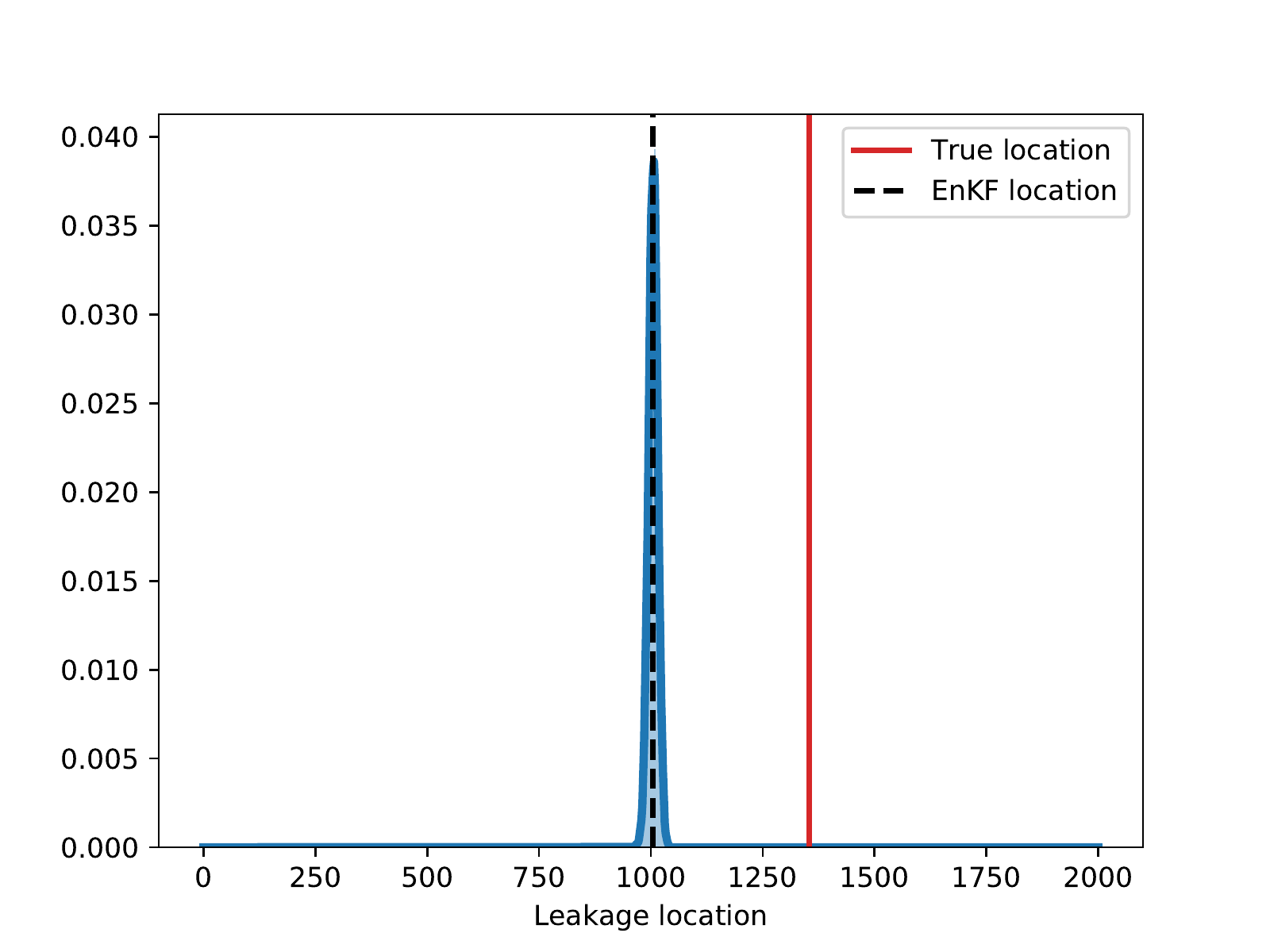}
         \caption{Leakage location.}
         \label{fig:location_histogram_kalman}
     \end{subfigure}
     \begin{subfigure}[t]{0.24\textwidth}
         \centering
         \includegraphics[width=\textwidth]{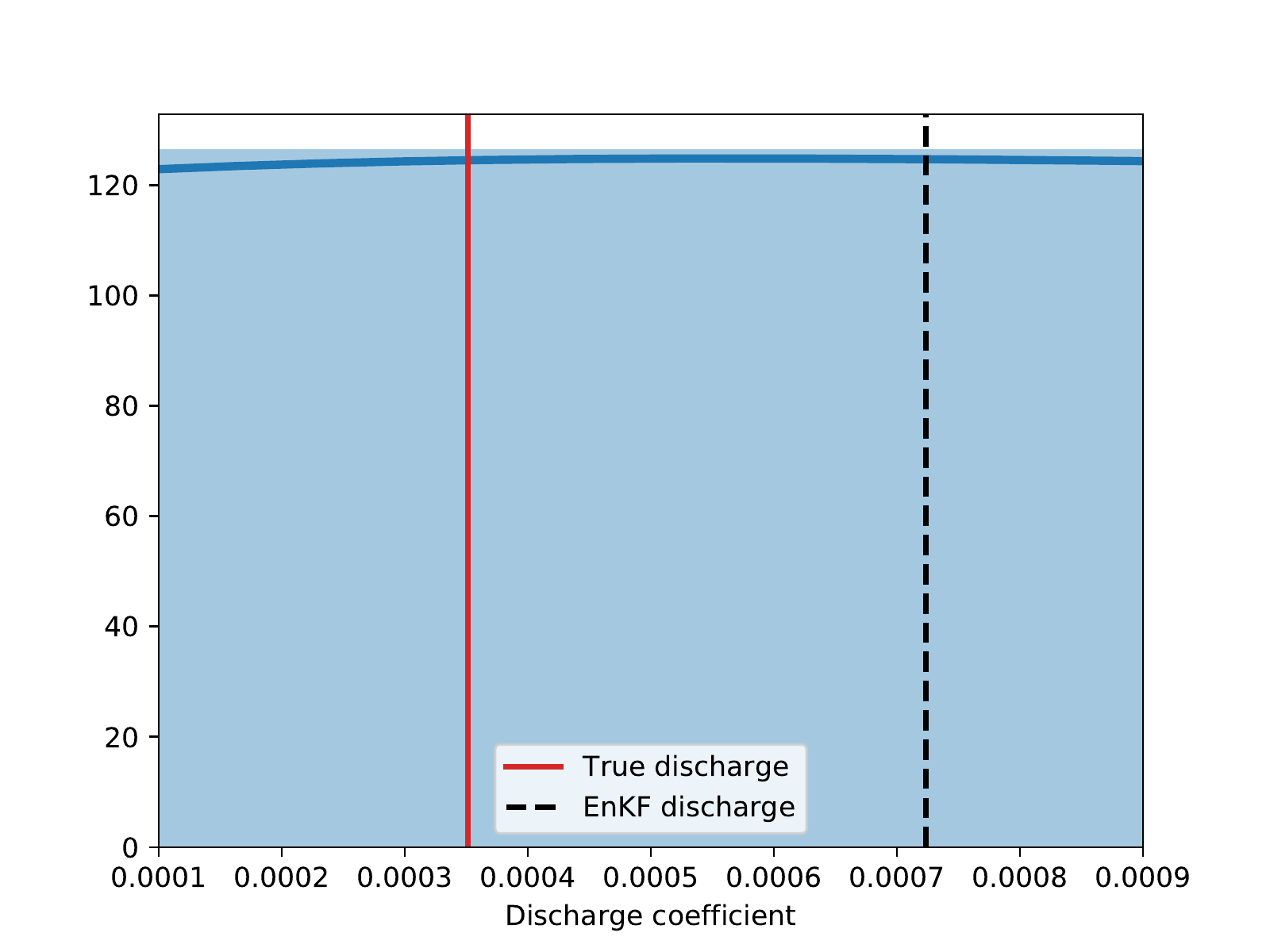}
         \caption{Discharge coefficient.}
         \label{fig:discharge_histogram_kalman}
     \end{subfigure}
    \caption{Results for EnKF method applied to the pipe flow with a leakage, Eq. \eqref{pipe_equations}. In (a) we see the space-time contour plots of the reconstructed velocity. In (b) we see the velocity reconstruction at $t=44$ with the shaded area denoting the standard deviation. In (c)-(d) we see the posterior distributions of the leakage location and discharge coefficient.}
    \label{fig:pipe_results_kalman}
\end{figure}

In the EnKF method we used an ensemble size of 2000. $\Gamma_q$ and $\Gamma_m$ are chosen to be identity matrices and $Q_q= \mathrm{diag}(0.01,0.001)^2$ and $Q_m=\mathrm{diag}(100,1\cdot 10^5)^2$. The ensemble is computed in parallel on 20 CPU cores. See Figure \ref{fig:pipe_results_kalman} for results.

%% file: main.bbl
\begin{thebibliography}{10}
\expandafter\ifx\csname url\endcsname\relax
  \def\url#1{\texttt{#1}}\fi
\expandafter\ifx\csname urlprefix\endcsname\relax\def\urlprefix{URL }\fi
\expandafter\ifx\csname href\endcsname\relax
  \def\href#1#2{#2} \def\path#1{#1}\fi

\bibitem{asch2016data}
M.~Asch, M.~Bocquet, M.~Nodet, Data assimilation: methods, algorithms, and
  applications, SIAM, 2016.

\bibitem{harlim2018data}
J.~Harlim, Data-driven computational methods: parameter and operator
  estimations, Cambridge University Press, 2018.

\bibitem{stuart2010inverse}
A.~M. Stuart, Inverse problems: a bayesian perspective, Acta numerica 19 (2010)
  451--559.

\bibitem{kaipio2006statistical}
J.~Kaipio, E.~Somersalo, Statistical and computational inverse problems, Vol.
  160, Springer Science \& Business Media, 2006.

\bibitem{kapteyn2021probabilistic}
M.~G. Kapteyn, J.~V. Pretorius, K.~E. Willcox, A probabilistic graphical model
  foundation for enabling predictive digital twins at scale, Nature
  Computational Science 1~(5) (2021) 337--347.

\bibitem{brooks2011handbook}
S.~Brooks, A.~Gelman, G.~Jones, X.-L. Meng, Handbook of markov chain monte
  carlo, CRC press, 2011.

\bibitem{gamerman2006markov}
D.~Gamerman, H.~F. Lopes, Markov chain Monte Carlo: stochastic simulation for
  Bayesian inference, CRC Press, 2006.

\bibitem{quarteroni2015reduced}
A.~Quarteroni, A.~Manzoni, F.~Negri, Reduced basis methods for partial
  differential equations: an introduction, Vol.~92, Springer, 2015.

\bibitem{xiu2010numerical}
D.~Xiu, Numerical methods for stochastic computations, Princeton university
  press, 2010.

\bibitem{wang2018adaptive}
H.~Wang, J.~Li, Adaptive gaussian process approximation for bayesian inference
  with expensive likelihood functions, Neural computation 30~(11) (2018)
  3072--3094.

\bibitem{brunton2020machine}
S.~L. Brunton, B.~R. Noack, P.~Koumoutsakos, Machine learning for fluid
  mechanics, Annual Review of Fluid Mechanics 52 (2020) 477--508.

\bibitem{baker2019workshop}
N.~Baker, F.~Alexander, T.~Bremer, A.~Hagberg, Y.~Kevrekidis, H.~Najm,
  M.~Parashar, A.~Patra, J.~Sethian, S.~Wild, et~al., Workshop report on basic
  research needs for scientific machine learning: Core technologies for
  artificial intelligence, Tech. rep., USDOE Office of Science (SC),
  Washington, DC (United States) (2019).

\bibitem{gribonval2021approximation}
R.~Gribonval, G.~Kutyniok, M.~Nielsen, F.~Voigtlaender, Approximation spaces of
  deep neural networks, Constructive Approximation (2021) 1--109.

\bibitem{mucke2021reduced}
N.~T. M{\"u}cke, S.~M. Boht{\'e}, C.~W. Oosterlee, Reduced order modeling for
  parameterized time-dependent pdes using spatially and memory aware deep
  learning, Journal of Computational Science (2021) 101408.

\bibitem{hesthaven2018non}
J.~S. Hesthaven, S.~Ubbiali, Non-intrusive reduced order modeling of nonlinear
  problems using neural networks, Journal of Computational Physics 363 (2018)
  55--78.

\bibitem{li2020fourier}
Z.~Li, N.~Kovachki, K.~Azizzadenesheli, B.~Liu, K.~Bhattacharya, A.~Stuart,
  A.~Anandkumar, Fourier neural operator for parametric partial differential
  equations, arXiv preprint arXiv:2010.08895 (2020).

\bibitem{kadeethum2021framework}
T.~Kadeethum, D.~O'Malley, J.~N. Fuhg, Y.~Choi, J.~Lee, H.~S. Viswanathan,
  N.~Bouklas, A framework for data-driven solution and parameter estimation of
  pdes using conditional generative adversarial networks, arXiv preprint
  arXiv:2105.13136 (2021).

\bibitem{ruthotto2021introduction}
L.~Ruthotto, E.~Haber, An introduction to deep generative modeling,
  GAMM-Mitteilungen (2021) e202100008.

\bibitem{goodfellow2014generative}
I.~J. Goodfellow, J.~Pouget-Abadie, M.~Mirza, B.~Xu, D.~Warde-Farley, S.~Ozair,
  A.~Courville, Y.~Bengio, Generative adversarial networks, arXiv preprint
  arXiv:1406.2661 (2014).

\bibitem{kingma2013auto}
D.~P. Kingma, M.~Welling, Auto-encoding variational bayes, arXiv preprint
  arXiv:1312.6114 (2013).

\bibitem{dhariwal2021diffusion}
P.~Dhariwal, A.~Nichol, Diffusion models beat gans on image synthesis, arXiv
  preprint arXiv:2105.05233 (2021).

\bibitem{rezende2015variational}
D.~Rezende, S.~Mohamed, Variational inference with normalizing flows, in:
  International conference on machine learning, PMLR, 2015, pp. 1530--1538.

\bibitem{goh2019solving}
H.~Goh, S.~Sheriffdeen, J.~Wittmer, T.~Bui-Thanh, Solving bayesian inverse
  problems via variational autoencoders, arXiv preprint arXiv:1912.04212
  (2019).

\bibitem{whang2021composing}
J.~Whang, E.~Lindgren, A.~Dimakis, Composing normalizing flows for inverse
  problems, in: International Conference on Machine Learning, PMLR, 2021, pp.
  11158--11169.

\bibitem{patel2020bayesian}
D.~V. Patel, D.~Ray, H.~Ramaswamy, A.~Oberai, Bayesian inference in
  physics-driven problems with adversarial priors, in: NeurIPS 2020 Workshop on
  Deep Learning and Inverse Problems, 2020.

\bibitem{xia2022bayesian}
Y.~Xia, N.~Zabaras, Bayesian multiscale deep generative model for the solution
  of high-dimensional inverse problems, Journal of Computational Physics 455
  (2022) 111008.

\bibitem{metropolis1953equation}
N.~Metropolis, A.~W. Rosenbluth, M.~N. Rosenbluth, A.~H. Teller, E.~Teller,
  Equation of state calculations by fast computing machines, The journal of
  chemical physics 21~(6) (1953) 1087--1092.

\bibitem{hastings1970monte}
W.~K. Hastings, Monte carlo sampling methods using markov chains and their
  applications (1970).

\bibitem{hoffman2014no}
M.~D. Hoffman, A.~Gelman, The no-u-turn sampler: adaptively setting path
  lengths in hamiltonian monte carlo., J. Mach. Learn. Res. 15~(1) (2014)
  1593--1623.

\bibitem{stuart2018posterior}
A.~Stuart, A.~Teckentrup, Posterior consistency for gaussian process
  approximations of bayesian posterior distributions, Mathematics of
  Computation 87~(310) (2018) 721--753.

\bibitem{jabbar2020survey}
A.~Jabbar, X.~Li, B.~Omar, A survey on generative adversarial networks:
  Variants, applications, and training, arXiv preprint arXiv:2006.05132 (2020).

\bibitem{arjovsky2017wasserstein}
M.~Arjovsky, S.~Chintala, L.~Bottou, Wasserstein generative adversarial
  networks, in: International conference on machine learning, PMLR, 2017, pp.
  214--223.

\bibitem{bremaud2020probability}
P.~Br{\'e}maud, Probability Theory and Stochastic Processes, Springer Nature,
  2020.

\bibitem{gulrajani2017improved}
I.~Gulrajani, F.~Ahmed, M.~Arjovsky, V.~Dumoulin, A.~Courville, Improved
  training of wasserstein gans, arXiv preprint arXiv:1704.00028 (2017).

\bibitem{liu2017approximation}
S.~Liu, O.~Bousquet, K.~Chaudhuri, Approximation and convergence properties of
  generative adversarial learning, arXiv preprint arXiv:1705.08991 (2017).

\bibitem{sanderse2021efficient}
B.~Sanderse, V.~V. Dighe, K.~Boorsma, G.~Schepers, Efficient bayesian
  calibration of aerodynamic wind turbine models using surrogate modeling, Wind
  Energy Science Discussions (2021) 1--34.

\bibitem{lu2015limitations}
F.~Lu, M.~Morzfeld, X.~Tu, A.~J. Chorin, Limitations of polynomial chaos
  expansions in the bayesian solution of inverse problems, Journal of
  Computational Physics 282 (2015) 138--147.

\bibitem{bogachev2007measure}
V.~I. Bogachev, Measure theory, Vol.~1, Springer Science \& Business Media,
  2007.

\bibitem{sprungk2020local}
B.~Sprungk, On the local lipschitz stability of bayesian inverse problems,
  Inverse Problems 36~(5) (2020) 055015.

\bibitem{panaretos2020invitation}
V.~M. Panaretos, Y.~Zemel, An invitation to statistics in Wasserstein space,
  Springer Nature, 2020.

\bibitem{colton1998inverse}
D.~L. Colton, R.~Kress, R.~Kress, Inverse acoustic and electromagnetic
  scattering theory, Vol.~93, Springer, 1998.

\bibitem{domesova2017solution}
S.~Domesov{\'a}, M.~Beres, Solution of inverse problems using bayesian approach
  with application to estimation of material parameters in darcy flow.,
  Advances in Electrical \& Electronic Engineering 15~(2) (2017).

\bibitem{ruchi2019transform}
S.~Ruchi, S.~Dubinkina, M.~Iglesias, Transform-based particle filtering for
  elliptic bayesian inverse problems, Inverse Problems 35~(11) (2019) 115005.

\bibitem{kumar2018multigrid}
P.~Kumar, P.~Luo, F.~J. Gaspar, C.~W. Oosterlee, A multigrid multilevel monte
  carlo method for transport in the darcy--stokes system, Journal of
  Computational Physics 371 (2018) 382--408.

\bibitem{cockburn2009superconvergent}
B.~Cockburn, J.~Guzm{\'a}n, H.~Wang, Superconvergent discontinuous galerkin
  methods for second-order elliptic problems, Mathematics of Computation
  78~(265) (2009) 1--24.

\bibitem{logg2012automated}
A.~Logg, K.-A. Mardal, G.~Wells, Automated solution of differential equations
  by the finite element method: The FEniCS book, Vol.~84, Springer Science \&
  Business Media, 2012.

\bibitem{ding2021ensemble}
Z.~Ding, Q.~Li, Ensemble kalman inversion: mean-field limit and convergence
  analysis, Statistics and Computing 31~(1) (2021) 1--21.

\bibitem{kundu2002fluid}
P.~K. Kundu, I.~M. Cohen, Fluid mechanics (2002).

\bibitem{hauge2007model}
E.~Hauge, O.~M. Aamo, J.-M. Godhavn, Model based pipeline monitoring with leak
  detection, IFAC Proceedings Volumes 40~(12) (2007) 318--323.

\bibitem{schetz1996handbook}
J.~A. Schetz, A.~E. Fuhs, Handbook of fluid dynamics and fluid machinery,
  Vol.~1, Wiley New York, 1996.

\bibitem{hesthaven2007nodal}
J.~S. Hesthaven, T.~Warburton, Nodal discontinuous Galerkin methods:
  algorithms, analysis, and applications, Springer Science \& Business Media,
  2007.

\bibitem{leveque2007finite}
R.~J. LeVeque, Finite difference methods for ordinary and partial differential
  equations: steady-state and time-dependent problems, SIAM, 2007.

\bibitem{horsholt2019spatial}
A.~H{\o}rsholt, L.~H. Christiansen, K.~Meyer, J.~K. Huusom, J.~B. J{\o}rgensen,
  Spatial discretization and kalman filtering for ideal packed-bed
  chromatography, in: 2019 18th European Control Conference (ECC), IEEE, 2019,
  pp. 2356--2361.

\bibitem{feinberg2015chaospy}
J.~Feinberg, H.~P. Langtangen, Chaospy: An open source tool for designing
  methods of uncertainty quantification, Journal of Computational Science 11
  (2015) 46--57.

\end{thebibliography}
